\documentclass{ipiof}
\usepackage{amsmath,amssymb,amsthm}
\usepackage[utf8]{inputenc}
\usepackage{enumitem} 
\usepackage{graphics}
\usepackage{epsfig}
\usepackage{graphicx}
\usepackage{epstopdf}
\usepackage[colorlinks=true]{hyperref}
\hypersetup{urlcolor=blue, citecolor=red}

\def\currentvolume{X}
 \def\currentissue{X}
  \def\currentyear{200X}
   \def\currentmonth{XX}
    \def\ppages{X--XX}
     \def\DOI{10.3934/ipi.2021042}

\theoremstyle{definition}

\newtheorem{remark}{Remark}

\title[Parameter identification for the LLG equation in MPI]{On numerical aspects of parameter identification for the Landau-Lifshitz-Gilbert equation in Magnetic Particle Imaging}

\author[Tram Thi Ngoc Nguyen and Anne Wald]{}


 \email{tram.nguyen@uni-graz.at}
 \email{a.wald@math.uni-goettingen.de}

\thanks{$^*$ Corresponding author: Tram Thi Ngoc Nguyen}

\def\R{\mathbb{R}}
\def\bF{\mathbb{F}}
\def\altil{\hat{\alpha}}
\def\mv{\mathbf{m}}
\def\mh{\hat{\mathbf{m}}}
\def\uv{\mathbf{u}}
\def\vv{\mathbf{v}}
\def\yv{\mathbf{y}}
\def\wv{\mathbf{w}}
\def\zv{\mathbf{z}}
\def\hv{\mathbf{h}}
\def\Tt{T}
\def\ker{\mathbf{K}}
\def\cU{\mathcal{U}}

  \newtheorem{theo}{Theorem}[section]

  \newtheorem{algo}[theo]{Algorithm}

\def\N{\mathbb{N}}
\def\cD{\mathcal{D}}
\def\cK{\mathcal{K}}
\def\cX{\mathcal{X}}
\def\cY{\mathcal{Y}}
\def\cB{\mathcal{B}}

\def\cW{\mathcal{W}}

\def\m{\textbf{m}}

\def\u{\textbf{u}}
\def\s{\textbf{s}}
\def\f{\textbf{f}}
\def\g{\textbf{g}}

\def\pii{\mathbf{\pi}}
\def\pp{\textbf{q}}
\def\p{\pp^z}
\def\pr{\textbf{p}_{\ell}^R}

\def\h{\textbf{h}}

\def\Ktil{\tilde{K}}
\def\KtilT{\tilde{K}_T}

\def\Om2R{\Omega,\R^3}

\newcommand{\pR}{\mathbf{p}^{\mathrm{R}}}

\newcommand{\blue}[1]{\textcolor{black}{#1}}

\begin{document}
\maketitle

\centerline{\scshape Tram Thi Ngoc Nguyen$^*$}
\medskip
{\footnotesize
 \centerline{Institute of Mathematics and Scientific Computing, University of Graz}
   \centerline{Heinrichstra{\ss}e 36, A-8010 Graz, Austria}
}

\medskip

\centerline{\scshape Anne Wald}
\medskip
{\footnotesize
 \centerline{Institute of Numerical and Applied Mathematics, University of G\"ottingen}
   \centerline{Lotzestraße 16-18, 37073 G\"ottingen, Germany}
}

\bigskip

 \centerline{(Communicated by Habib Ammari)}

\begin{abstract}
The Landau-Lifshitz-Gilbert equation yields a mathematical model to describe the evolution of the magnetization of a magnetic material, particularly in response to an external applied magnetic field. It allows one to take into account various physical effects, such as the exchange within the magnetic material itself. In particular, the Landau-Lifshitz-Gilbert equation encodes relaxation effects, i.e., it describes the time-delayed alignment of the magnetization field with an external magnetic field. These relaxation effects are an important aspect in magnetic particle imaging, particularly in the calibration process. In this article, we address the data-driven modeling of the system function in magnetic particle imaging, where the Landau-Lifshitz-Gilbert equation serves as the basic tool to include relaxation effects in the model. We formulate the respective parameter identification problem both in the all-at-once and the reduced setting, present reconstruction algorithms that yield a regularized solution and discuss numerical experiments. Apart from that, we propose a practical numerical solver to the nonlinear Landau-Lifshitz-Gilbert equation, not via the classical finite element method, but through solving only linear PDEs in an inverse problem framework.

\end{abstract}

\section{Introduction} \label{sec:introduction}

In \cite{KNSW}, a parameter identification problem for the nonlinear Landau-Lifshitz-Gilbert equation was introduced and analyzed, in particular in view of a model-based calibration in magnetic particle imaging (MPI). 
MPI is a novel medical imaging technique, see \cite{Gleich2005, tkngmm17magnetic} that aims at, for example, the imaging of blood vessels, see, e.g., \cite{Borgert2012}. To this end, magnetic nano-particles are injected into the blood stream. A strong dynamic external magnetic field with a field-free point is then applied to change the magnetization of the particles, which induces an electric voltage in the receive coils of the scanner. Particularly those particles close to the field-free point undergo an abrupt change of their magnetization, which ensures that spatial information is encoded in the measured signals. 
During data generation, the field-free point is driven through the entire field-of-view on a pre-defined trajectory, resulting in a scan of the region of interest. A detailed description of this imaging technique can be found in \cite{tkngmm17magnetic}. 

\vspace*{2ex}

The forward problem in MPI is essentially described by an integral equation
\begin{displaymath}
 v(t) = \int_{\Omega} c(\mathbf{x}) s(\mathbf{x},t) \,\mathrm{d}\mathbf{x},
\end{displaymath}
where the imaging consists in reconstructing the particle concentration $c$ from the measured data $v$. The calibration process, on the other hand, is the determination of the system function $s$ from measurements $v$ for known calibration concentrations $c$. 

In \cite{KNSW} as well as in this article, we are concerned with the determination of the system function $s$. More precisely, we aim at a physical model for $s$, for which we identify some underlying parameters. The system function consists of the time-derivative of the magnetization of the magnetic material as well as some quantities describing the measurement process. The model is based on the Landau-Lifshitz-Gilbert equation, which is a partial differential equation that describes the evolution of the magnetization in response to an external magnetic field. The Landau-Lifshitz-Gilbert equation allows us to include a range of physical effects into the mathematical model, in particular it admits the inclusion of relaxation effects, which highly influence the measured data \cite{lcpgsc12, tk18}. The use of the Landau-Lifshitz-Gilbert equation is motivated by its applications in micromagnetism \cite{kruzikprohl}, since it not only incorporates relaxation effects, but allows to include further physical effects such as the exchange within the magnetic material. 

Our main goal in this article is thus a data-driven determination of the underlying mathematical model to enable a model-based approach for the system function, which particularly includes relaxation effects. 

In parameter identification, one usually aims at the reconstruction of parameters that appear in the coefficients of a (partial) differential equation from the knowledge of the solution or the state of the system, or even boundary measurements thereof. Examples are electric impedance tomography \cite{borcea02}, photoacoustic tomography, optical coherence tomography \cite{Elbau_et_al_2017} or dynamic load monitoring \cite{Binder_2015}.
The inverse problem we are concerned with in this article does not only aims at the reconstruction of parameters, but also requires the determination of the respective solution or state function. This motivated the formulation in both the all-at-once and the reduced version, which have been presented and analyzed in \cite{Kaltenbacher} and specifically in \cite{KNSW} for the identification of the system function in MPI via the Landau-Lifshitz-Gilbert equation. While \cite{KNSW} contains a mathematical analysis of this inverse problem, we are now aiming at numerical experiments to evaluate the proposed reconstruction approaches from \cite{KNSW}.

Furthermore, we propose a numerical solver to the nonlinear Landau-Lifshitz-Gilbert (LLG) equation, not via the classical approaches of finite difference or finite element discretizations, but through formulating it as an inverse problem and solving it via a new all-at-once approach. This is an important problem with a wide range of applications, and attracts much attention from physicists and mathematicians \cite{kruzikprohl,cimrak}. Several classical strategies exist;  some well-established algorithms include \blue{\cite{BartelsProhl, BartelsProhl1} by Bartels and Prohl in 2006, 2008}, \cite{alouges, alougesKritsikisSteinerToussaint} by Alouges et.~al.~in 2008, 2014, and \cite{BanasEtAl} by Ba\v nas et.~al.~in 2014. These authors investigated several forms of the LLG equation, e.g.,  $p$-harmonic heat flow, equations with magnetostriction effect,  Maxwell-LLG equations etc. Our contribution consists of a practical numerical solver for the highly nonlinear LLG, based on solving only linear PDEs in an inverse problem framework. A detailed discussion is presented in Section \ref{sec-algorithm-AAO}.

The outline is as follows. In Section \ref{Sec:theory}, we summarize the physical model as well as the theoretical foundations from \cite{KNSW} that are needed for a numerical solution of the addressed parameter identification problem. The algorithms that are used for the implementation are presented and discussed in Section \ref{topci4-sec:algorithm}, the respective numerical experiments are to be found in Section \ref{topci4-sec:example}. We conclude this article with a short discussion on the use of the proposed methods to calibrate an MPI scanner.

\section{Parameter identification for the Landau-Lifshitz-Gilbert equation} \label{Sec:theory}
In this section, we want to introduce the modelling aspects and recap some of the results from \cite{KNSW}, which yield the theoretical basis for the numerical results presented in this work. 

\subsection{The forward problem}
Let $\Omega \subseteq \mathbb{R}^n$, $n=1,2,3$, denote the field-of-view. The particle magnetization $\mathbf{M}^{\mathrm{P}}(x,t)$ at a point $x \in \Omega$ and time point $t \in I:=[0,T]$, $T>0$, is given by 
\begin{displaymath}
 \mathbf{M}^{\mathrm{P}}(x,t) = c(x,t)\mathbf{m}(x,t), 
\end{displaymath}
where $c$ is a dimensionless quantity with values in $[0,1]$ corresponding to the concentration of magnetic material in $(x,t)$ and $\mathbf{m}$ is the respective magnetization of the magnetic material in $(x,t)$. In contrast to $\mathbf{M}^{\mathrm{P}}$, the magnetization has a fixed length $m_{\mathrm{s}} := \lvert \mathbf{m} \rvert$ that only depends on the magnetic material. 

\vspace*{2ex}

In the presence of a dynamic external field $\mathbf{H}_{\mathrm{ext}}$, the article magnetization changes in response to temporal changes of $\mathbf{H}_{\mathrm{ext}}$. According to Faraday's law of induction, the changes in the external field as well as in the particle magnetization can be measured by a receive coil, where they induce an electric voltage. The signal generated in the $l$-th coil, $l=1,...,L$, by the external field is filtered out from the measured signal by means of the transfer function $a_l: [0,T] \to \mathbb{R}$ resp.~its periodic continuation $\widetilde{a}_l: \mathbb{R} \to \mathbb{R}$. Together with the sensitivity $\mathbf{p}_l^{\mathrm{R}}$, which is a geometrical property of the receive coil, we obtain 
\begin{equation}
 v_l(t) = \int_0^T \int_{\Omega} \mathbf{K}_{l}(t,\tau,x) \cdot \frac{\partial}{\partial \tau} \mathbf{m}(x,\tau) \, \mathrm{d} x \, \mathrm{d}\tau
\end{equation}
with
\begin{displaymath}
 \mathbf{K}_{l}(t,\tau,x) := -\mu_0 \widetilde{a}_l(t-\tau) c(x,\tau) \, \mathbf{p}_l^{\mathrm{R}}(x)
\end{displaymath}
The physical constant $\mu_0$ is called the magnetic permeability.
We assume that $\mathbf{p}^{\mathrm{R}}$ is known. 
The temporal derivative of the magnetization $\mathbf{m}$ is described by the Landau-Lifshitz-Gilbert equation, which is a microscopic model (see, e.g., \cite{kruzikprohl}), as follows:
\begin{align}
 \mathbf{m}_t &= - \alpha_1 \mathbf{m} \times \left( \mathbf{m} \times (\Delta \mathbf{m} + \mathbf{h}) \right) + \alpha_2 \mathbf{m} \times (\Delta \mathbf{m} + \mathbf{h}) \ && \text{in} \ \Omega \times [0,T], \label{llg1}\\
 0 &= \partial_{\nu} \mathbf{m}  && \text{on} \ \partial\Omega \times [0,T], \label{llg1_bc}\\
 \mathbf{m}_0 &= \mathbf{m}(t=0), \ \lvert \mathbf{m}_0 \rvert = m_{\mathrm{S}} && \text{in} \ \Omega, \label{llg1_abs}
\end{align}
with
\begin{displaymath}
\alpha_1 := \frac{2A\gamma \alpha_{\mathrm{D}}}{m_{\mathrm{S}}(1+\alpha_{\mathrm{D}}^2)} > 0, \quad \alpha_2 := \frac{2A\gamma}{(1+\alpha_{\mathrm{D}}^2)} > 0
\end{displaymath}
and the scaled external field
\begin{displaymath}
 \mathbf{h} = \frac{\mu_0 m_{\mathrm{S}}}{2A} \blue{\mathbf{H}_{\mathrm{ext}}}.
\end{displaymath}
The saturation magnetization is denoted by $m_{\mathrm{S}}$. By $\gamma$ we denote the gyromagnetic constant, the parameter $\alpha_D$ is a damping parameter and $A$ is the exchange stiffness constant. The latter two constants as well as $m_{\mathrm{S}}$ are material-dependent.  
The term including the Laplacian in \eqref{llg1} describes the interaction within the magnetic material.

For the reconstruction, various known particle concentrations $c_k$, $k=1,...,K$, are used to increase the amount of data points $v_{kl}(t)$. In this case, we have
\begin{equation}
 v_{kl}(t) = \int_0^T \int_{\Omega} \mathbf{K}_{kl}(t,\tau,x) \cdot \frac{\partial}{\partial \tau} \mathbf{m}(x,\tau) \, \mathrm{d} x \, \mathrm{d}\tau,
\end{equation}
where
\begin{displaymath}
 \mathbf{K}_{kl}(t,\tau,x) := -\mu_0 \widetilde{a}_l(t-\tau) c_k(x,\tau) \, \mathbf{p}_l^{\mathrm{R}}(x)
\end{displaymath}

Since the vector field $\mathbf{m}$ has a fixed length $m_{\mathrm{S}}$ in $\Omega$, the Landau-Lifshitz-Gilbert equation can be reformulated (see \cite{KNSW}, Section 2.2) as 
\begin{alignat}{3}
\hat{\alpha}_1 m_{\mathrm{S}}^2 \mathbf{m}_t - \hat{\alpha}_2 \mathbf{m} \times \mathbf{m}_t -m_{\mathrm{S}}^2\Delta \mathbf{m} &= \lvert \nabla \mathbf{m} \rvert^2 \mathbf{m} + m_{\mathrm{S}}^2\mathbf{h} - 
\langle \mathbf{m}, \mathbf{h} \rangle \mathbf{m}  \,\,&& \text{in } [0,T] \times \Omega \label{llg3}\\
 0 &= \partial_{\nu} \mathbf{m}  && \text{on } 
  [0,T] \times \partial\Omega \label{llg3_bc}\\
 \mathbf{m}_0 &= \mathbf{m}(t=0), \ \lvert \mathbf{m}_0 \rvert = m_{\mathrm{S}} && \text{in } \Omega \label{llg3_abs}\,.
\end{alignat}
with
\begin{displaymath}
 \hat{\alpha}_1 = \frac{\alpha_1}{m_{\mathrm{S}}^2\alpha_1^2+\alpha_2^2}, \quad\hat{\alpha}_2=\frac{\alpha_2}{m_{\mathrm{S}}^2\alpha_1^2+\alpha_2^2}.
\end{displaymath}

For the all-at-once setting, we use $\widetilde{a} \in L^2(0,T)$, $c \cdot \mathbf{p}^{\mathrm{R}} \in L^2(\Omega,\mathbb{R}^3)$, $\mathbf{m}_0 \in H^1(\Omega,\mathbb{R}^3)$ and $\mathbf{h} \in L^2(0,T;L^p(\Omega,\mathbb{R}^3))$, $p\geq 2$. For the reduced setting, we assume $\mathbf{h} \in L^2(0,T;H^1(\Omega,\mathbb{R}^3))$ with $\partial_{\nu}\mathbf{h}=0$ on $\partial\Omega$ or $\mathbf{h} \in L^2(0,T;H^2(\Omega,\mathbb{R}^3))$, $c \cdot \mathbf{p}^{\mathrm{R}} \in H^1(\Omega,\mathbb{R}^3)$, and $\mathbf{m}_0 \in H^2(\Omega,\mathbb{R}^3)$ (see also \cite{KNSW} for more details).  


\subsection{The inverse problem}
The inverse problem of identifying the time-derivative of the magnetization $\mathbf{m}$ along with the two constants $\hat{\alpha}_1$ and $\hat{\alpha}_2$ has been formulated in \cite{KNSW}, Sections 3, 4.1 and 4.2, in both the all-at-once as well as in the reduced setting, see also \cite{Kaltenbacher}.

\subsubsection{The reduced formulation}
Since $\mathbf{m}$ satisfies \eqref{llg3} - \eqref{llg3_abs}, which has a unique solution (see \cite{KNSW}, Section 4.2.3), it is sufficient to reconstruct the two parameters $\hat{\alpha}_1$ and $\hat{\alpha}_2$ and calculate the respective solution to \eqref{llg3} - \eqref{llg3_abs}. This allows us to formulate the calibration problem in the reduced setting
\begin{equation}
 F(\hat{\alpha}) = F\big((\hat{\alpha}_1, \hat{\alpha}_2)\big)= v
\end{equation}
with data $v = (v_{kl})_{k=1,...,K,\, l=1,...,L} \in \mathcal{Y} := L^2([0,T])^{KL}$. The forward operator is then defined by
\begin{displaymath}
 F: \mathcal{D}(F) \subseteq \mathcal{X} \to \mathcal{Y}, \quad F(\hat{\alpha}) = \mathcal{K} \frac{\partial}{\partial t} S(\hat{\alpha})
\end{displaymath}
with the parameter-to-state map $S: \mathcal{D}(F) \subseteq \mathcal{X} := \mathbb{R}^2 \to \widetilde{\mathcal{U}}$, mapping $\hat{\alpha}$ to the corresponding solution of the LLG equation \eqref{llg3} - \eqref{llg3_abs} (see also Remark 4 in \cite{KNSW}). The observation operator $\mathcal{K}$ is defined by
\begin{displaymath}
 \mathcal{K} \mathbf{u} := \left(\int_0^{\infty} \int_{\Omega} \mathbf{K}_{kl}(t,\tau,x) \mathbf{u}(t,x) \, \mathrm{d}x \, \mathrm{d}\tau\right)_{k=1,...,K, \, l=1,...,L}.
\end{displaymath}
Here, we set
\begin{displaymath}
 \widetilde{\mathcal{U}} := H^1\big( 0,T;L^2(\Omega,\mathbb{R}^3) \big).
\end{displaymath}

\subsubsection{The all-at-once formulation}
In contrast to the reduced setting, we now aim at an all-at-once determination of the parameters $\hat{\alpha}$ and the respective magnetization $\mathbf{m}$, which avoids the direct solution of the LLG equation. \\
We split 
\begin{displaymath}
 \m := \hat{\m} + \m_0
\end{displaymath}
into the initial value $\m_0$ and the unknown rest $\hat{\m}$. The inverse problem is now formulated as
\begin{equation}
 \mathbb{F}: \mathcal{U} \times \mathcal{X} \to \mathcal{W} \times \mathcal{Y},  \ \mathbb{F}(\hat{\mathbf{m}},\hat{\alpha}) = \begin{pmatrix}
                             \mathbb{F}_0(\hat{\mathbf{m}},\hat{\alpha}) \\
                             \big(\mathbb{F}_{kl}(\hat{\mathbf{m}},\hat{\alpha})\big)_{k=1,...,K,\, l=1,...,L}
                            \end{pmatrix}
                          = \begin{pmatrix}
                             0 \\ y
                            \end{pmatrix}
                          =: \mathbf{y}  
\end{equation}
with
\begin{equation} \label{F_0aao}
\begin{split}
 \mathbb{F}_0(\hat{\mathbf{m}},\hat{\alpha}) &:= \hat{\alpha}_1 m_{\mathrm{S}}^2 \hat{\mathbf{m}}_t - \hat{\alpha}_2 \big(\hat{\mathbf{m}} + \m_0\big) \times \hat{\mathbf{m}}_t -m_{\mathrm{S}}^2\Delta_N \big(\hat{\mathbf{m}} + \m_0\big)  \\
  &\quad -\lvert \Delta_N \big(\hat{\mathbf{m}} + \m_0 \big) \rvert^2 \big(\hat{\mathbf{m}} + \m_0 \big) 
 - m_{\mathrm{S}}^2\mathbf{h} + 
\langle (\hat{\mathbf{m}} + \m_0), \mathbf{h} \rangle (\hat{\mathbf{m}}+\m_0),
\end{split}
\end{equation}
which we obtain from the LLG equation \eqref{llg3} and where $\Delta_N : H^2_N(\Omega) \to L^2(\Omega)$ with $H_N^2(\Omega) = \left\lbrace u\in H^2(\Omega) \, : \, \partial_{\nu} u = 0 \text{ on } \partial\Omega \right\rbrace$ is the Neumann-Laplacian. The observations are given by
\begin{displaymath}
 \mathbb{F}_{kl}(\mathbf{m},\hat{\alpha}) = \mathcal{K}_{kl}\mathbf{m}_t.
\end{displaymath}
We set
\begin{displaymath}
 \mathcal{U} := \left\lbrace \mathbf{u} \in  L^2\big(0,T; H_N^2(\Omega,\mathbb{R}^3) \big) \cap H^1\big( 0,T;L^2(\Omega,\mathbb{R}^3) \big) \, : \, \mathbf{u}(0)=0 \right\rbrace,
\end{displaymath}
equip $\mathcal{U}$ with the inner product
\begin{displaymath}
 (\uv_1,\uv_2)_\cU:=\int_0^T\int_\Omega \Bigl( (-\Delta_N\uv_1)\cdot(-\Delta_N\uv_2) +\uv_{1t}\cdot\uv_{2t}\Bigr)\, dx\, dt 
+\int_\Omega \nabla\uv_1(T) \colon \nabla\uv_2(T)\, dx\,.
\end{displaymath}
Furthermore, we set
\begin{displaymath}
 \mathcal{W} := H^1\big(0,T;H^1(\Omega,\mathbb{R}^3)\big)^*, 
\end{displaymath}
which is equipped with the inner product
\begin{displaymath}
\begin{aligned}
(\wv_1,\wv_2)_\cW
&:=\int_0^T\int_\Omega \Bigl( I_1[\nabla(-\Delta_N+\mbox{id})^{-1}\wv_1](t)\colon I_1[\nabla(-\Delta_N+\mbox{id})^{-1}\wv_2](t)\\
&\qquad\qquad+I_1[(-\Delta_N+\mbox{id})^{-1}\wv_1](t)\cdot I_1[(-\Delta_N+\mbox{id})^{-1}\wv_2](t)
\, dx\, dt\,
\end{aligned}
\end{displaymath}
with the isomorphism $-\Delta_N+\mbox{id} : H^1(\Omega) \to \big(H^1(\Omega)\big)^*$ and
\begin{equation} \label{I1I2}
 \begin{aligned}
I_1[w](t)&:=\int_0^t w(s)\, ds-\frac{1}{T}\int_0^T(T-s)w(s)\, ds\,,\\
I_2[w](t)&:=-\int_0^t (t-s)w(s)\, ds + \frac{t}{T}\int_0^T(T-s)w(s)\, ds\,.
\end{aligned}
\end{equation}
The well-definedness of $\mathbb{F}$ has been shown in \cite{KNSW}, Section 4.1.1.

\vspace*{1ex}

The basis for the numerical solution of the above problems is the Landweber iteration in combination with Kaczmarz' method. In the following we give a short introduction:\\
As an iterative method, the Landweber technique \cite{HankeNeubauerScherzer} initializes from a starting point $x_0$ and then runs the successive iterations
\begin{align}
 x_{k+1}= x_k - \mu_k F'(x_k)^*(F(x_k)-y^\delta) \qquad k\in\N_0,
\end{align}
where $F'(x)$ is the derivative of $F$ at $x$ and $F'(x)^*$ is its adjoint.
In case we have a collection of operators $F=(F_0,\ldots, F_{n-1}):\bigcap\limits^{n-1}_{i=0} \cD(F_i) \subset \cX \rightarrow \cY^n$ as well as data $y^\delta=(y^\delta_0,\ldots,y^\delta_{n-1})$, the Landweber-Kaczmarz method reads 
\begin{align} \label{KaczmarzIntro}
x_{k+1}= x_k -\mu_k F'_{j(k)}(x_k)^*(F_{j(k)}(x_k)-y_{j(k)}^\delta)\qquad k\in\N_0
\end{align}
with $j(k)=k-n\lfloor k/n\rfloor$, the largest integer lower or equal to $k$. We refer to \cite{HaltmeierKowarLeitaoScherzer, HaltmeierLeitaoScherzer, KowarScherzer} for a more comprehensive look on this method.

\cite{BarbaraNeubauerScherzer} suggests to terminate the Landweber iteration  according to a discrepancy principle, namely, if we let
\begin{equation} \label{mu}
 x_{k+1}= x_k - w_k h_k,
\qquad\text{where}\qquad
h_k = \mu_kF'_{j(k)}(x_k)^*(F_{j(k)}(x_k)-y_{j(k)}^\delta), 
\end{equation}
and
\begin{align} \label{w}
w_k=\begin{cases}
1 \qquad \text{if } \|F_{j(k)}(x_k)-y_{j(k)}^\delta\| \geq \tau\delta_{j(k)},\\
0 \qquad \text{otherwise},
\end{cases}
\end{align}
then the iteration stops at the first index $k_*$, for which $w_k=0$ is fulfilled in a full cycle, i.e., 
\begin{align} \label{k}
w_{k_*-i}=0, \quad i=0,\ldots, n-1 \qquad\text{and}\qquad w_{k_*-n}=1.
\end{align}
Since the data is always contaminated by noise, this is a reasonable stopping rule since we require the residual $\|F_{j(k)}(x_k)-y_{j(k)}^\delta\|$ to be of the order of the data error via some sufficiently large $\tau>2$ rather than hoping for the residual to be smaller than $\delta_k$, the noise level in the $j(k)$-th equation.

In \eqref{w}, the constant $\tau$ is chosen subject to the tangential cone condition
\begin{align*}
\|F(\tilde{x})-F(x)-F'(x)(\tilde{x}-x)\|_{\cY}\leq c_{tc}\|F(\tilde{x})-F(x)\|_{\cY} \qquad \forall x,\tilde{x}\in\cB_\rho(x_0).
\end{align*}
Additionally, the step size $\mu_k$ in \eqref{mu} is derived from uniform boundedness of the derivative of the forward operator
\begin{align*}
\|F'(x)\|_{L(\cX,\cY)}\leq C \qquad \forall x\in\cB_\rho(x_0)),
\end{align*}
namely $\mu_k \in \big(0,\frac{1}{\| F'(x_k)\|^2}\big].$\\
These conditions on $F$ are supposed to be satisfied locally in $\cB_\rho(x_0)$, the ball of center $x_0$ and radius $\rho$,  in order to ensure applicability of the iterative regularization methods \cite[Theorem 3.26]{BarbaraNeubauerScherzer}. Locally uniform smallness of $\|F'(x)\|$ has been verified for each setting in \cite{KNSW}, meanwhile the tangential cone condition or sufficient conditions for it, e.g., range invariance or adjoint range invariance \cite[Section 4.3]{BarbaraNeubauerScherzer}, are yet to be confirmed.

\subsection{Derivatives and adjoints}
In order to solve the inverse problem as formulated in the two preceding sections, we need the respective Fr\'echet derivatives as well as their adjoints, see also \cite{BarbaraNeubauerScherzer}. In this chapter, we present the algorithm for the unit length magnetization, i.e. $m_{\mathrm{S}}=|\mv|=1.$

\subsubsection{The all-at-once setting}

The Fr\'echet derivative of $\mathbb{F}$ (i.e., the forward operator in the all-at-once setting) is given by
\begin{displaymath}
 \begin{aligned}
&\bF'(\mh,\altil_1,\altil_2)(\uv,\beta_1,\beta_2)\\
&=\left(\begin{array}{l}
\beta_1 \mh_t - \beta_2 (\mv_0+\mh)\times \mh_t \\
\hspace*{1.5cm}
+\altil_1 \uv_t-\Delta_N \uv - \altil_2 \uv\times \mh_t-\altil_2 (\mv_0+\mh)\times \uv_t
 \\
\hspace*{1.5cm}-2(\nabla (\mv_0+\mh)\colon\nabla\uv)(\mv_0+\mh)- |\nabla (\mv_0+\mh)|^2 \uv
 \\
\hspace*{1.5cm}
+((\mv_0+\mh)\cdot\hv) \uv 
+ (\uv\cdot\hv) (\mv_0+\mh)
\\[2ex]
\Bigl(\int_0^\Tt \int_\Omega \ker_{k\ell}(t,\tau,x) \cdot\uv_t(x,\tau)\, dx\, d\tau\Bigr)_{k=1,\ldots,K\,, \ \ell=1,\ldots,L}
\end{array}\right)\\
&=\left(\begin{array}{ccc}
\frac{\partial\bF_0}{\partial\mh}(\mh,\altil)&\frac{\partial\bF_0}{\partial\altil_1}(\mh,\altil)&\frac{\partial\bF_0}{\partial\altil_2}(\mh,\altil)\\
(\frac{\partial\bF_{k\ell}}{\partial\mh}(\mh,\altil))_{k=1,\ldots,K,\ell=1,\ldots,L}&0&0
\end{array}\right)
\left(\begin{array}{c} \uv \\ \beta_1\\ \beta_2\end{array}\right).
\end{aligned}
\end{displaymath}
Its adjoint is calculated componentwise. First of all, we have that
\begin{displaymath}
 \frac{\partial\bF_0}{\partial\mh}(\mh,\altil)^*\yv=:\zv
\end{displaymath}
is the solution of the two auxiliary problems
\begin{equation}\label{auxprob1}
\left\{\begin{array}{rcl}
\zv_t-\Delta \zv&=&\vv\mbox{ in }(0,T)\times\Omega\\
\partial_\nu\zv&=&0\mbox{ on }(0,T)\times\partial\Omega\\
\zv(0)&=&0\mbox{ in }\Omega
\end{array}\right.
\end{equation}
and
\begin{equation}\label{auxprob2}
\left\{\begin{array}{rcl}
-\vv_t-\Delta \vv&=&\mathbf{f}\mbox{ in }(0,T)\times\Omega\\
\partial_\nu\vv&=&0\mbox{ on }(0,T)\times\partial\Omega\\
\vv(T)&=&\mathbf{g}\mbox{ in }\Omega
\end{array}\right.
\end{equation}
with
\begin{equation} \label{fygy}
\begin{split}
 \mathbf{f} = \mathbf{f}^{\mathbf{y}} &:= -\altil_1 \yv_t+ (-\Delta\yv) - \altil_2 \mh_t\times\yv+\altil_2 \yv_t\times(\mv_0+\mh)+\altil_2 \yv\times\mh_t \\
 &\quad-2((\mv_0+\mh)\cdot\yv)\, (-\Delta_N(\mv_0+\mh)) +2((\nabla (\mv_0+\mh)^T(\nabla \yv))\, (\mv_0+\mh)\\
 &\quad +2((\nabla (\mv_0+\mh)^T(\nabla (\mv_0+\mh)))\, \yv -|\nabla (\mv_0+\mh)|^2 \yv \\
 &\quad+((\mv_0+\mh)\cdot\hv) \,\yv + ((\mv_0+\mh)\cdot\yv) \, \hv \Bigr) \\
 \mathbf{g} = \mathbf{g}^{\mathbf{y}}_T &:= \altil_1\yv(T) -\altil_2 \yv(T)\times(\mv_0+\mh(T))
\end{split}
\end{equation}
and $\mathbf{y} = I_2[\widetilde{y}]$, where $\widetilde{y}(t)$ solves
\begin{equation} \label{pde_ytil}
 \left\{\begin{array}{rcl}
-\Delta \widetilde{y}(t)+\widetilde{y}(t)&=&\wv(t)\mbox{ in }\Omega\\
\partial_\nu\widetilde{y}&=&0\mbox{ on }\partial\Omega
\end{array}\right.
\end{equation}
for each $t \in (0,T)$. For the same $\mathbf{y}$ we have
\begin{align}
\frac{\partial\bF_0}{\partial\altil_1}(\mh,\altil)^*\wv &= \int_0^T\int_\Omega \mh_t\cdot \yv\, dx\, dt\,, \label{Daltil1}\\
\frac{\partial\bF_0}{\partial\altil_2}(\mh,\altil)^*\wv &= -\int_0^T\int_\Omega ((\mv_0+\mh)\times\mh_t)\cdot \yv\, dx\, dt\,. \label{Daltil2}
\end{align}
Finally, the adjoint
\begin{displaymath}
 (\frac{\partial\bF_{k\ell}}{\partial\mh}(\mh,\altil))_{k=1,\ldots,K,\ell=1,\ldots,L}^*y=\zv
\end{displaymath}
is calculated by solving \eqref{auxprob1}, \eqref{auxprob2} with
\begin{equation} \label{frgr}
 \begin{split}
\mathbf{f}(x,\tau)&=-\int_0^T\sum_{k=1}^K\sum_{\ell=1}^L \frac{\partial}{\partial\tau}\ker_{k\ell}(t,\tau,x) y_{k\ell}(t)\ dt,\\
\mathbf{g}(x)&=\int_0^T\sum_{k=1}^K\sum_{\ell=1}^L \ker_{k\ell}(t,\Tt,x) y_{k\ell}(t)\ dt\,.
\end{split}
\end{equation}

\subsubsection{The reduced setting}
The Fr\'echet derivative of $F$ in $\hat{\alpha}$ is given by $F'(\hat{\alpha})\beta = \mathcal{K} \mathbf{u}_t$, where $\uv$ solves the \emph{linearized Landau-Lifshitz-Gilbert equation}
\begin{alignat*}{3}
\altil_1\u_t &- \altil_2\m\times\u_t - \altil_2\u\times\m_t - \Delta\u -2(\nabla\u:\nabla\m)\m\\
&\qquad+\u(-|\nabla\m|^2+(\m\cdot \h)) + (\u\cdot \h)\m &&\\
&=-\beta_1\m_t+\beta_2\m\times\m_t \qquad &&\text{in } (0,T)\times\Omega\\
 \partial_\nu\u &=0 && \text{on } (0,T)\times\partial\Omega\\
 \u(0)&=0 && \text{in } \Omega.
\end{alignat*}
The respective Hilbert space adjoint 
\begin{displaymath}
 F'(\altil)^*:L^2(0,T)^{KL}\rightarrow\R^2
\end{displaymath}
has been shown to be of the form
\begin{align} \label{red-adjoint}
F'(\altil)^*z=\left( \int_0^T\int_\Omega -\m_t\cdot\p \,dx\,dt , \int_0^T\int_\Omega(\m\times\m_t)\cdot\p \,dx\,dt\right),
\end{align}
where $\p$ solves the adjoint equation
\begin{alignat}{3}
&-\altil_1\p_t - \altil_2\m\times\p_t - 2\altil_2\m_t\times\p - \Delta\p \nonumber\\
&\quad +2\left((\nabla\m)^\top\nabla\m\right)\p + 2\left((\nabla\m)^\top\nabla\p\right)\m \nonumber\\
&\quad+(-|\nabla\m|^2+(\m\cdot \h))\p + (\m\cdot\p)(\h + 2\Delta\m)=\Ktil z  \qquad&&\text{in } (0,T)\times\Omega \label{red-adjoint-eq0-1}\\ 
& \partial_\nu\p=0 && \text{on } (0,T)\times\partial\Omega \label{red-adjoint-eq0-2}\\
& \altil_1\p(T)+\altil_2(\m\times\p)(T)=\KtilT z && \text{in } \Omega. \label{red-adjoint-eq0-3}
\end{alignat}
with
\begin{align}
&(\Ktil r)(x,t)=\sum_{k=1}^K\sum_{\ell=1}^L-\mu_0c_k(x)\pr(x)\int_0^T \widetilde{a}_{\ell\,\tau}(\tau-t)r_{k\ell}(\tau)\,d\tau, \quad t\in(0,T), \label{eqKtil1}\\
&(\KtilT r)(x)=\sum_{k=1}^K\sum_{\ell=1}^L-\mu_0c_k(x)\pr(x)\int_0^T \widetilde{a}_{\ell}(\tau)r_{k\ell}(\tau)\,d\tau.\label{eqKtil2}
\end{align}

\section{Algorithms: parameter identification for the LLG equation}\label{topci4-sec:algorithm}
Let us formulate the problem in both settings into one general form $F:\cX\rightarrow\cY, \ F(x)=y$, where in practice, only an approximation $y^\delta$ of $y$ is measurable. Regarding the all-at-once setting this means $F:=\bF, x:=(\mh,\altil_1,\altil_2)$, and for the reduced setting one has $F:=F, x:=(\altil_1,\altil_2)$.

\vspace*{1ex}

We now explicitly present all steps comprised in each of these methods in the all-at-once version as well as its counterpart reduced version. For the all-at-once setting, the algorithm relies on Section 4.1, and for the reduced setting, it results from Section 4.2, both sections are in \cite{KNSW}.

\subsection{All-at-once Landweber}\label{sec-algorithm-AAO}
In the following, we formulate the algorithm corresponding to the Landweber iteration for the all-at-once setting.

\begin{algo} \label{algo_aaoL}
Starting from an initial guess $(\mh, \altil_1, \altil_2)_{j=0}=(\m-\m_0, \altil_1, \altil_2)_{j=0}$, run the following steps:
\begin{enumerate}[label=S.\arabic*.]
\item  Set argument to adjoint equations (see \eqref{F_0aao}): \label{S1}\\[2ex]
Let $\mh:=\mh_j, \altil_1:=\altil_{1\,j}, \altil_2:=\altil_{2\,j}$.\\
Compute $\wv = \mathbb{F}_0(\mh,\altil)$ and the residual $r(t)$ by
\begin{align*}
\wv&= \altil_1 \mh_t-\Delta_N (\m_0+\mh) - \altil_2 (\m_0+\mh)\times \mh_t\\ 
&\qquad\quad-|\nabla (\m_0+\mh)|^2(\m_0+\mh) - \hv +((\m_0+\mh)\cdot\hv)(\m_0+\mh)\\
r(t) &=\left(\int_0^T \int_\Omega \ker_{k\ell}(x,t,\tau) \mv_t(x,\tau)\, dx\, d\tau\right)_{k=1,\ldots,K, \ell=1,\ldots,L} - y^\delta(t).
\end{align*}
Check: Stopping rule according to discrepancy principle.

\item  Compute the adjoints: \label{S2}
\begin{enumerate}[label=\Alph*.]

\item Compute $\zv= \frac{\partial\bF_0}{\partial\mh}(\mh,\altil)^*\wv$ \label{A}
\begin{enumerate}[label=A.\arabic*.]
\item Input: $\wv$ \label{A.1}\\
Solve \eqref{pde_ytil}, i.e.,
\begin{alignat*}{3}
-\Delta \widetilde{y}(t)+\widetilde{y}(t)&=\wv(t) \quad &&\mbox{ in }\Omega\\
\partial_\nu\widetilde{y}&=0 &&\mbox{ on }\partial\Omega.
\end{alignat*}
Output: $\tilde{y}$
\item Input: $\tilde{y}$ \label{A.2}\\
Compute (according to \eqref{I1I2})
\[
\yv(t)=I_2[\tilde{y}](t)=-\int_0^t (t-s)\tilde{y}(s)\, ds + \frac{t}{T}\int_0^T(T-s)\tilde{y}(s)\, ds.\
\]
Output: $\yv$
\item Input: $\yv$ \label{A.3}\\
Compute (see \eqref{fygy})
\[
\begin{aligned}
\f^{\yv}
&=-\altil_1 \yv_t+ (-\Delta_N\yv) - \altil_2 \mh_t\times\yv+\altil_2 \yv_t\times(\m_0+\mh)+\altil_2 \yv\times\mh_t\\
&\qquad -2((\m_0+\mh)\cdot\yv)\, (-\Delta_N(\m_0+\mh))\\
&\qquad +2((\nabla (\m_0+\mh)^T(\nabla \yv))\, (\m_0+\mh)\\
&\qquad +2((\nabla (\m_0+\mh)^T(\nabla (\m_0+\mh)))\, \yv
-|\nabla (\m_0+\mh)|^2 \yv
 \\
&\qquad+((\m_0+\mh)\cdot\hv) \,\yv 
+ ((\m_0+\mh)\cdot\yv) \, \hv\\
\g^{\yv}_T&=\altil_1\yv(T) -\altil_2 \yv(T)\times(\m_0+\mh(T)).
\end{aligned}
\]
Output: $\f^{\yv}, \g_T^{\yv}$
\item Input: $\f^{\yv}, \g_T^{\yv}$ \label{A.4}\\
Solve \eqref{auxprob2}
\begin{alignat*}{3}
-\vv_t-\Delta \vv&=\f^{\yv} \quad&&\mbox{ in }(0,T)\times\Omega\\
\partial_\nu\vv&=0 &&\mbox{ on }(0,T)\times\partial\Omega\\
\vv(T)&=\g^{\yv}_T &&\mbox{ in }\Omega.
\end{alignat*}
Output: $\vv$
\item Input: $\vv$ \label{A.5}\\
Solve \eqref{auxprob1}
\begin{alignat*}{3}
\zv_t-\Delta \zv&=\vv \quad&&\mbox{ in }(0,T)\times\Omega\\
\partial_\nu\zv&=0 &&\mbox{ on }(0,T)\times\partial\Omega\\
\zv(0)&=0 &&\mbox{ in }\Omega.
\end{alignat*}
Output: $\zv$
\end{enumerate}

\item Compute $\s=\left(\frac{\partial\bF_{k\ell}}{\partial\mh}(\mh,\altil)_{k=1,\ldots,K,\ell=1,\ldots,L}\right)^*r$ \label{B}
\begin{enumerate}[label=B.\arabic*.]
\item Input: $r$ from Step \ref{S1}\\
Compute
\begin{align*}
&\mathbf{f}^r(x,\tau)=-\int_0^T\sum_{k=1}^K\sum_{\ell=1}^L \frac{\partial}{\partial\tau}\ker_{k\ell}(t,\tau,x) r_{k\ell}(t)\ dt,\\
&\mathbf{g}^r(x)=\int_0^T\sum_{k=1}^K\sum_{\ell=1}^L \ker_{k\ell}(t,\Tt,x) r_{k\ell}(t)\ dt,
\end{align*}
according to \eqref{frgr} with
\[
\mathbf{K}_{k\ell}(t,\tau,x) := -\mu_0\widetilde{a}_\ell(t - \tau ) c_k(x) \pR_\ell(x).
\]
Output: $\f^r,\g^r$
\item Step \ref{A.4} with input: $\f^r,\g^r$
\item Step \ref{A.5} with output: $\s$
\end{enumerate}

\item Compute $\beta_{1}=\frac{\partial\bF_0}{\partial\altil_1}(\mh,\altil)^*\wv$ \label{C} according to \eqref{Daltil1}.
\begin{enumerate}[label=C.\arabic*.]
\item Input: $\yv$ from Step \ref{A.2}\\
Compute
\begin{align*}
\beta_{1}= \int_0^T\int_\Omega \mh_t\cdot \yv\, dx\, dt
\end{align*}
Output: $\beta_1$
\end{enumerate}

\item Compute $\beta_{2}=\frac{\partial\bF_0}{\partial\altil_2}(\mh,\altil)^*\wv$ \label{D} according to \eqref{Daltil2}.
\begin{enumerate}[label=D.\arabic*.]
\item Input: $\yv$ from Step \ref{A.2}\\
Compute
\begin{align*}
\beta_{2} = -\int_0^T\int_\Omega ((\mv_0+\mh)\times\mh_t)\cdot \yv\, dx\, dt\,. 
\end{align*}
Output: $\beta_2$
\end{enumerate}

\end{enumerate}

\item Update $\mh, \altil_1, \altil_2$ with step size $\mu$: \label{S3}
\begin{align}
\mh_{j+1}&=\mh_j - \mu (\zv + \s) \label{S3_1} \\
\altil_{1\,j+1}&=\altil_{1\,j} - \mu \beta_1 \label{S3_2}\\
\altil_{2\,j+1}&=\altil_{2\,j} - \mu \beta_2.	\label{S3_3}
\end{align}

\end{enumerate}
\end{algo}

\begin{figure}[!h] 
\setlength{\belowcaptionskip}{-5pt}
\centering
\includegraphics[width=8cm]{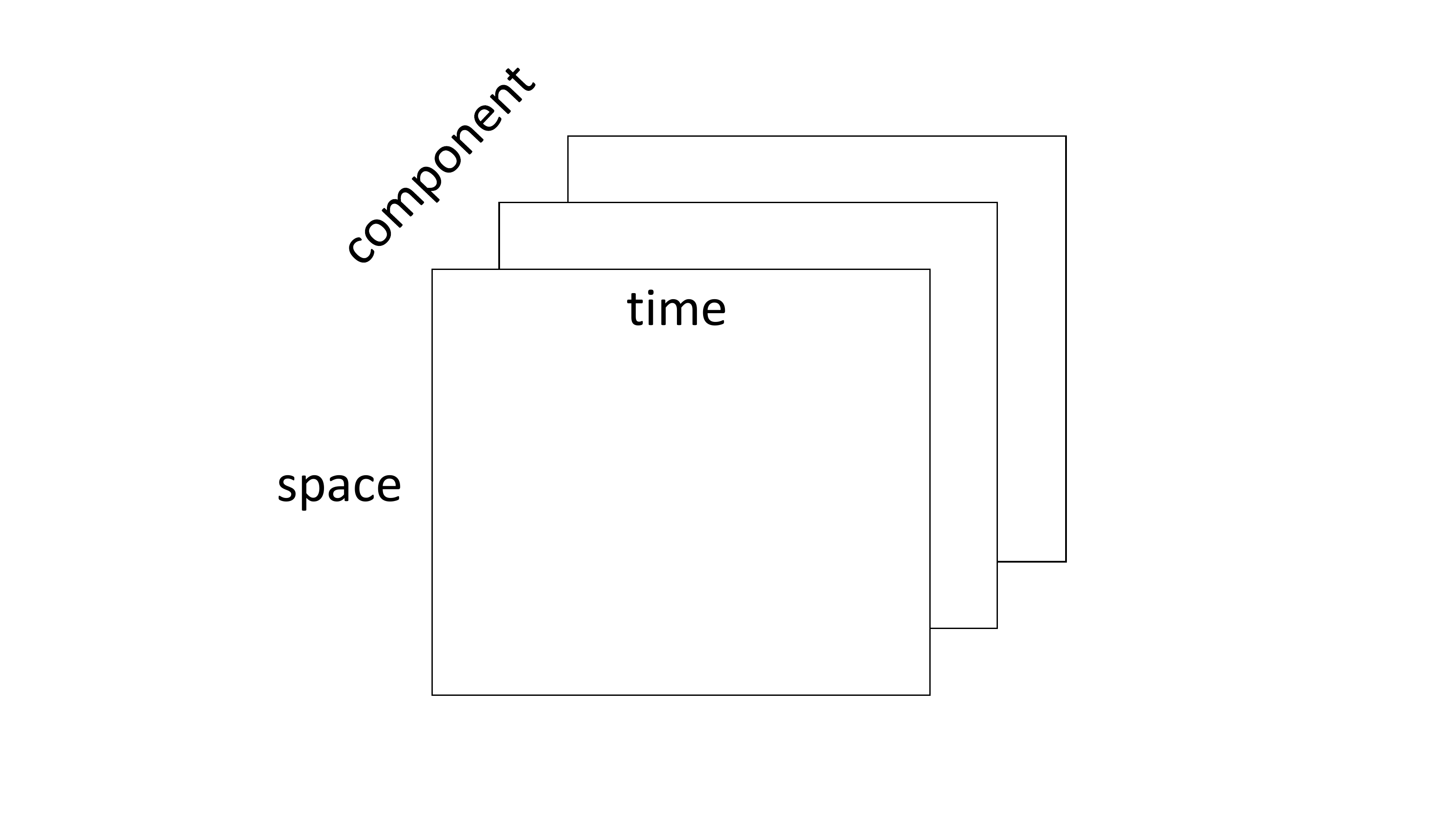}
\caption{Matrix representation for a vector field in the all-at-once setting.}
\label{matrix-aao}
\end{figure}

In the implementation, one can consider each of the vector fields as a three-dimensional matrix as illustrated in Figure \ref{matrix-aao}. Steps \ref{A.1}, \ref{A.2}, \ref{A.4}, \ref{A.5} then operate on each time-space slice, meanwhile Step \ref{A.3} needs to be calculated among 3D-matrices.

\begin{remark} \label{LLGsolver}
In case the values of the parameters $\altil_1, \altil_2$ are known, one can omit Steps \ref{B}, \ref{C}, \ref{D} and Steps \eqref{S3_2}-\eqref{S3_3} in \ref{S3} to obtain a numerical solution to $\mv$. This procedure can be seen as an LLG solver.
\end{remark}

As far as a numerical solution of the Landau-Lifshitz-Gilbert equation equation is concerned, numerous investigations have been recently published. \blue{In the seminal work \cite{BartelsProhl} in 2006, Bartels and Prohl considered the LLG equation involving a 2-harmonic heat flow, and established a numerical implementation strategy based on solving a nonlinear system at each time step. This algorithm integrates a fixed-point method with stopping criteria (to handle the nonlinear system) into a lowest order conforming finite element method.
In \cite{BartelsProhl1}, these authors extended the results for the case of $p$-harmonic heat flow with $1<p<\infty$.} The work of Alouges et.~al.~in 2008, see \cite{alouges}, describes a new implicit finite element scheme, which avoids solving nonlinear systems. This $\theta$-scheme introduces the new term $\vv:=\m_t$ to form a linear equation in $\vv$. At each time step $n$, $\theta\in[0,1]$ involves $\nabla(\m_n+\theta k\vv_n)$ in the variation formula (hence implicit), where $k$ is the time step size. After solving $\vv_n$ by an implicit finite element scheme, $\m$ at the next time step $(n+1)$ is updated by $\m_{n+1}:=\frac{\m_n+k\vv_n}{|\m_n+k\vv_n|}$. Inspired by the $\theta$-scheme, Ba\v nas et.~al.~in 2014, see \cite{BanasEtAl}, studied the more general model including the magnetostriction effect instead of the magnetostatic simplification. The phenomenon is governed by a coupled problem of an LLG equation and a second time-dependent PDE representing the magnetostrictive contribution. The authors later dealt with the full Maxwell-LLG equation in \cite{BanasPagePraetorius} and proposed a fully decoupled scheme. In 2014, Alouges et.~al.~\cite{alougesKritsikisSteinerToussaint} upgraded the original first order $\theta$-scheme by replacing the tangential update for $\m$ by a higher approximation via $\text{Proj}_{\m^\perp}$ and $\m_{tt}$, creating  a new (almost) order two $\theta$-scheme. The algorithm initializes with differentiating with respect to time and the LLG equation then ends up with a variation formula linear in $\vv$. The approximation in space is still of order one ($P^1$ Lagrange finite elements), and the convergence of the scheme does not hold for higher order elements.

\begin{remark} \label{LLGsolver2}
Independently of these existing works, our algorithm is  attractive from a practical point of view since we also require to solve only linear PDEs per iteration step. Our contribution as well as the advantage of the method can be summarized as follows:
\begin{itemize}
\item 
Unlike the mentioned novel schemes that find $\mv$ at each time point per iteration step, our scheme calculates $\mv$ at all time instances per iteration. The loop in our scheme is dedicated to the Landweber(-Kaczmarz) iteration, which improves the whole $\mv([0,T]\times\Omega)$ gradually. This shall grant access to space-time adaptive discretization.
\item
In each step, only three separate and conventional linear PDEs, i.e., \ref{A.1}, \ref{A.4}, and \ref{A.5} are required to be solved. Our method, therefore, does not need to derive new theory for proving unique existence and convergence of finite element solutions, which are the main results in most of the current literature.
\item
Also for this reason, the suggested method is favorable in implementation since one can make use of existing standard finite difference or finite element codes.
\item
We also remark that our method is able to be upgraded to higher order through
\begin{itemize}
\item higher order standard FD/FE for solving the conventional PDEs in \ref{A.1}, \ref{A.4}, \ref{A.5} (feasible).
\item higher order numerical approximation for the integrals and derivatives in \ref{A.2}, \ref{A.3} (feasible).
\item a higher convergence rate of the Landweber. This fact depends on how smooth the exact solution is (source condition, see, e.g., \cite[Section 2.3]{BarbaraNeubauerScherzer}). And this turns out to be the smoothness of $\m_0$ and $\hv$, which is feasible if one proceeds similarly as in the proof of regularity of $\mv$ in the reduced setting (c.f., \cite[Section 4.2.3]{KNSW}).
\end{itemize}
\item
Although not being presented here, our method can be extended to the case when anisotropy, represented by a term $\pii(\mv)$ in the effective magnetic field, is involved, under certain conditions on $\pii$. This is due to the fact that after calculating $\bF'(\mh,\altil_1,\altil_2)^*\wv$, only the term $\pi'(\mh)^*\wv$ is added to $\f^{\yv}$ in Step \ref{A.3}, while all other steps remain unchanged.
\item
As the examples illustrated in this paper are in a simple one-dimensional domain, a finite difference discretization is well adapted. On the other hand, there is no reason preventing the proposed method to be implemented by finite elements to better suit complex geometries. 
\item
Concerning memory requirement, in each of the Landweber iterations - if using an implicit Euler time stepping scheme - the proposed method demands storage for computing the cross product of two matrices of size nx$\times$3 (number of space grids $\times$ number of components) or for the multiplication of a matrix of size nx$\times$nx (finite difference matrices or stiffness matrix) with a vector of nx elements. This yields similarity in the memory requirement with the other existing schemes.
\end{itemize}
However:
\begin{itemize}
\item
Due to the nature of nonlinear inverse problems, our scheme is able to solve it just locally, i.e., a sufficiently good initial guess needs to be known.
\item
Concerning memory requirement, our algorithm demands more memory than the others as the whole $\mv([0,T]\times\Omega)$ needs to be allocated in RAM for each Landweber iteration.
\end{itemize}
\end{remark}

\subsection{Reduced Landweber}
We now present the respective Landweber algorithm for the numerical solution of our inverse problem in the reduced setting.

\begin{algo} \label{algo_redL}
Starting from an initial guess $\altil_{j=0}=(\altil_1, \altil_2)_{j=0}$, run the following steps:
\begin{enumerate}[label=S.\arabic*.]
\item Compute the state $\m:=S(\altil_j)$ according to the LLG equation.\label{red-S1}

\item Set argument to the adjoint equation.\label{red-S2}\\
Compute the residual
\begin{align*}
r(t) &=\left(\int_0^T \int_\Omega \ker_{k\ell}(x,t,\tau) \mv_\tau(x,\tau)\, dx\, d\tau\right)_{k=1,\ldots,K, \ell=1,\ldots,L} - y^\delta(t).
\end{align*}
Check: Stopping rule according to discrepancy principle.

\item Compute the adjoint state $\p=F'(\altil)^*r$ as in \eqref{red-adjoint-eq0-1} - \eqref{red-adjoint-eq0-3} according to \label{red-S3}
\begin{alignat}{3}
&-\altil_1\p_t - \altil_2\m\times\p_t - 2\altil_2\m_t\times\p - \Delta\p \nonumber\\
&\quad +2\left((\nabla\m)^\top\nabla\m\right)\p + 2\left((\nabla\m)^\top\nabla\p\right)\m \nonumber\\
&\quad+(-|\nabla\m|^2+(\m\cdot \h))\p + (\m\cdot\p)(\h + 2\Delta\m)=\Ktil r  \,&&\text{ in } (0,T)\times\Omega \label{topic4-red-adjoint-eq0-1}\\
& \partial_\nu\p=0 && \text{ on } (0,T)\times\partial\Omega \label{topic4-red-adjoint-eq0-2}\\
& \altil_1\p(T)+\altil_2(\m\times\p)(T)=\KtilT r && \text{ in } \Omega \label{topic4-red-adjoint-eq0-3}
\end{alignat}
with 
\begin{align*}
&(\Ktil r)(x,t)=\sum_{k=1}^K\sum_{\ell=1}^L-\mu_0c_k(x)\pr(x)\int_0^T \widetilde{a}_{\ell\,\tau}(\tau-t)r_{k\ell}(\tau)\,d\tau, \quad t\in(0,T),\\
&(\KtilT r)(x)=\sum_{k=1}^K\sum_{\ell=1}^L-\mu_0c_k(x)\pr(x)\int_0^T \widetilde{a}_{\ell}(\tau)r_{k\ell}(\tau)\,d\tau
\end{align*}
from \eqref{eqKtil1}, \eqref{eqKtil2}.
\item Update $\altil_1, \altil_2$ with step size $\mu$ according to \eqref{red-adjoint}: \label{red-S4}
\begin{align}
\altil_{1\,j+1}&=\altil_{1\,j} - \mu\int_0^T\int_\Omega -\m_{t}\cdot\p \,dx\,dt \label{red-S4-1}\\
\altil_{2\,j+1}&=\altil_{2\,j} - \mu\int_0^T\int_\Omega(\m\times\m_{t})\cdot\p \,dx\,dt.	\label{red-S4-2}
\end{align}
\end{enumerate}
\end{algo}

\begin{figure}[!h] 
\centering
\includegraphics[width=9cm]{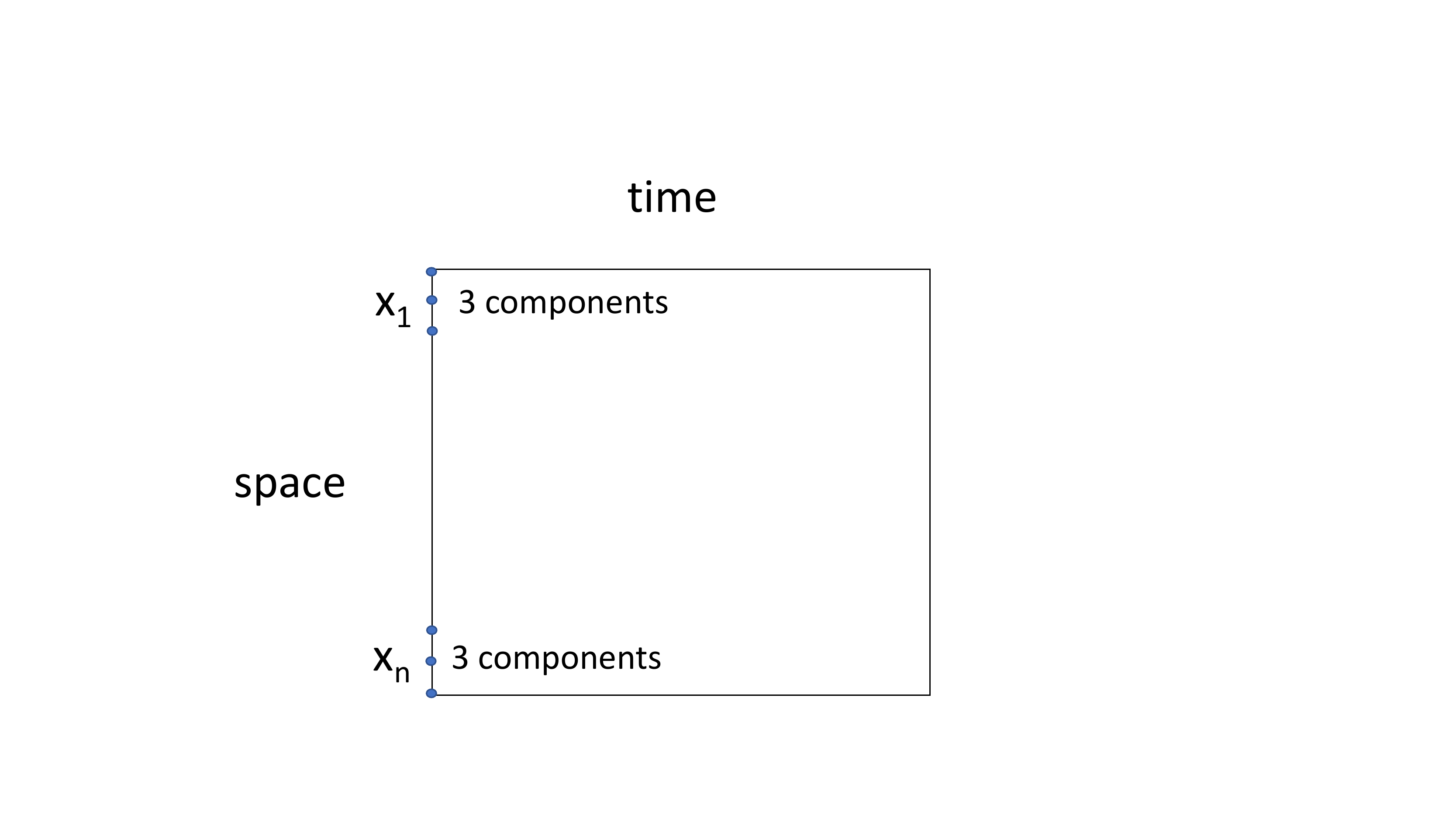}
\caption{Matrix representation for a vector field in the reduced setting.}
\label{matrix-red}
\end{figure}

In contrast to the all-at-once setting, which solves the conventional PDEs on time-space slices, the PDE in the adjoint equation (Step \ref{red-S3}) of the reduced one involves cross products and moduli, it therefore interacts between three components. Figure \ref{matrix-red} illustrates an example for assigning a vector field to a matrix when dealing with implementing the reduced version. One might reshape the matrix in Figure \ref{matrix-red} into one vector of $nx\times3\times nt$ elements.

For Step \ref{red-S2} here, one can employ the LLG solver introduced in Remark \ref{LLGsolver}. In order to search for $\m$, the program needs an initial guess, which can be chosen as the initial state $\m_0$. Then the computed $\m_j$ of the Landweber iteration $j$, besides being the input to the next steps \ref{red-S3}-\ref{red-S4} in the stream, on the other hand, plays the role of an initial guess for \ref{red-S2} in the next iteration $(j+1)$.

The reduced setting is supposed to run more slowly than the all-at-once one, if using the same step size, as each of the Landweber iteration calls an LLG solver leading to an additional inner loop.

Given the time-dependence of the data as well as its dependence on the receive coil $l$ and the sample $k$, it is possible to formulate our inverse problem in a semi-discrete setting, i.e., a collection of sub-problems related to the dependencies in the data. This is especially attractive for problems with large data sets or high dimensional model operator and calls for an application of Kaczmarz' method (see also \cite{natterer}, \cite{Nguyen}, or \cite{bhw20}), which successively sweeps through each of those sub-problems in each iteration. \\
Here, we address two independent approaches: One is based on time-segmenting, i.e., the sub-problems correspond to time-segments of the time interval during which the data acquisition takes place. The other option is based on data selection, i.e., for each combination of receive coil $l$ and calibration concentration $c_k$, we obtain one sub-problem (and $K\cdot L$ sub-problems in total). \\
The algorithms for those schemes are detailed in the following remarks.

\begin{algo}{\bf Kaczmarz based on time segmenting}\\
Starting from $\altil_{j=0}=(\altil_1, \altil_2)_{j=0}$, run:
\begin{enumerate}[label=S.\arabic*.]
\item Compute the state $\m:=S(\altil_j)$ according to the LLG equation
\item Set argument to adjoint equation\\
Compute
\begin{align*}
r(t)=\left(\int_{t^j}^{t^{j+1}}\int_\Omega \ker_{k\ell}(x,t,\tau) \mv_\tau(x,\tau)\, dx\, d\tau\right)_{k=1,\ldots,K, \ell=1,\ldots,L} - y^\delta(t)\chi_{[t^j,t^{j+1}]}(t)
\end{align*}
with $t^j=\left\lfloor \dfrac{j}{n} \right\rfloor$  and $ 0=t^0<\ldots<t^{n-1}=T$.\\[1ex]
Check: Stopping rule according to discrepancy principle.
\item Compute the adjoint state $\p=F'(\altil)^*r$ according to \eqref{topic4-red-adjoint-eq0-1}-\eqref{topic4-red-adjoint-eq0-3}
with
\begin{align*}
&(\Ktil r)(x,t)=\sum_{k=1}^K\sum_{\ell=1}^L-\mu_0c_k(x)\pr(x)\int_{t^j}^{t^{j+1}} \widetilde{a}_{\ell\,\tau}(\tau-t)r_{k\ell}(\tau)\,d\tau \qquad t\in(0,T)\\
&(\KtilT r)(x)=\sum_{k=1}^K\sum_{\ell=1}^L-\mu_0c_k(x)\pr(x)\int_{t^j}^{t^{j+1}} \widetilde{a}_{\ell}(\tau)r_{k\ell}(\tau)\,d\tau.
\end{align*}

\item Update $\altil_1, \altil_2$ as \eqref{red-S4-1}-\eqref{red-S4-2}.
\end{enumerate}
\end{algo}

\begin{algo}{\bf Kaczmarz based on data selection}\\
Starting from $\altil_{j=0}=(\altil_1, \altil_2)_{j=0}$, run:
\begin{enumerate}[label=S.\arabic*.]
\item Compute the state $\m:=S(\altil_j)$ according to the LLG equation
\item Set argument to adjoint equation\\
Compute
\begin{align*}
r(t)=\int_0^T\int_\Omega \ker_{k\ell}(x,t,\tau) \mv_\tau(x,\tau)\, dx\, d\tau - y_{k\ell}^\delta(t)
\end{align*}
with $k,\ell$ satisfying $k\ell=\left\lfloor \dfrac{j}{KL}\right\rfloor$.\\[1ex]
Check: Stopping rule according to discrepancy principle.
\item Compute the adjoint state $\p=F'(\altil)^*r$ according to \eqref{topic4-red-adjoint-eq0-1}-\eqref{topic4-red-adjoint-eq0-3}
with
\begin{align*}
&(\Ktil r)(x,t)=-\mu_0c_k(x)\pr(x)\int_0^T \widetilde{a}_{\ell\,\tau}(\tau-t)r(\tau)\,d\tau \qquad t\in(0,T)\\
&(\KtilT r)(x)=-\mu_0c_k(x)\pr(x)\int_0^T \widetilde{a}_{\ell}(\tau)r(\tau)\,d\tau.
\end{align*}

\item Update $\altil_1, \altil_2$ as \eqref{red-S4-1}-\eqref{red-S4-2}.
\end{enumerate}
\end{algo}

\section{Numerical experiments}\label{topci4-sec:example}
This section is dedicated to a range of numerical experiments. We start by presenting a numerical analysis of the LLG solver we obtain from Algorithm \ref{algo_aaoL} by implementing Remark \ref{LLGsolver} for examples with a known ground truth (Section \ref{sec-LLGsolver}). Subsequently, in Section \ref{sec:LLGsyn}, we use this Algorithm to simulate the evolution of the magnetization in a scenario inspired by the actual physical setting. The full reconstruction algorithms \ref{algo_aaoL} and \ref{algo_redL} are applied to synthetic data in Section \ref{Sec:recos} for an evaluation of the algorithms' performance.

\subsection{LLG solver}\label{sec-LLGsolver}
In this section, we examine the performance of the proposed solver for the following LLG equation
\begin{equation}\label{LLG-dropMs}
\begin{aligned}
\hat{\alpha}_1  \mathbf{m}_t - \hat{\alpha}_2 \mathbf{m} \times \mathbf{m}_t -\Delta \mathbf{m} &= \lvert \nabla \mathbf{m} \rvert^2 \mathbf{m} + \mathbf{h} - 
\langle \mathbf{m}, \mathbf{h} \rangle \mathbf{m}  && \text{in} \ [0,T] \times \Omega, \\
 0 &= \partial_{\nu} \mathbf{m}  && \text{on} \ [0,T] \times \partial\Omega, \\
 \mathbf{m}_0 &= \mathbf{m}(t=0), \ \lvert \mathbf{m}_0 \rvert = m_{\mathrm{S}} && \text{in} \ \Omega \,.
\end{aligned}
\end{equation}

When using Algorithm \ref{algo_aaoL} in the way mentioned in Remark \ref{LLGsolver} as an LLG solver, there is no noise involved in the process as the exact data $(y=0)$ is known. Hence, the algorithm should be run with a number of iterations as large as possible in order to reach an acceptable accuracy for the reconstructed $\mh$.

In the following, we shall numerically test this algorithm using the finite difference method for Steps \ref{A.1}, \ref{A.3}-\ref{A.5}. In particular, central difference quotients were employed to approximate time and space derivatives. The numerical integration in Step \ref{A.2} runs with the trapezoidal rule. For computing the respective $L^2$-norms, we use also the trapezoidal rule. For time discretization, the interval $[0,0.2]$ is partitioned into 51 time steps, and for the space domain $[0,2\pi]$, we impose a discretization of 101 grid points. The method in use is Landweber. 

We now analyze the numerical performance of the LLG solver by means of three test cases specified in Table \ref{tab-1}. Test 1 with the corresponding parameter set is indeed a true solution to the LLG equation \ref{LLG-dropMs}. Test 2, however, does not completely fulfill the original initial boundary value problem for the LLG equation since in the LLG model, $\m_0$ is supposed to have constant length over $\Omega$. For the same reason, initial data in Test 3 - although fulfilling the requirement of constant length ($|\m_0(x)|=2$) - does not have a homogeneous Neumann boundary. Nevertheless, those test cases still reflect the all-at-once model with unknown $\mh$, while $\m_0$ is being considered as just an additive term, and thus are still recognized as meaningful examples helping to increase the diversity of the experiments. Please note that the chosen quantities have no physical background here, which is why we do not use any units.

\begin{figure}[p] 
\vspace{2cm}
\centering
\includegraphics[width=4.3cm]{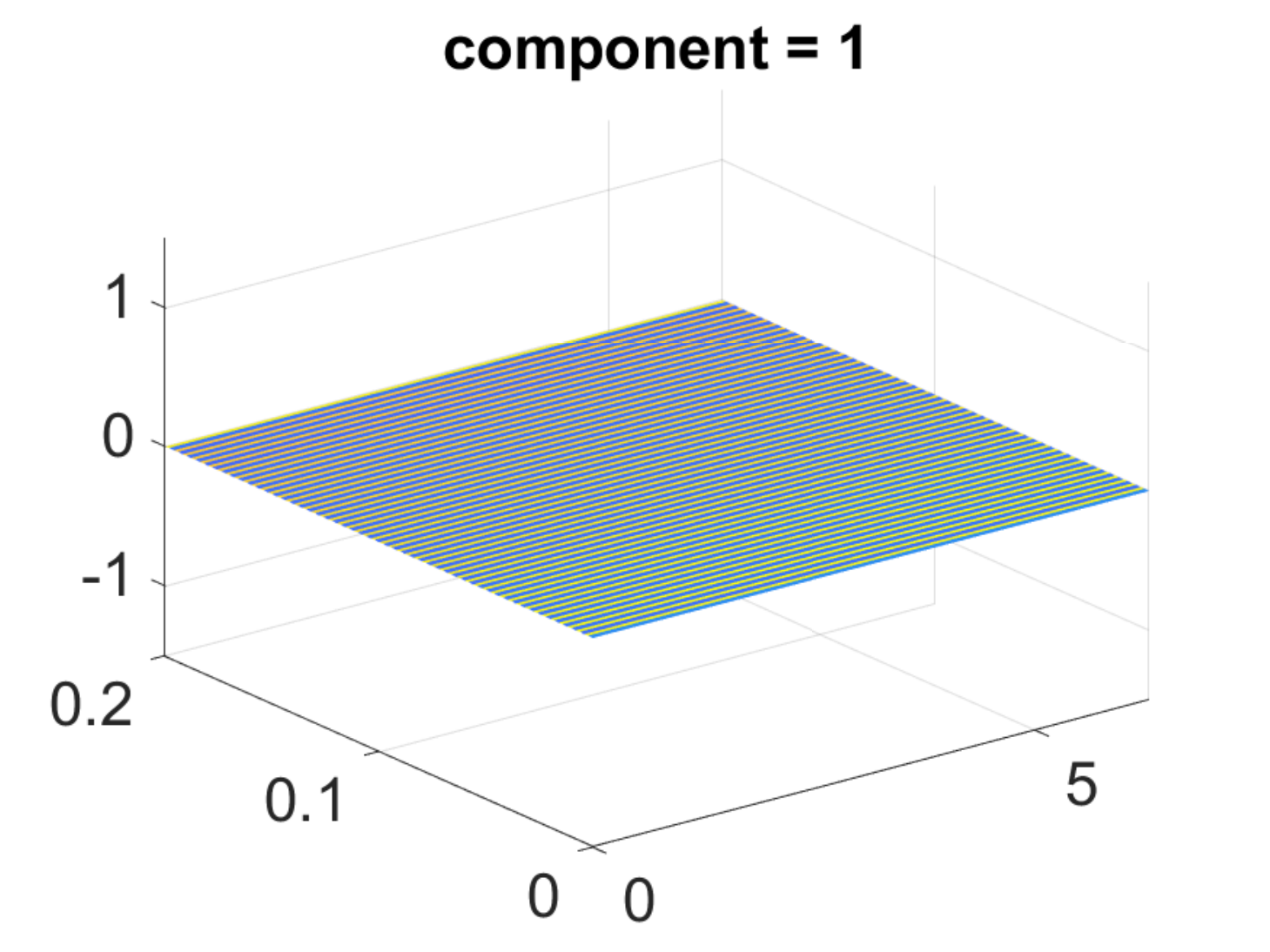}
\hspace{-0.5cm}
\includegraphics[width=4.3cm]{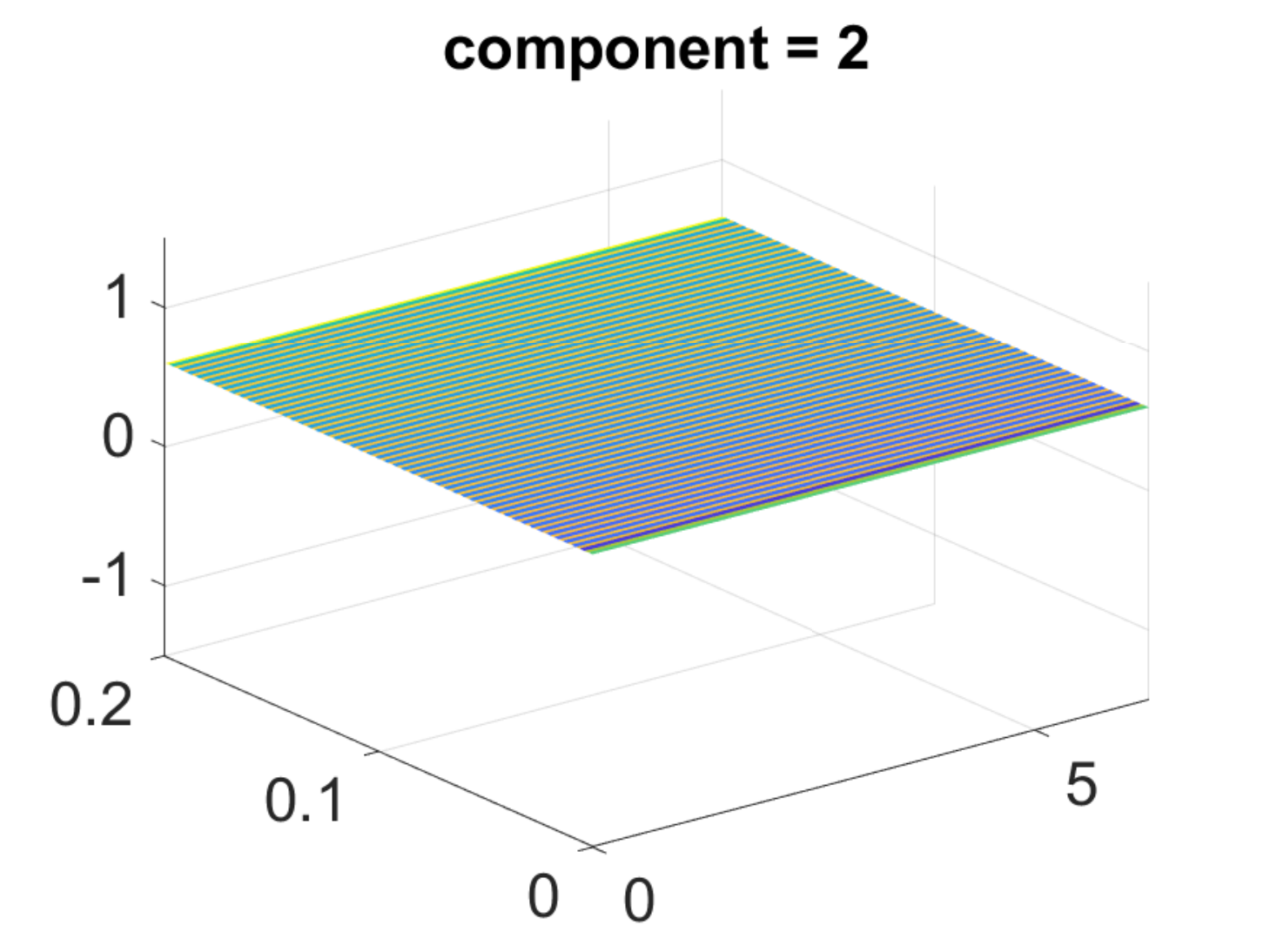}
\hspace{-0.5cm}
\includegraphics[width=4.3cm]{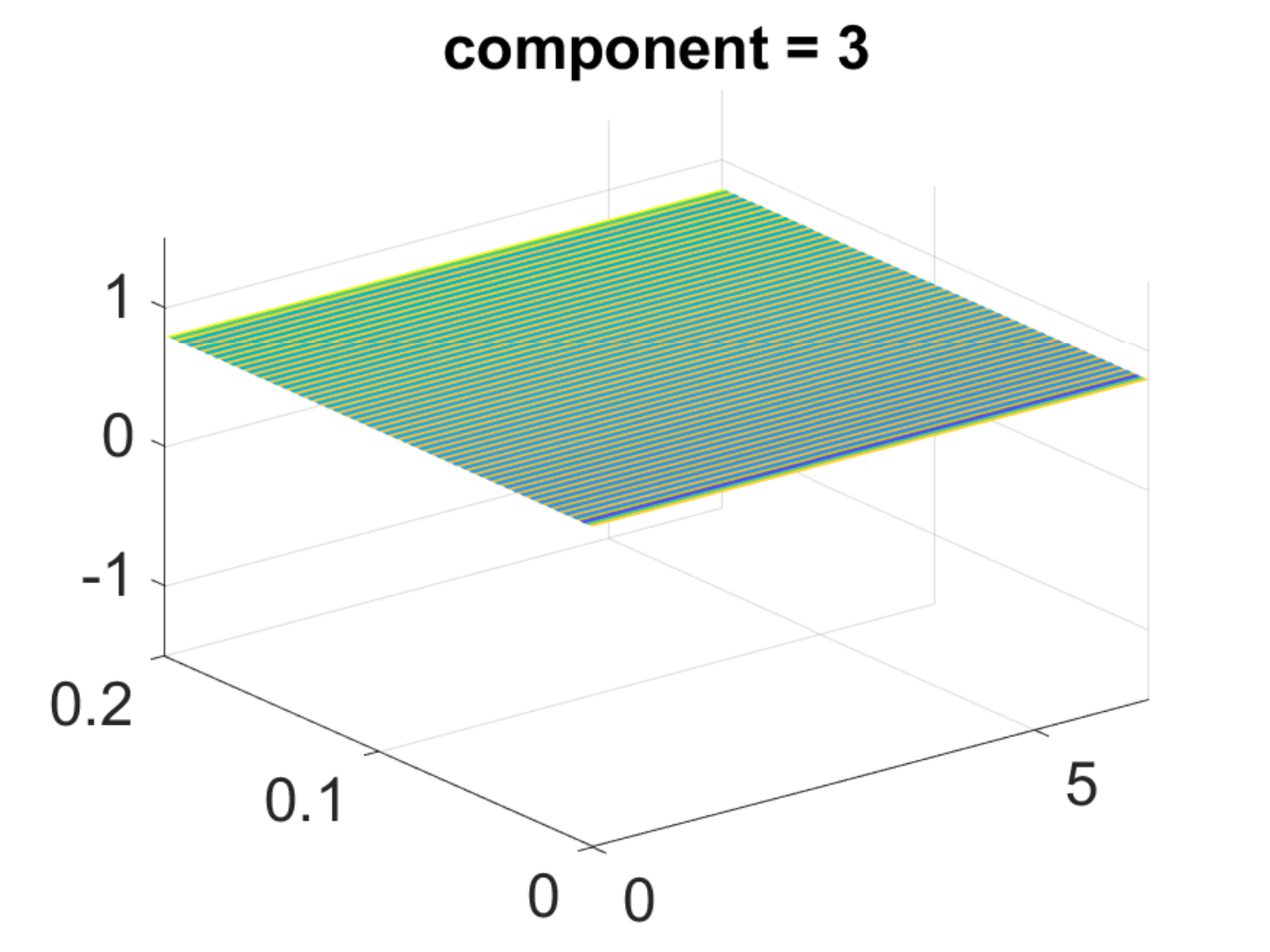}\\[1ex]
\includegraphics[width=4.3cm]{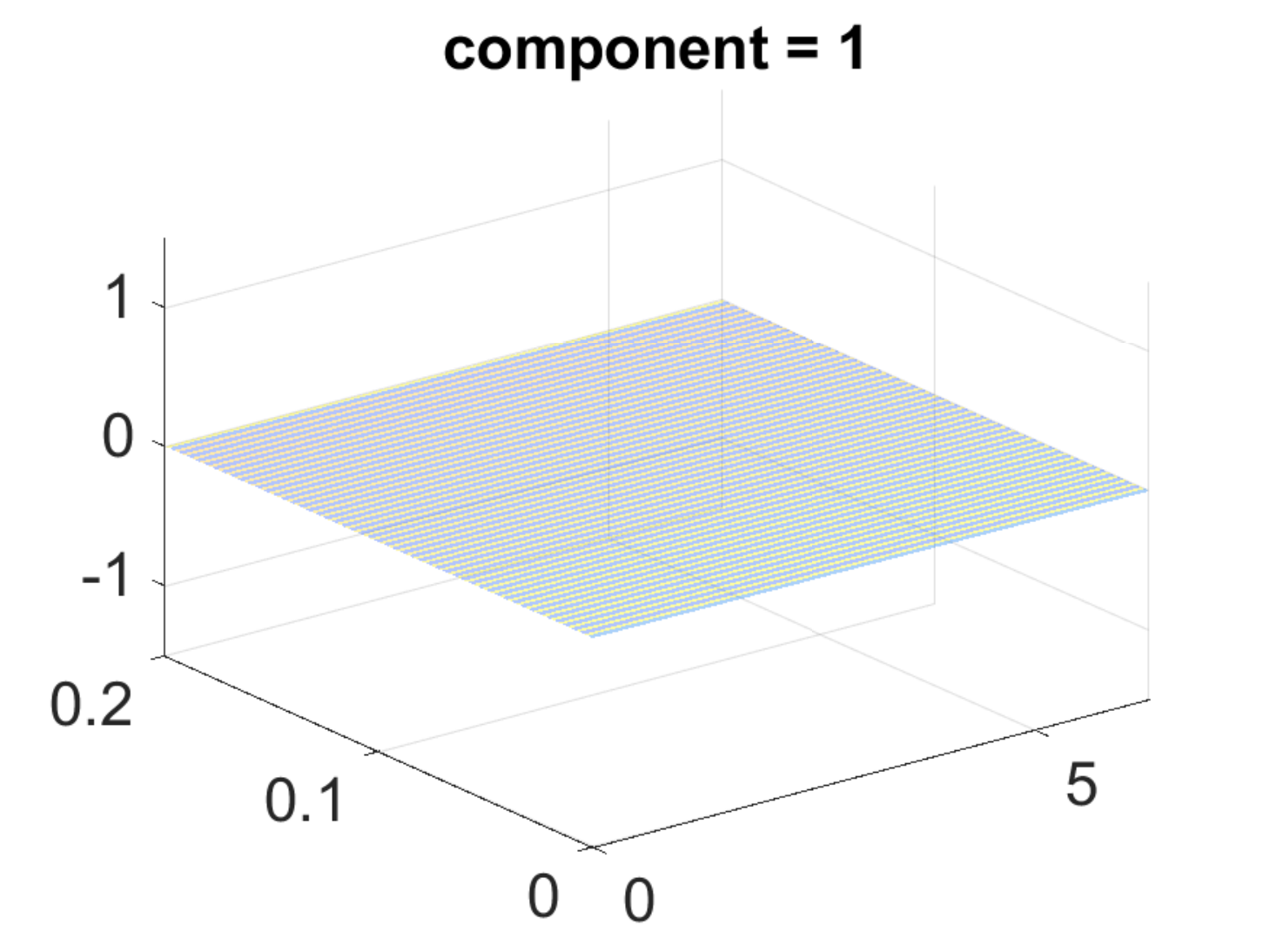}
\hspace{-0.5cm}
\includegraphics[width=4.3cm]{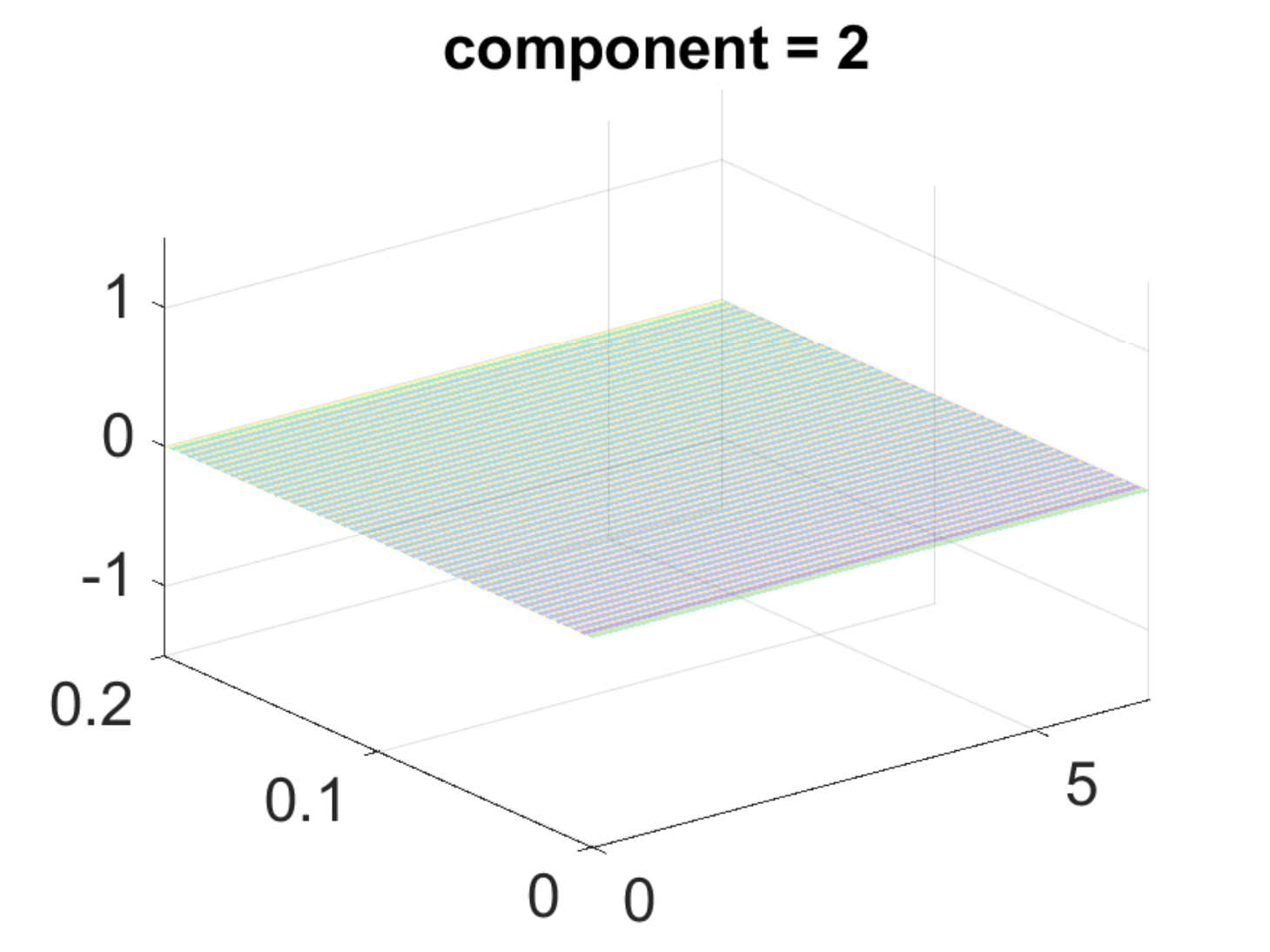}
\hspace{-0.5cm}
\includegraphics[width=4.3cm]{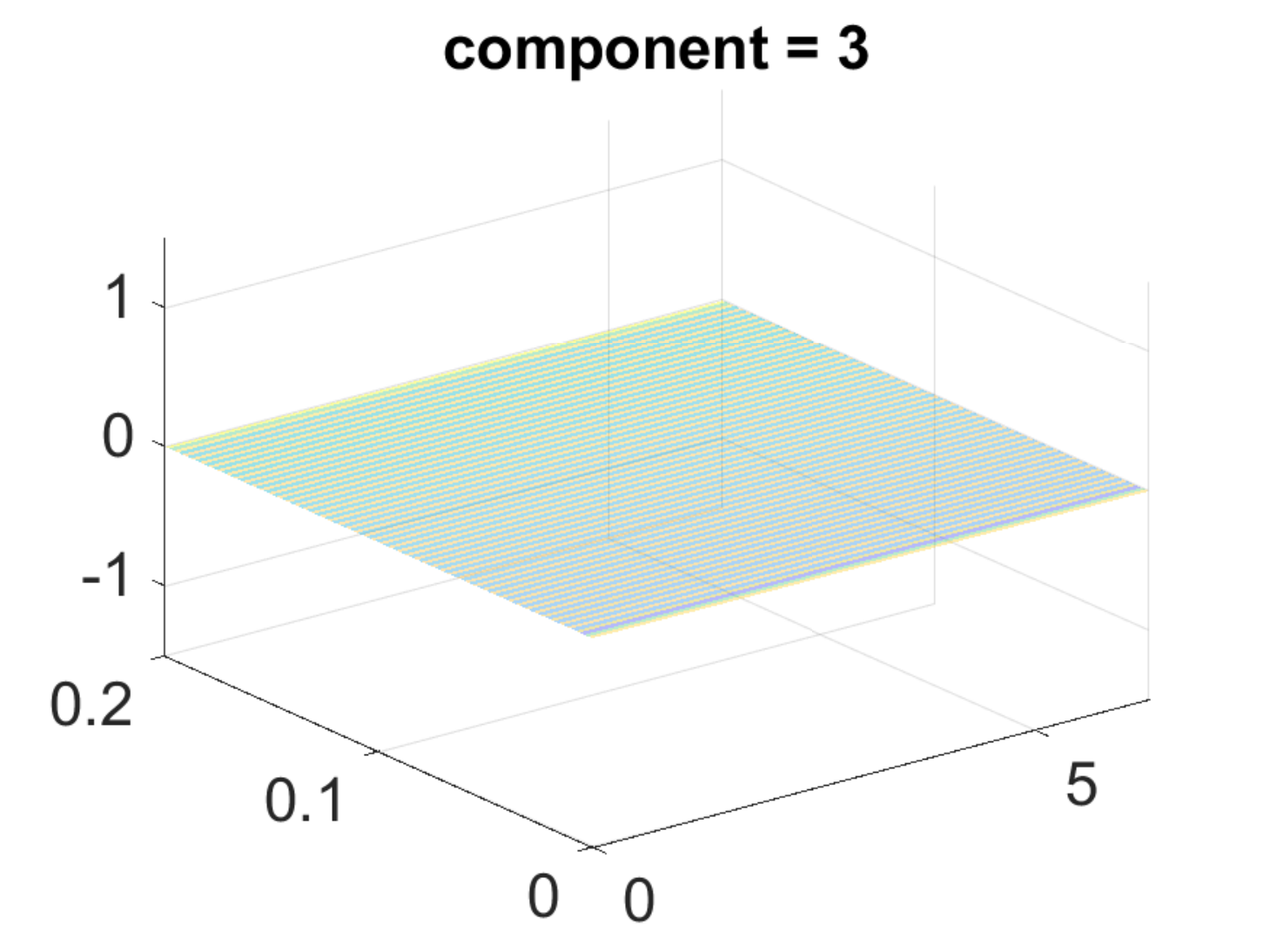}
\caption{Test 1. Reconstructed $\mv$ (top) and $\mv-\mv_{\text{exact}}$ (bottom) plotted against space (x-axis) and time (y-axis).}
\label{test1}

\vspace{3cm}
\includegraphics[width=6cm]{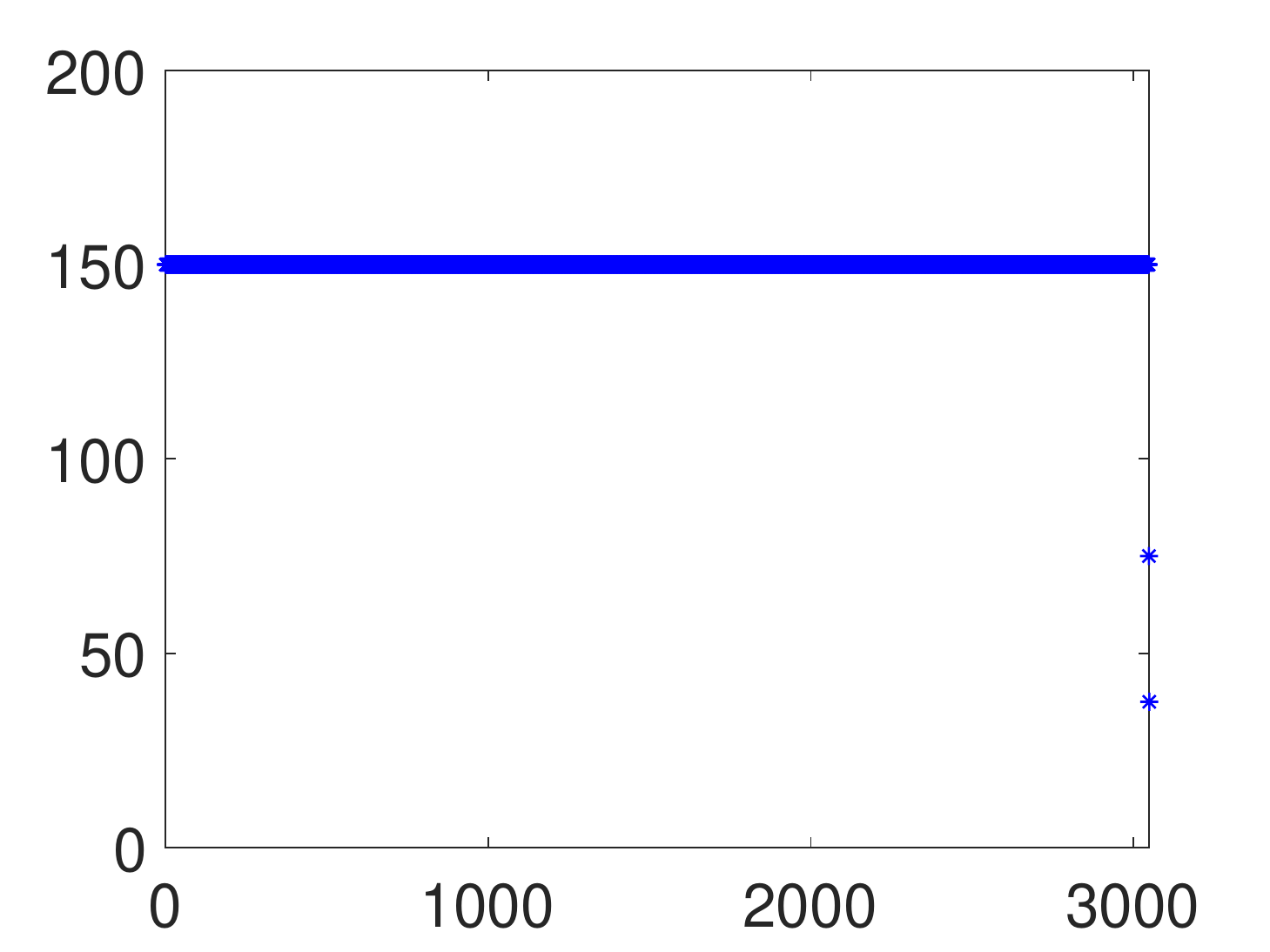}
\includegraphics[width=6cm]{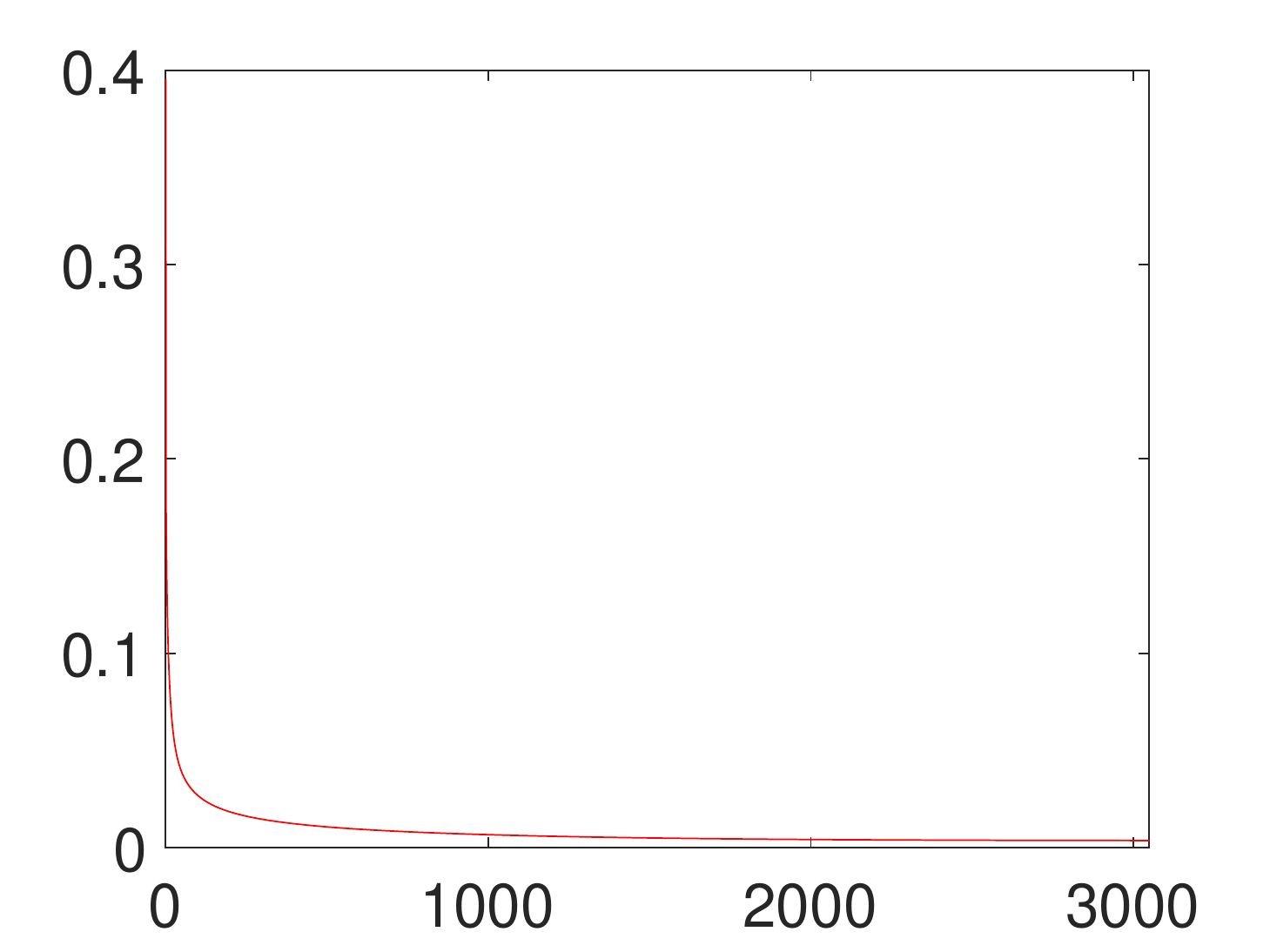}
\caption{Test 1. Plots of step size $\mu$ (left) and relative error $\frac{\|\m_k-\m_\text{ex}\|}{\|\m_\text{ex}\|}$ (right) over iteration index.}
\label{test11}
\end{figure}

\begin{figure}[p] 
\vspace{2cm}
\centering
\includegraphics[width=4.3cm]{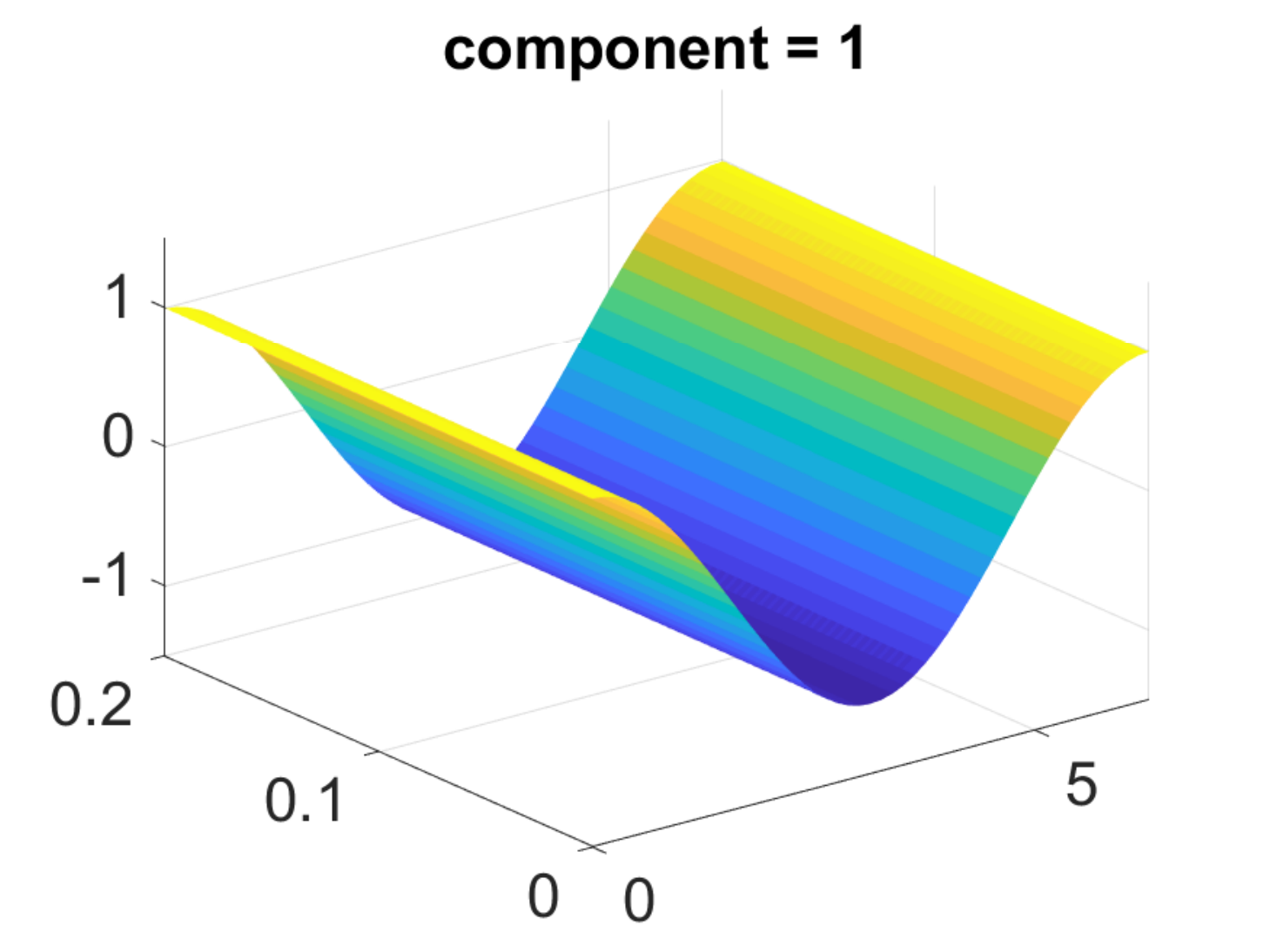}
\hspace{-0.5cm}
\includegraphics[width=4.3cm]{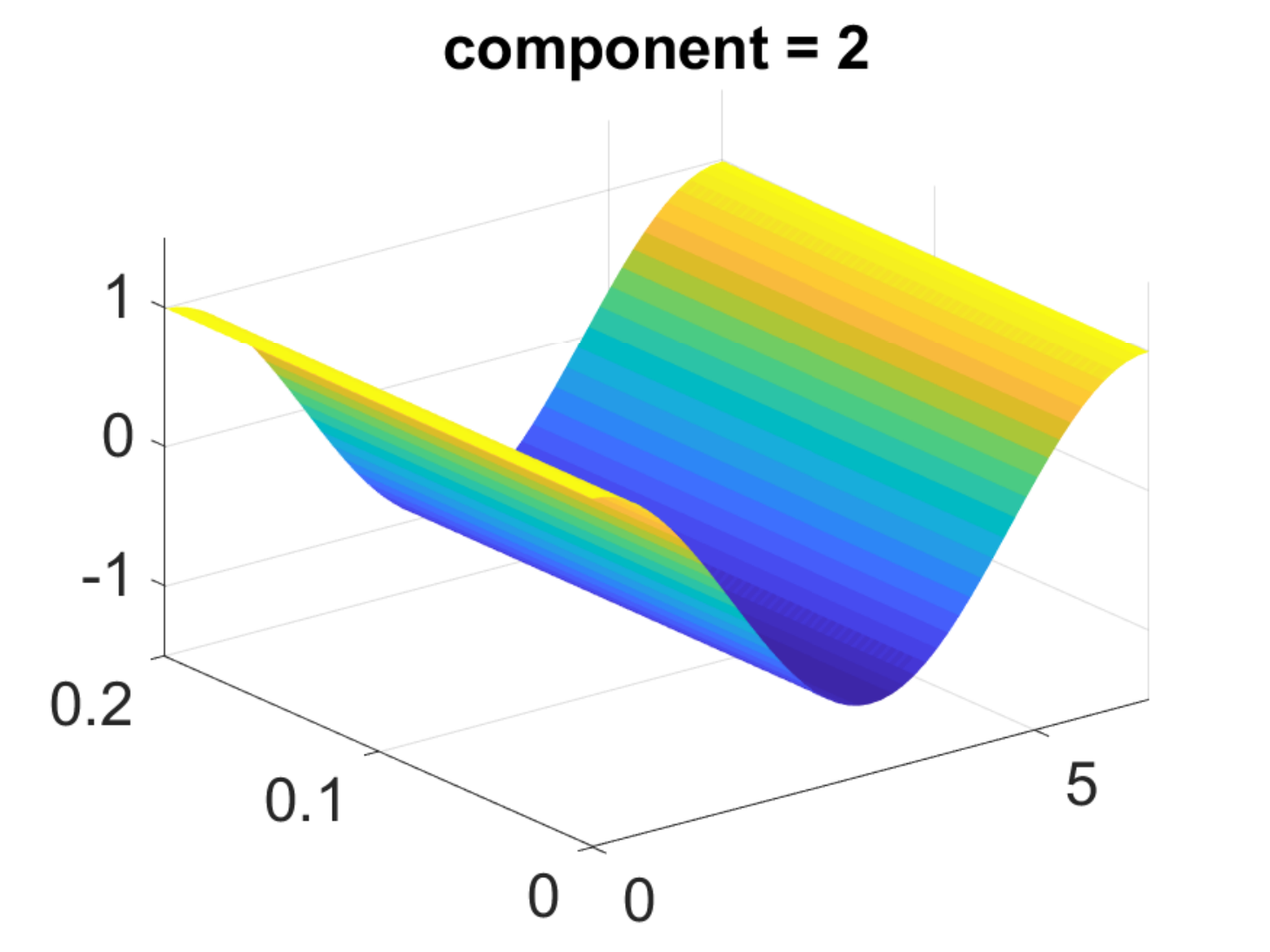}
\hspace{-0.5cm}
\includegraphics[width=4.3cm]{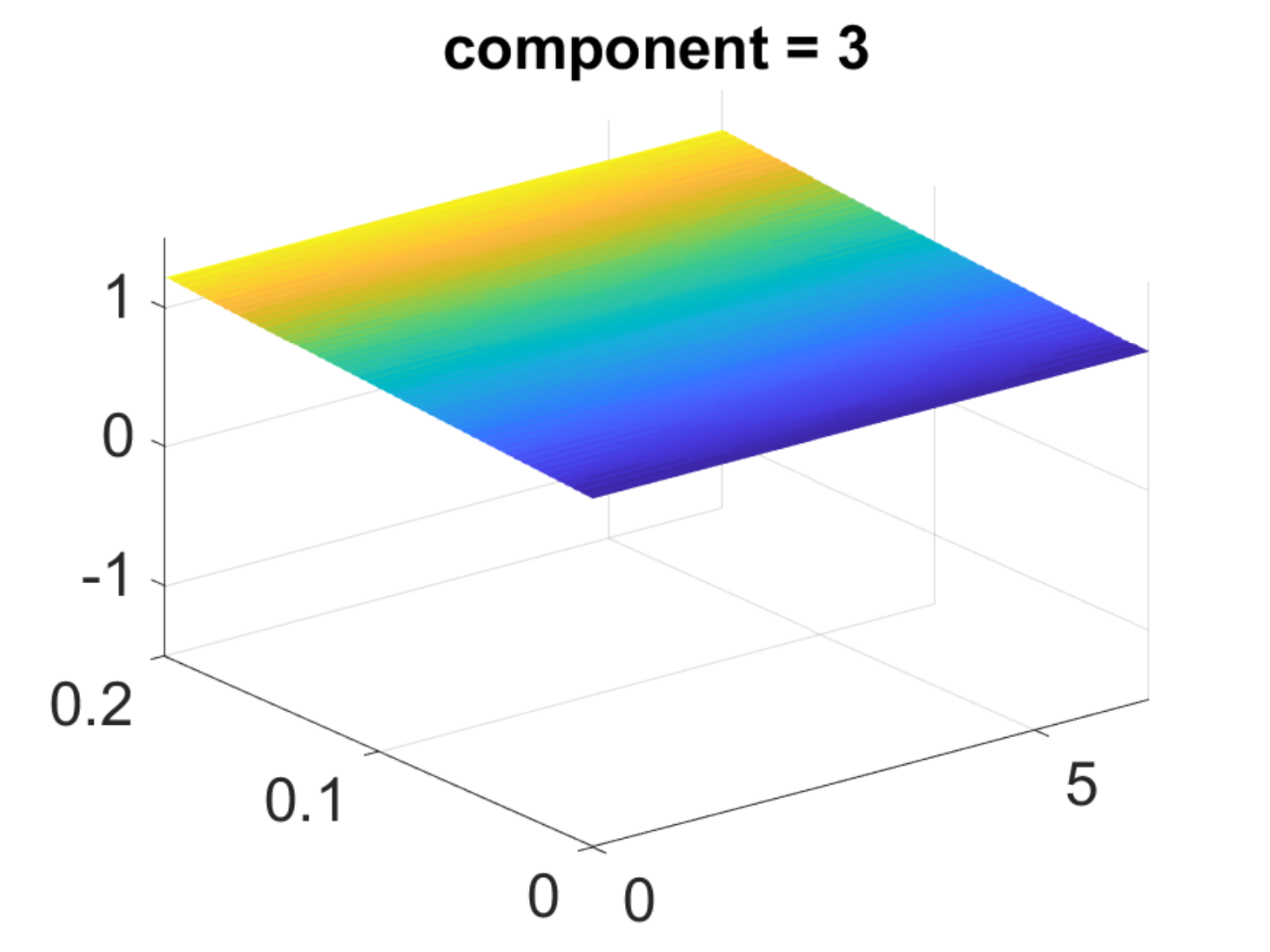}\\[1ex]
\includegraphics[width=4.3cm]{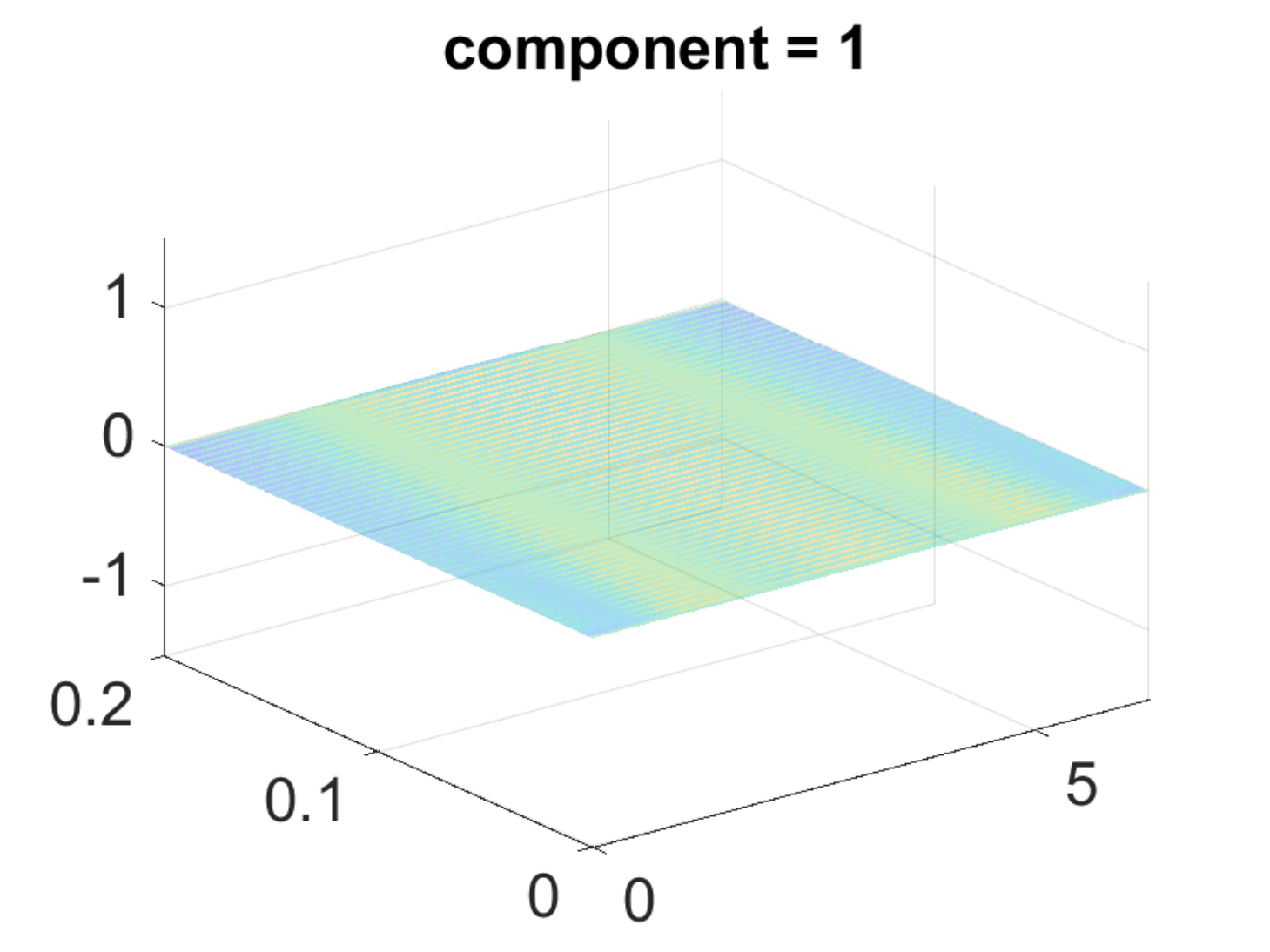}
\hspace{-0.5cm}
\includegraphics[width=4.3cm]{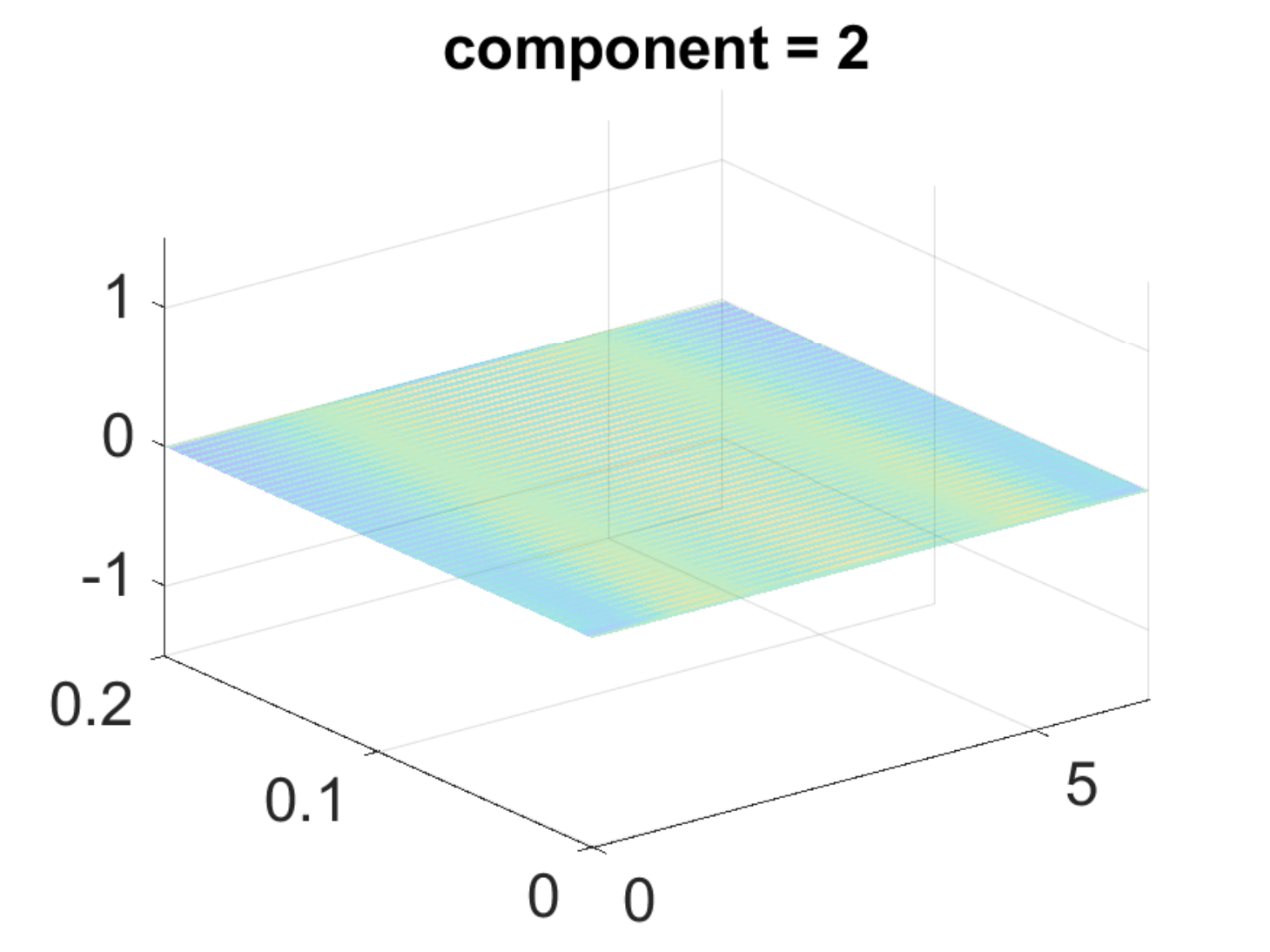}
\hspace{-0.5cm}
\includegraphics[width=4.3cm]{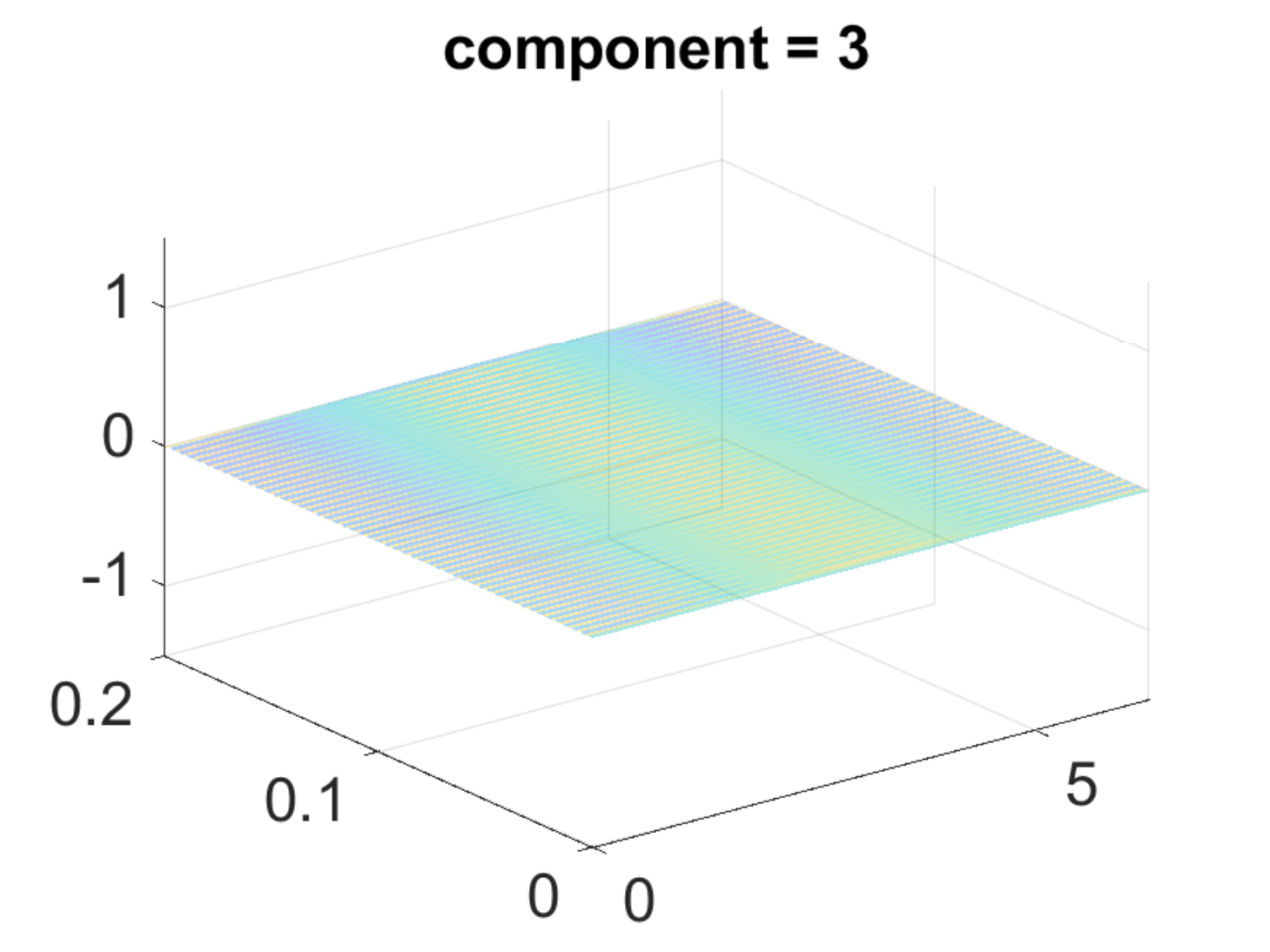}
\caption{Test 2. Reconstructed $\mv$ (top) and $\mv-\mv_{\text{exact}}$ (bottom). Left to right: each component.}
\label{test2}

\vspace{3cm}
\includegraphics[width=6cm]{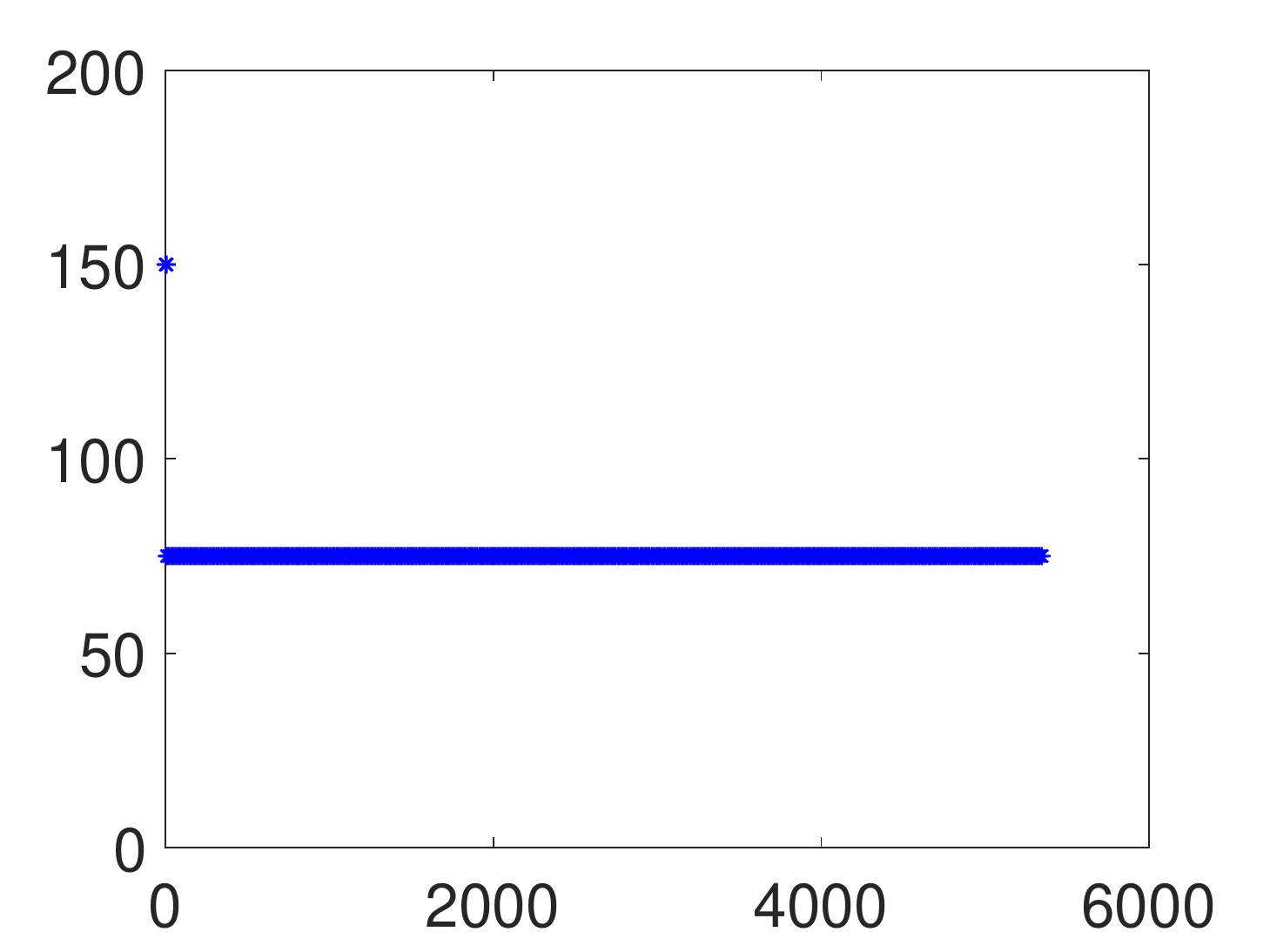}
\includegraphics[width=6cm]{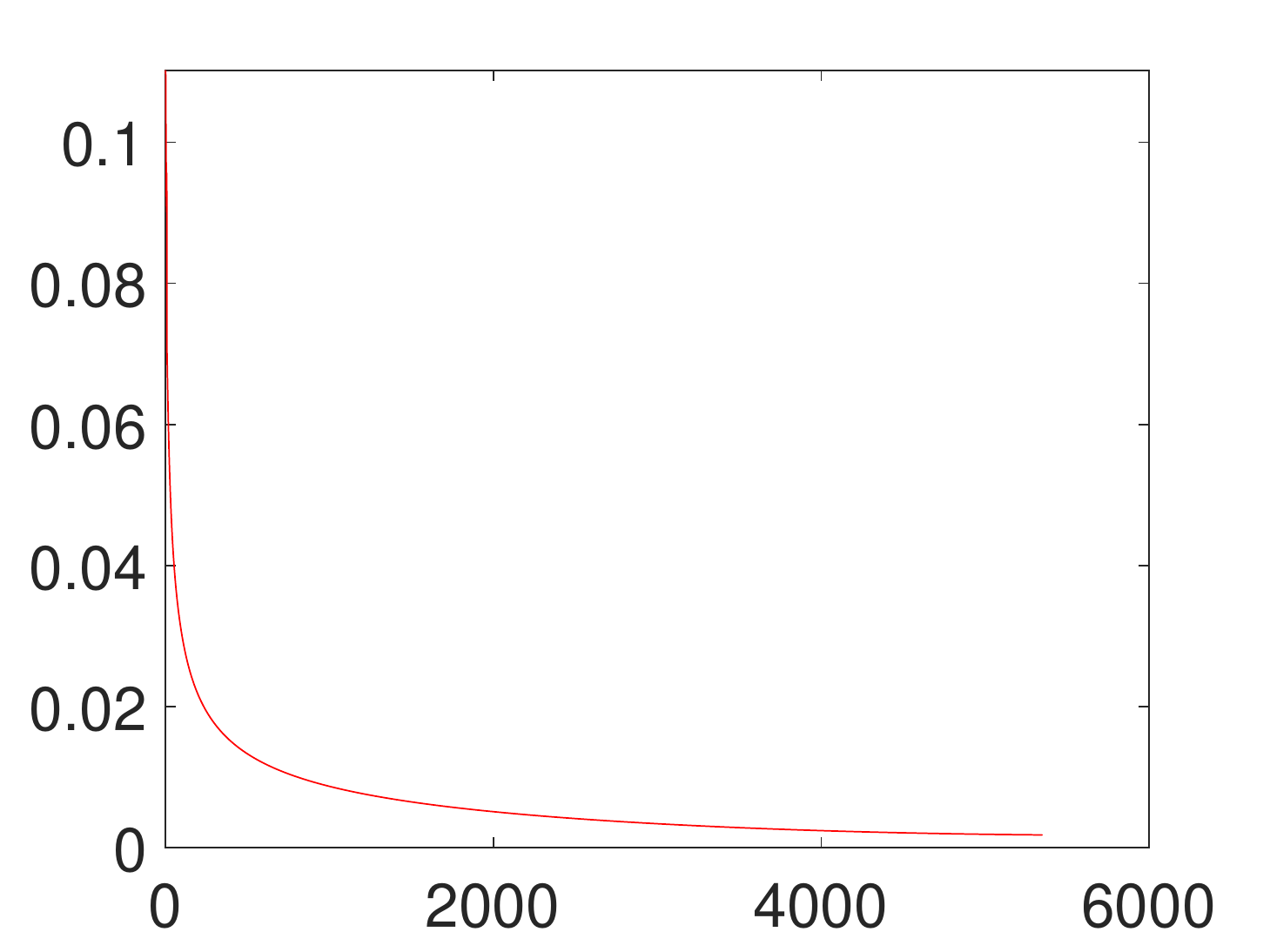}
\caption{Test 2. Plots of step size $\mu$ (left) and relative error $\frac{\|\m_k-\m_\text{ex}\|}{\|\m_\text{ex}\|}$ (right) over iteration index.}
\label{test22}
\end{figure}

\begin{figure}[p] 
\vspace{2cm}
\centering
\includegraphics[width=4.3cm]{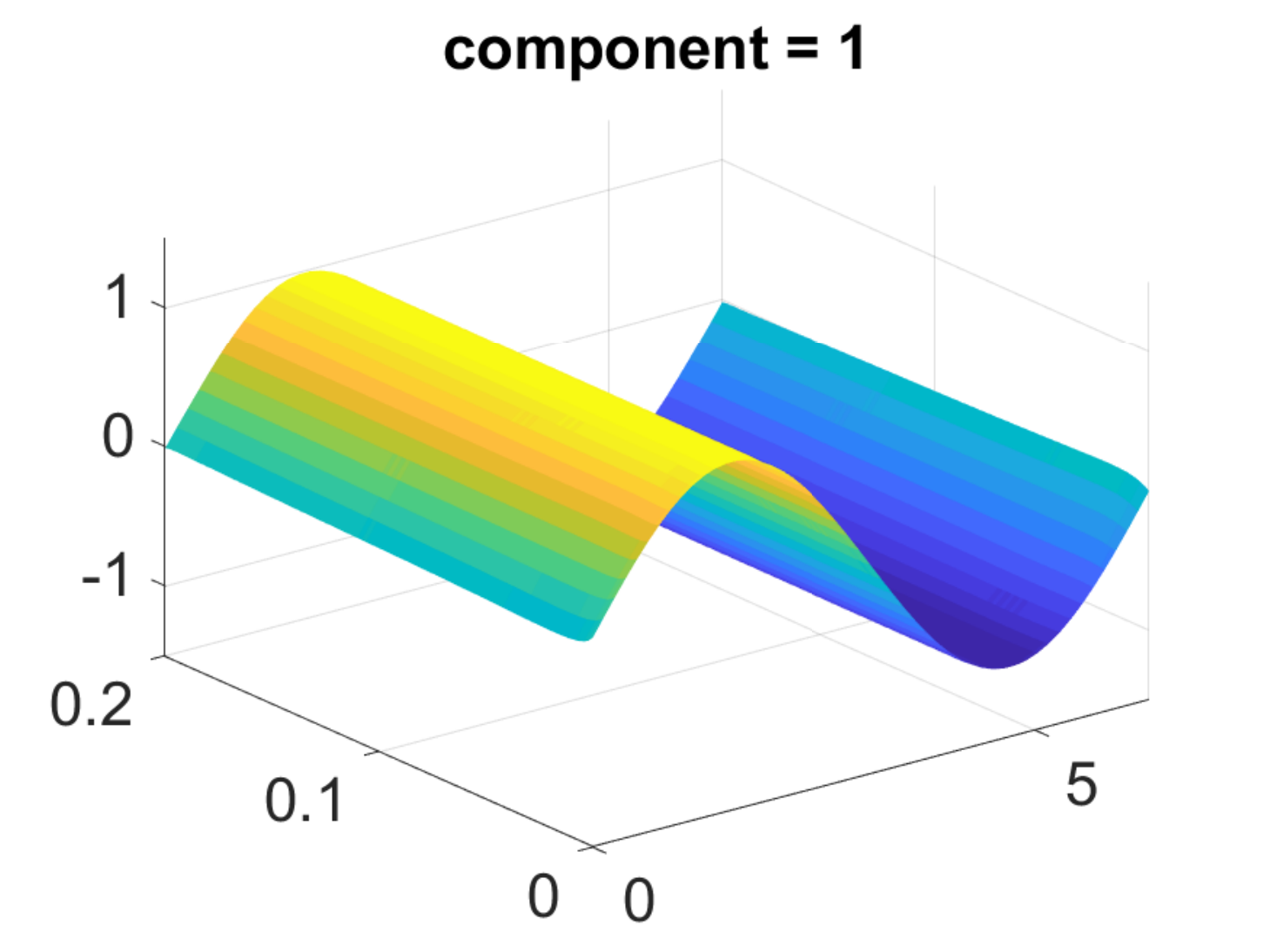}
\hspace{-0.5cm}
\includegraphics[width=4.3cm]{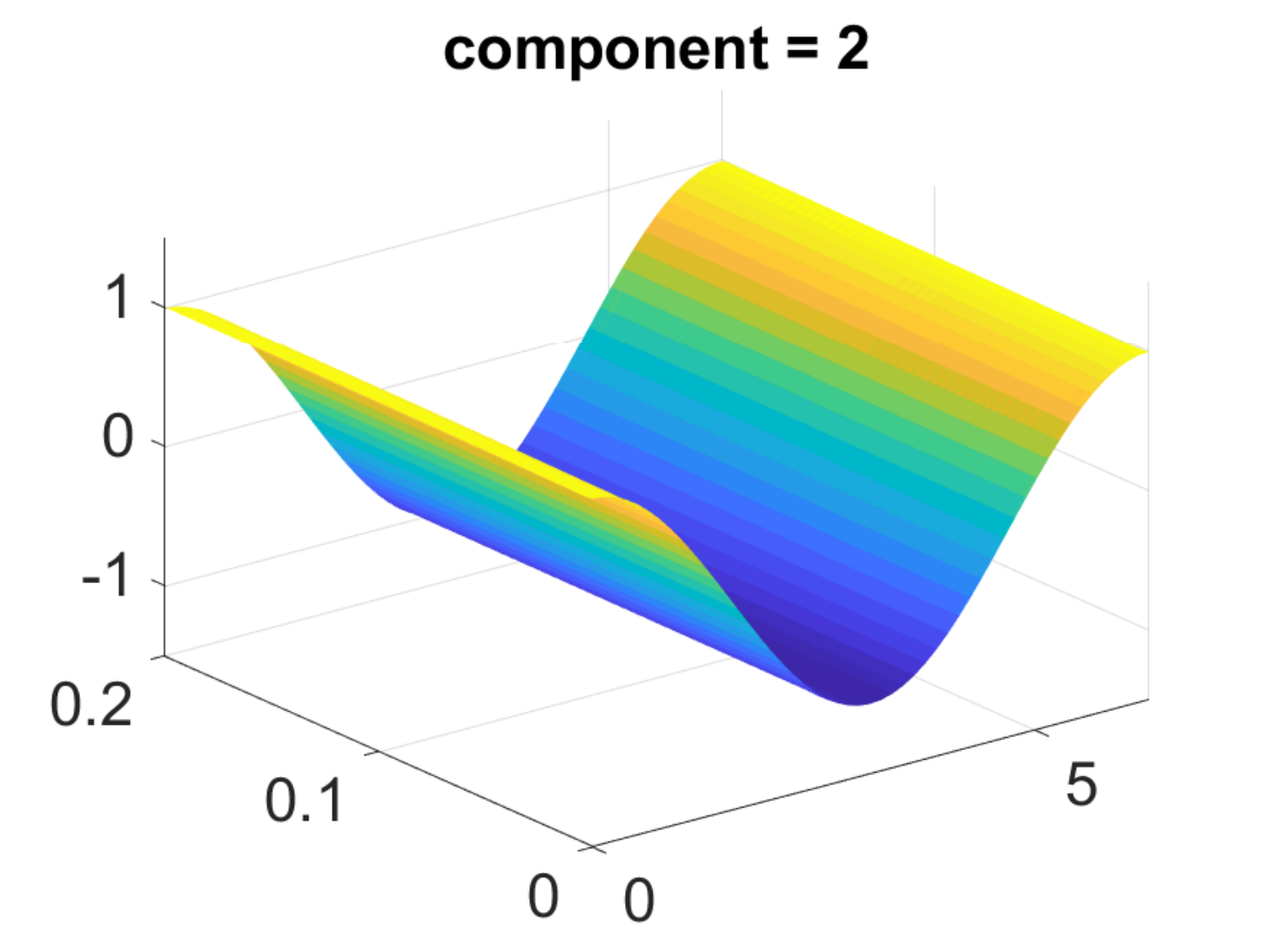}
\hspace{-0.5cm}
\includegraphics[width=4.3cm]{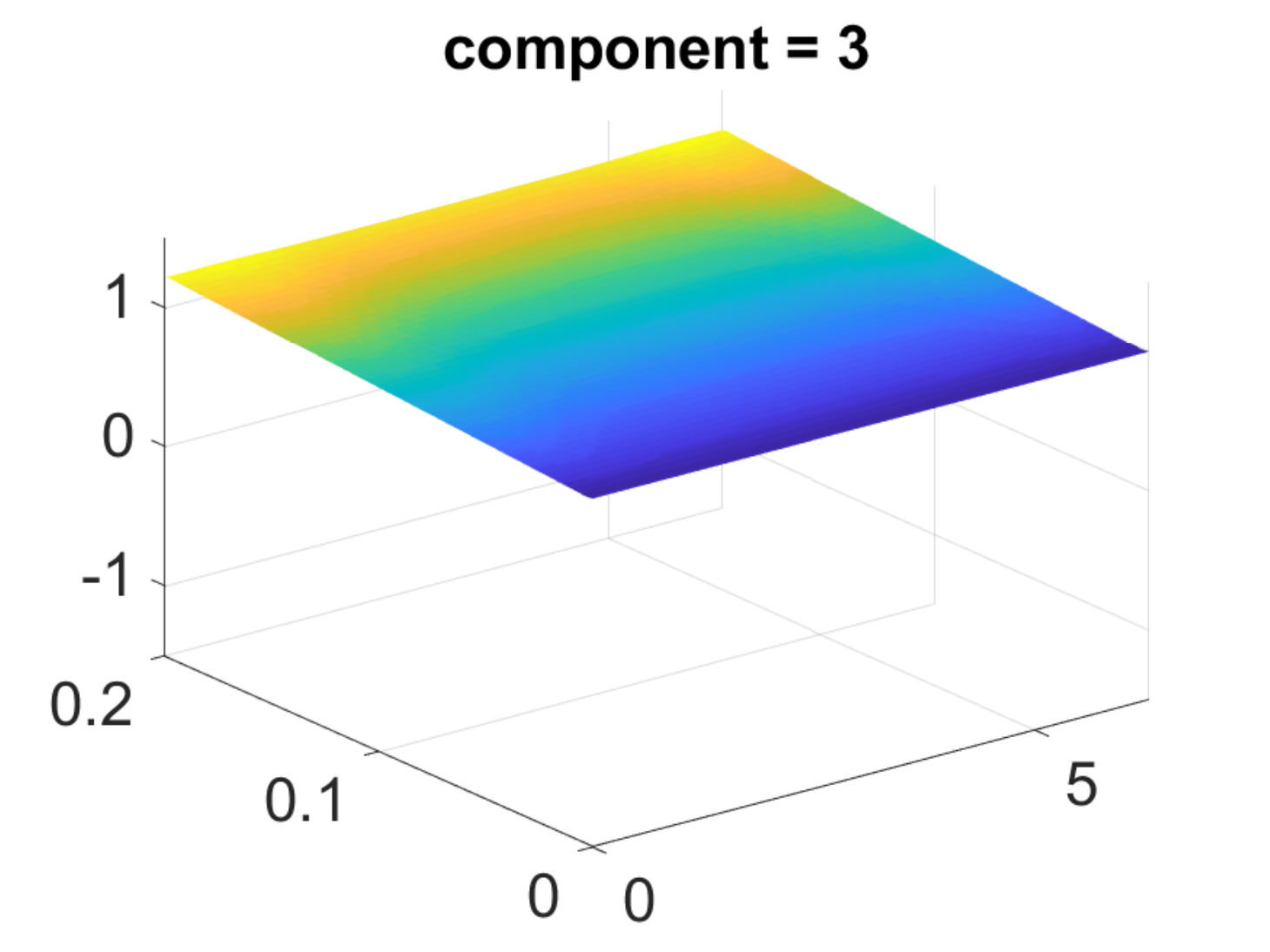}\\[1ex]
\includegraphics[width=4.3cm]{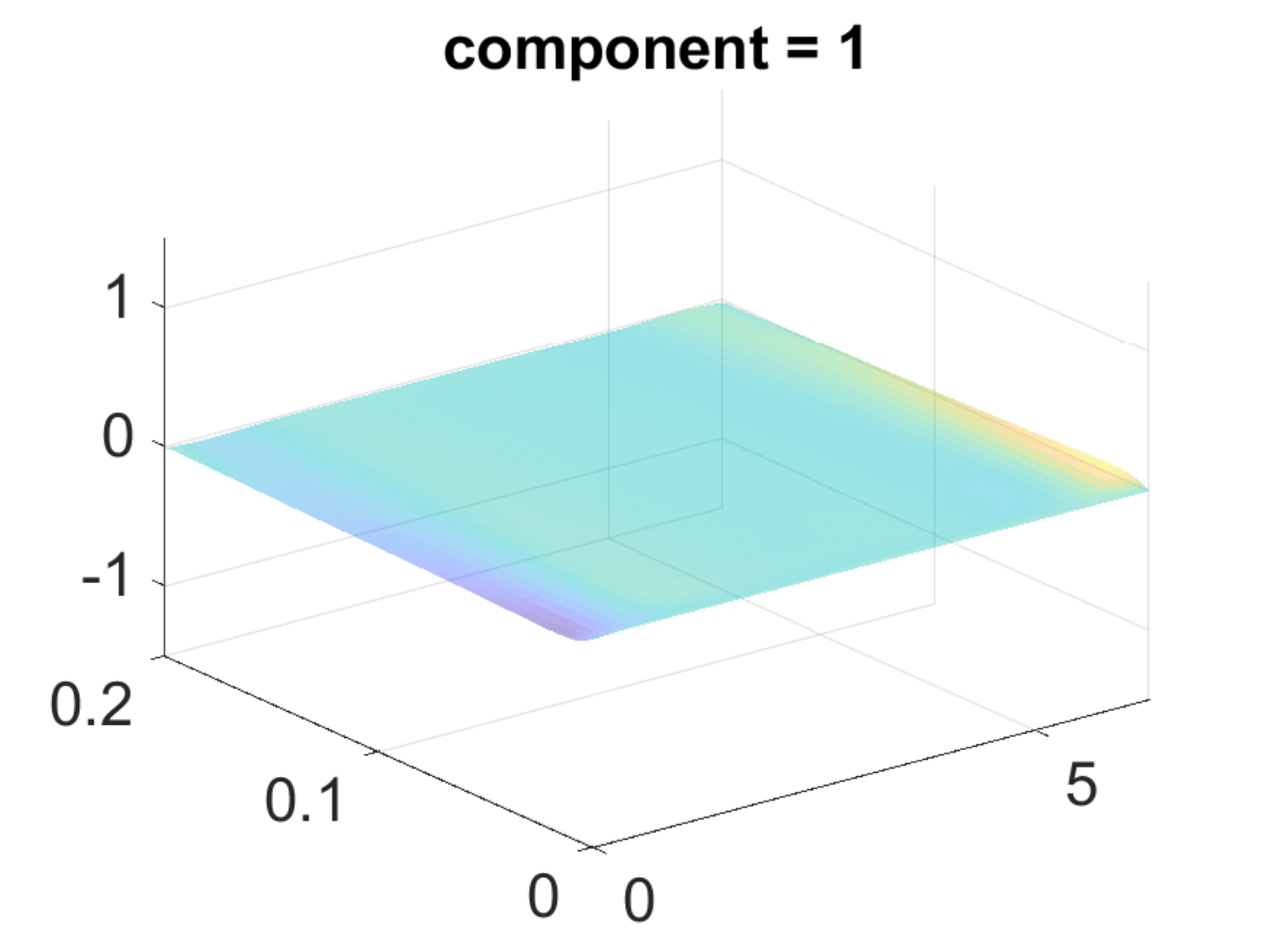}
\hspace{-0.5cm}
\includegraphics[width=4.3cm]{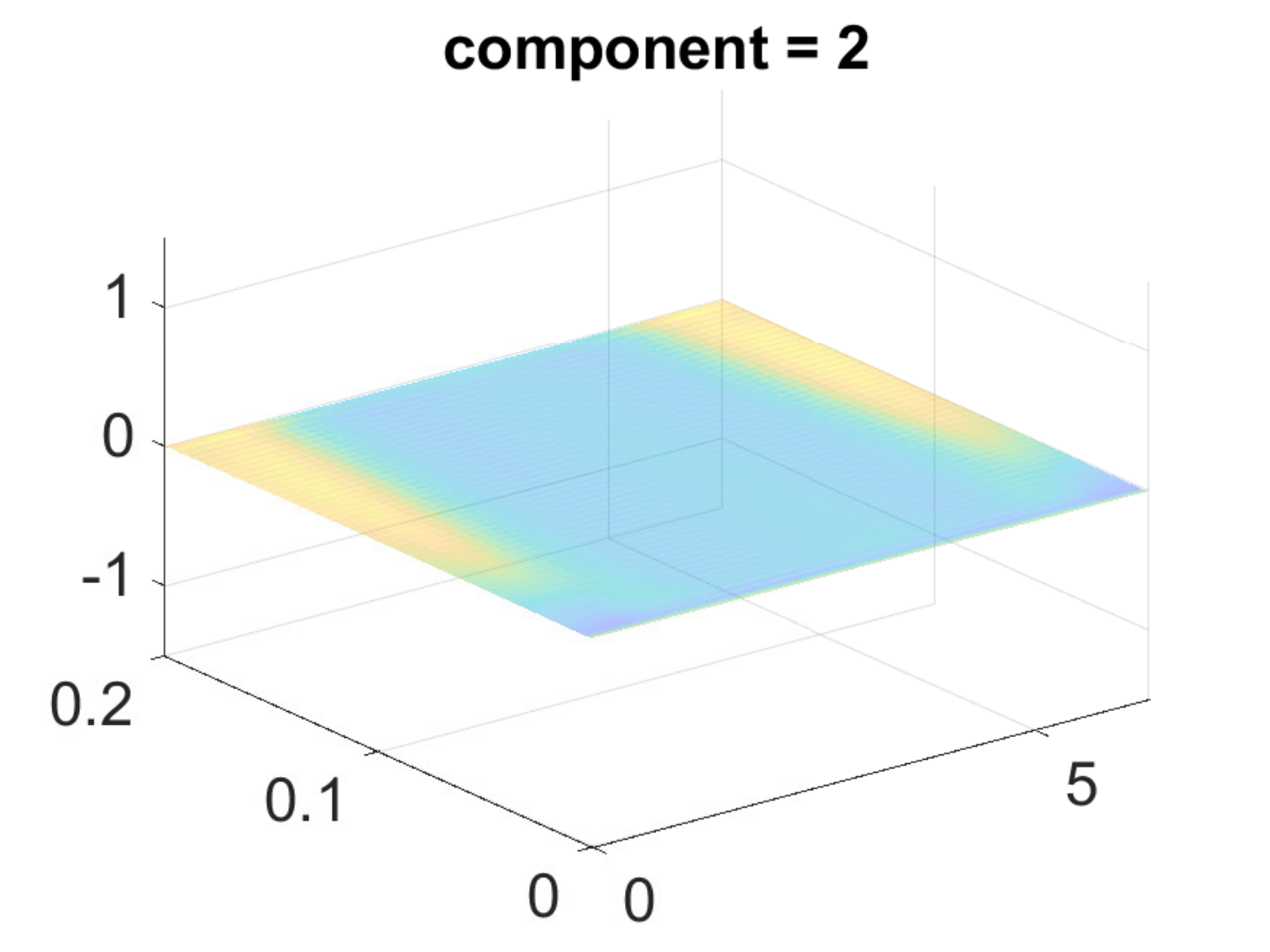}
\hspace{-0.5cm}
\includegraphics[width=4.3cm]{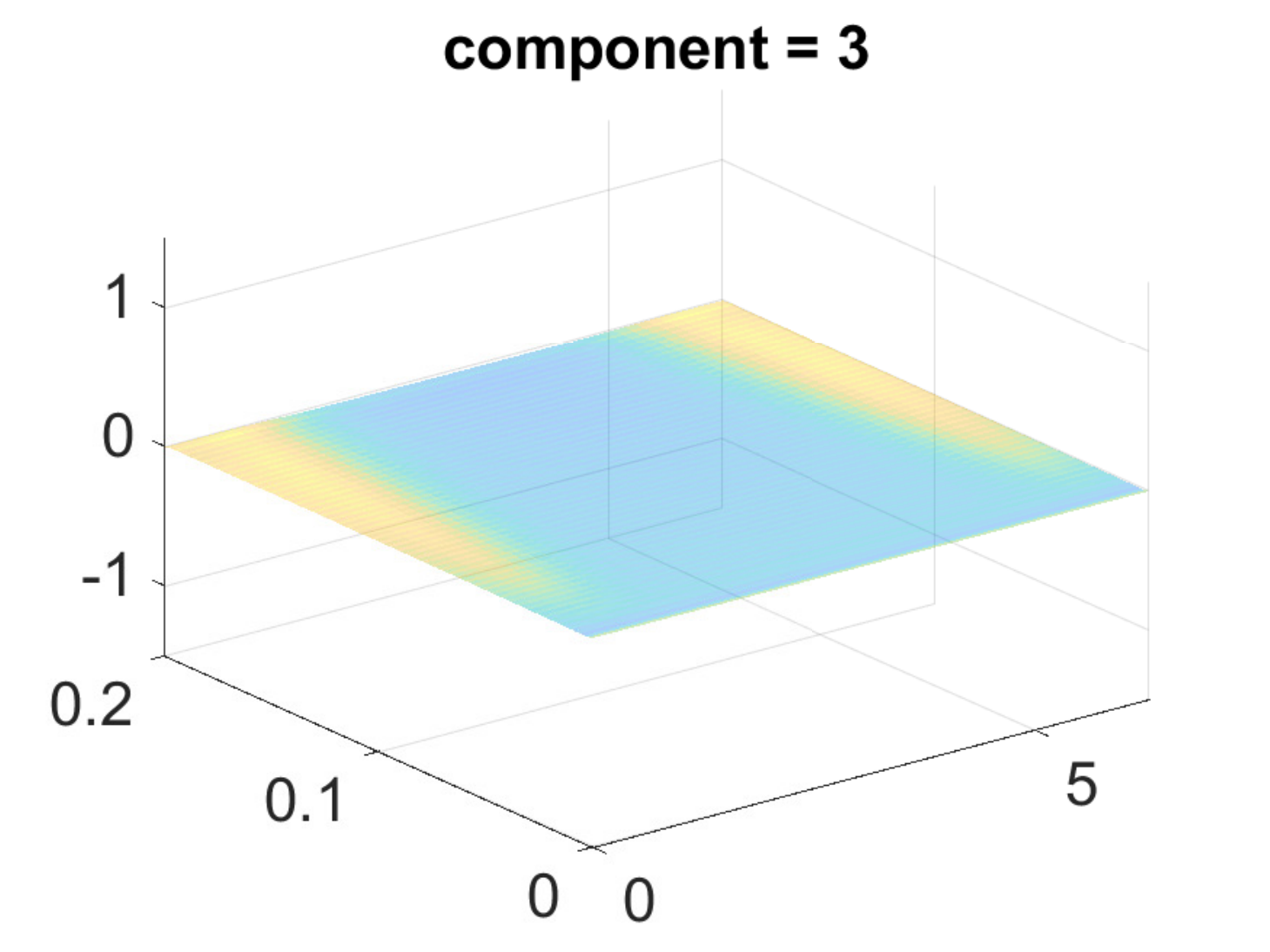}
\caption{Test 3. Reconstructed $\mv$ (top) and $\mv-\mv_{\text{exact}}$ (bottom) plotted against space (x-axis) and time (y-axis).}
\label{test3}

\vspace{3cm}
\includegraphics[width=6cm]{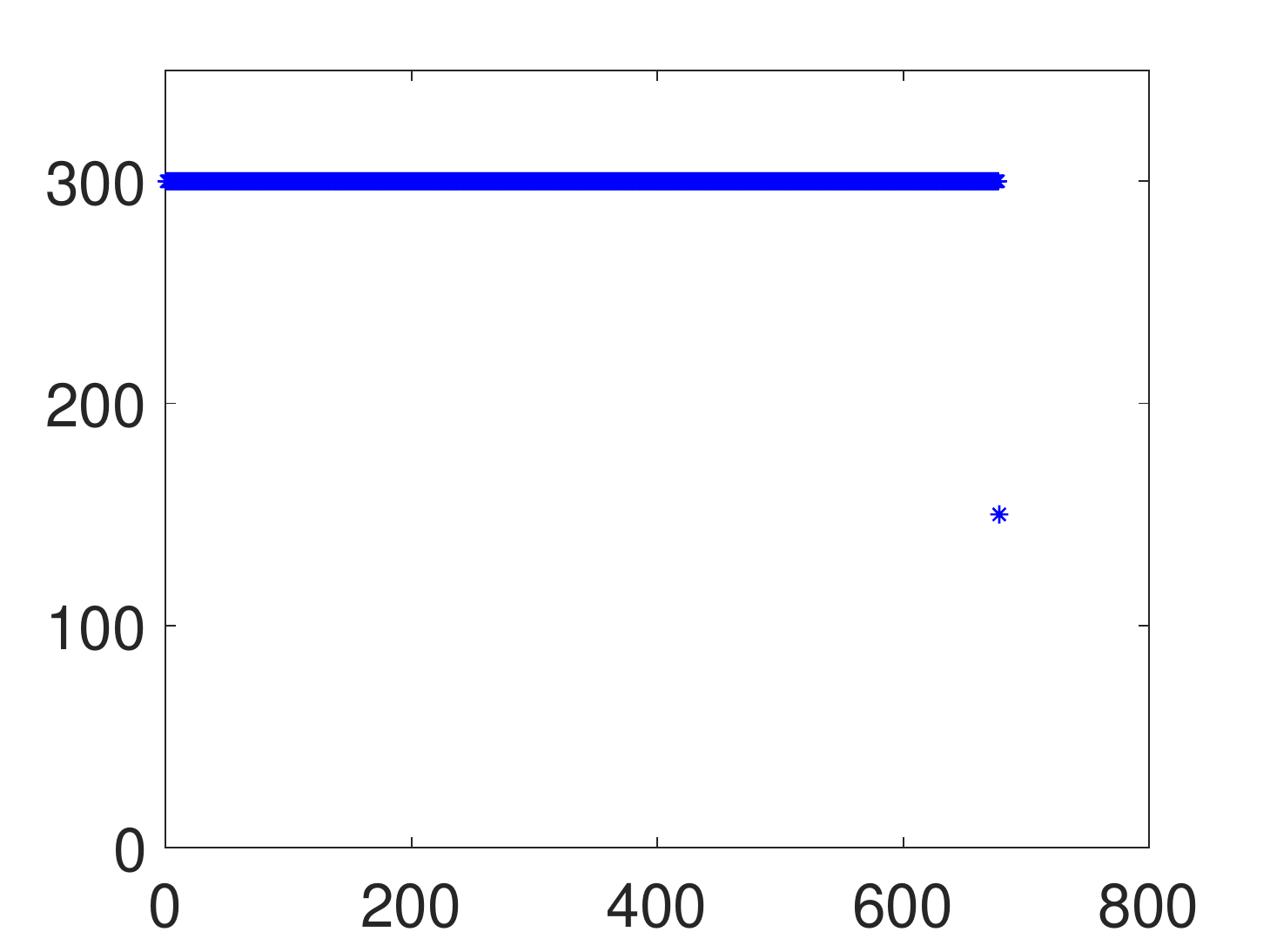}
\includegraphics[width=6cm]{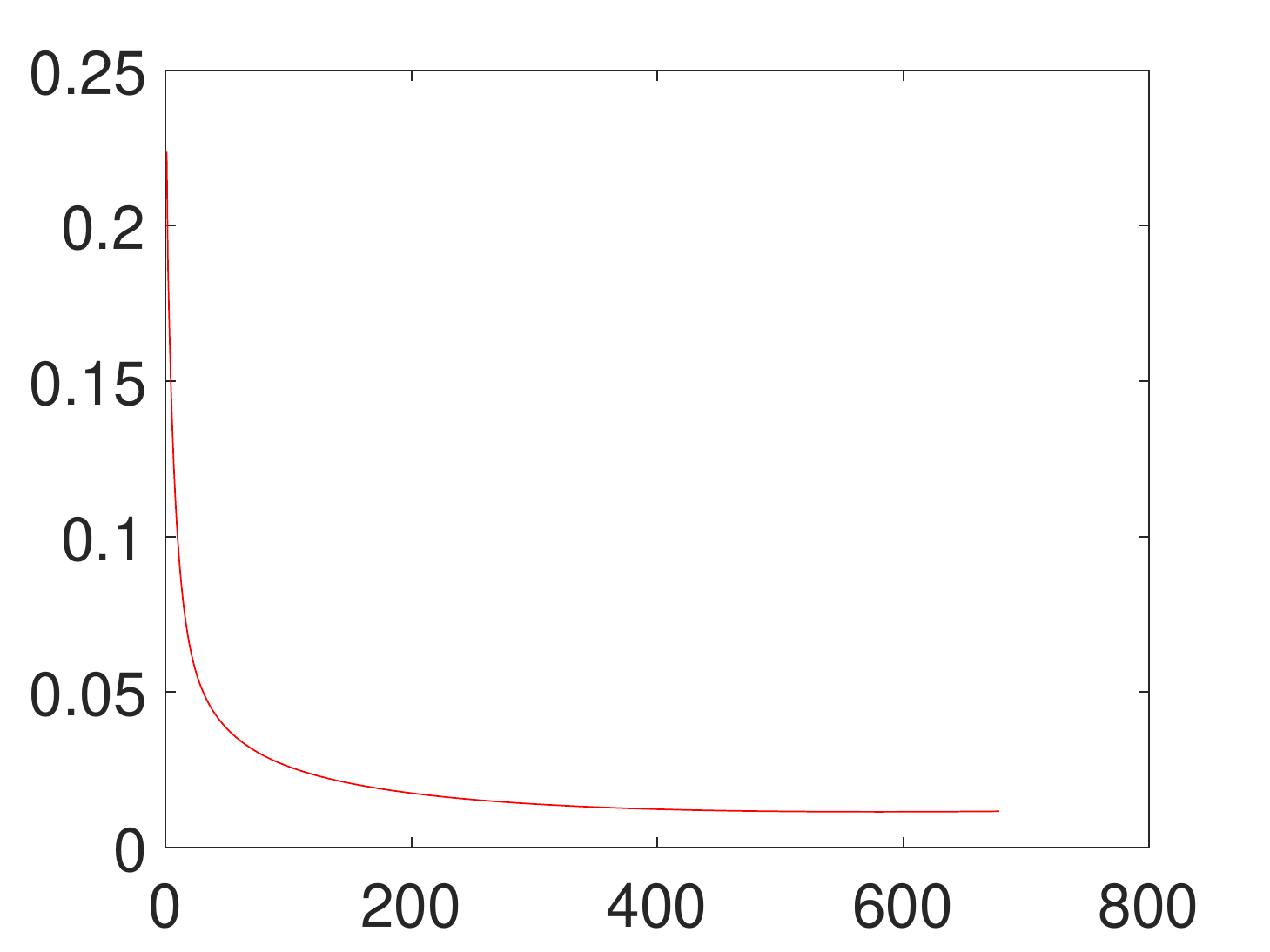}
\caption{Test 3. Plots of step size $\mu$ (left) and relative error $\frac{\|\m_k-\m_\text{ex}\|}{\|\m_\text{ex}\|}$ (right) over iteration index.}
\label{test33}
\end{figure}

Figures \ref{test1}, \ref{test2} and \ref{test3} display the results in time and space of the reconstructed states together with comparisons to the true ones. In Figure \ref{test11}, the initial error in the $L^2(0,T,L^2(\Om2R))$-norm measuring the distance between $\mv$ and $\mv_{\text{exact}}$ declines from $40\%$ to $0.4\%$ after 3050 iterations.
Now consider the test cases 2 and 3 presented in Figures \ref{test22} and \ref{test33}. Starting from $11\%$ and $22\%$, the errors drop to $0.2\%$ and $1.1\%$ after 5350 and 1850 iterations, respectively. In those tests, we create initial guesses for $\mh$, thus $\mv$, by perturbing the exact ones by different amounts to closely inspect the convergence of the method. In practice, one can choose the initial guess for the LLG solver as $\m_0$, which means just $(0,0,0)$ in Test 1 and $(0,0,1)$ in the latter two.

Relying on the monotonicity of the residual sequence, we implement an adaptive Landweber step size scheme in order to search for an appropriate one (Figures \ref{test11}, \ref{test22}, \ref{test33}, left). In particular, in each iteration, a residual comparison with the previous step takes place. If the current residual shows a decrease, the current step size $\mu$ is accepted, otherwise it is bisected. The iterations are terminated after reaching a certain level in smallness of the step size, alternatively speaking, the residual is not able to get significantly smaller. One can stop the iterations earlier by checking the smallness of the step size together with the residual tolerance. The runtime reports: 289 seconds, 512 seconds and 205 seconds respectively for the three tests.

\begin{table}[!htb] 
\bigskip\bigskip
\caption{Test cases and run parameters}
\bigskip\bigskip
\centering
\begin{tabular}{ |c|c|c|c| }
\hline
\rule{0pt}{15pt}
Test & 1 & 2 & 3 \\[1ex]
\hline
\rule{0pt}{13pt}
$\altil_1$ & 1 & 2 & 1\\
\rule{0pt}{15pt}
$\altil_2$ & -1 & 0 & 0\\
\rule{0pt}{13pt}
$\hv$  & $\cfrac{2}{5}(0,3,4)$ & $-(\cos(x),\cos(x),0)$ & (0,0,0)\\
\rule{0pt}{13pt}
$\mv_{\text{exact}}$ & $\cfrac{1}{5}(0,3,4)$ & $(\cos(x),\cos(x),e^t)$ & $(\sin(x),\cos(x),e^t)$\\
\rule{0pt}{13pt}
$\m_0$ & $\cfrac{1}{5}(0,3,4)$ & $(\cos(x),\cos(x),1)$ & $(\sin(x),\cos(x),1)$\\
\rule{0pt}{13pt}
$\hat{\mv}_{exact}$ & (0,0,0) & (0,0,$e^t-1$) & (0,0,$e^t-1$)\\
\rule{0pt}{13pt}
Initial guess $\mh$ & $-5t(1,1,1)$ & $-5t\cos(x)(1,1,1)$ & $-\cfrac{\sin(30t)}{5}(1,1,1)$\\
\rule{0pt}{13pt}
Step size $\mu$ & 150 & 75 & 300\\
\rule{0pt}{13pt}
\# iterations & 3050 &5350 & 680\\
\hline
\end{tabular}
\label{tab-1}
\end{table}

\clearpage

\clearpage

\clearpage

\subsection{LLG solver with physical parameters}
\label{sec:LLGsyn}
In this section, we perform some simulations to illustrate the behaviour of the magnetization vector as a response to an external field using physically relevant parameters. Table \ref{tab-physic} gives an overview on common parameters, that can be found, e.g., in \cite{BanasPagePraetorius, tk18}. The length of the magnetization vector $\m$ is specified by $m_\text{S}=474 000$ \blue{J $\text{m}^{-3} \text{T}^{-1}$};  $\altil_1$ and $\altil_2$ differ by a factor of $m_\text{S}$. For better numerical computing, we shall do the following scaling for the LLG equation: \\
Let $\widetilde{\m}:=\frac{\m}{m_\text{S}}, \widetilde{\alpha}_1:=m_\text{S}\widetilde{\alpha}_1$, i.e.,
\begin{displaymath}
 \widetilde{\alpha}_1 := \frac{\gamma \alpha_{\mathrm{D}}}{1+\alpha_{\mathrm{D}}^2} > 0, \qquad \widetilde{\alpha}_2 := \frac{\gamma}{1+\alpha_{\mathrm{D}}^2} > 0
\end{displaymath}
and let
\begin{eqnarray*}
\altil_1=\frac{\widetilde{\alpha}_1}{\widetilde{\alpha}_1^2+\widetilde{\alpha}_2^2}, \qquad \altil_2=\frac{\widetilde{\alpha}_2}{\widetilde{\alpha}_1^2+\widetilde{\alpha}_2^2}, \qquad\h=\mu_0 m_\text{S}\textbf{H}_\text{ext},
\end{eqnarray*}
we obtain the LLG equation for $\widetilde{\m}$ with $|\widetilde{\m}|=1$
\begin{align*}
\hat{\alpha}_1  \widetilde{\m}_t - \hat{\alpha}_2 \widetilde{\m} \times \widetilde{\m}_t -2Am_{\mathrm{S}}\Delta\widetilde{\m} &= 2Am_{\mathrm{S}}| \nabla \widetilde{\m} |^2\widetilde{\m} + \h - 
( \widetilde{\m}\cdot\h)\widetilde{\m}  && \text{in} \ [0,T] \times \Omega \\
  \partial_{\nu} \widetilde{\m} &= 0 && \text{on} \ [0,T] \times \partial\Omega \\
 \widetilde{\m}(t=0) &=  \widetilde{\m}_0, \ | \widetilde{\m}_0 | = 1 && \text{in} \ \Omega.
\end{align*}

\begin{figure}[p] 
\vspace{1cm}
\centering
\includegraphics[width=5.4cm]{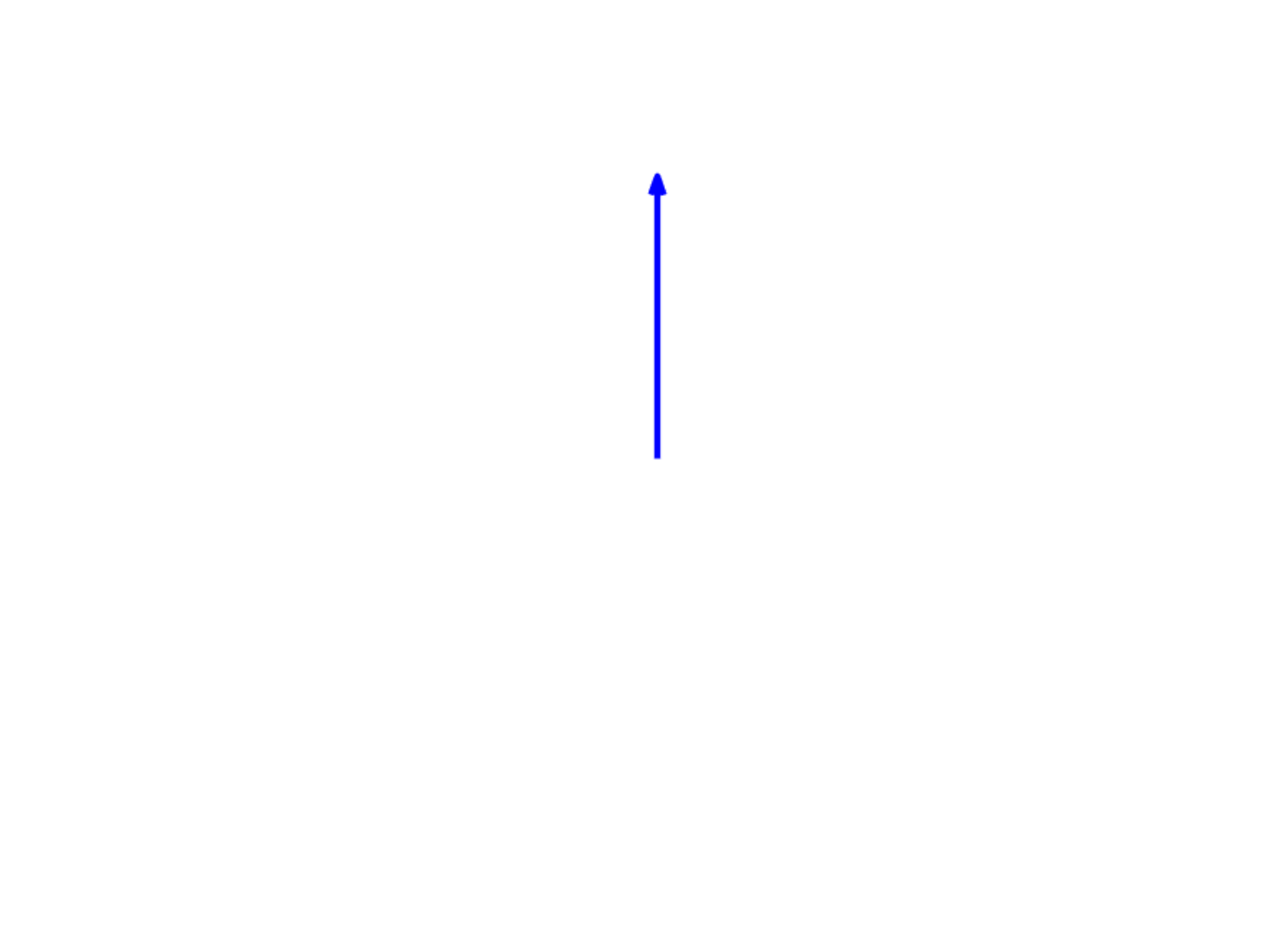}
\hspace{-2cm}
\includegraphics[width=5.4cm]{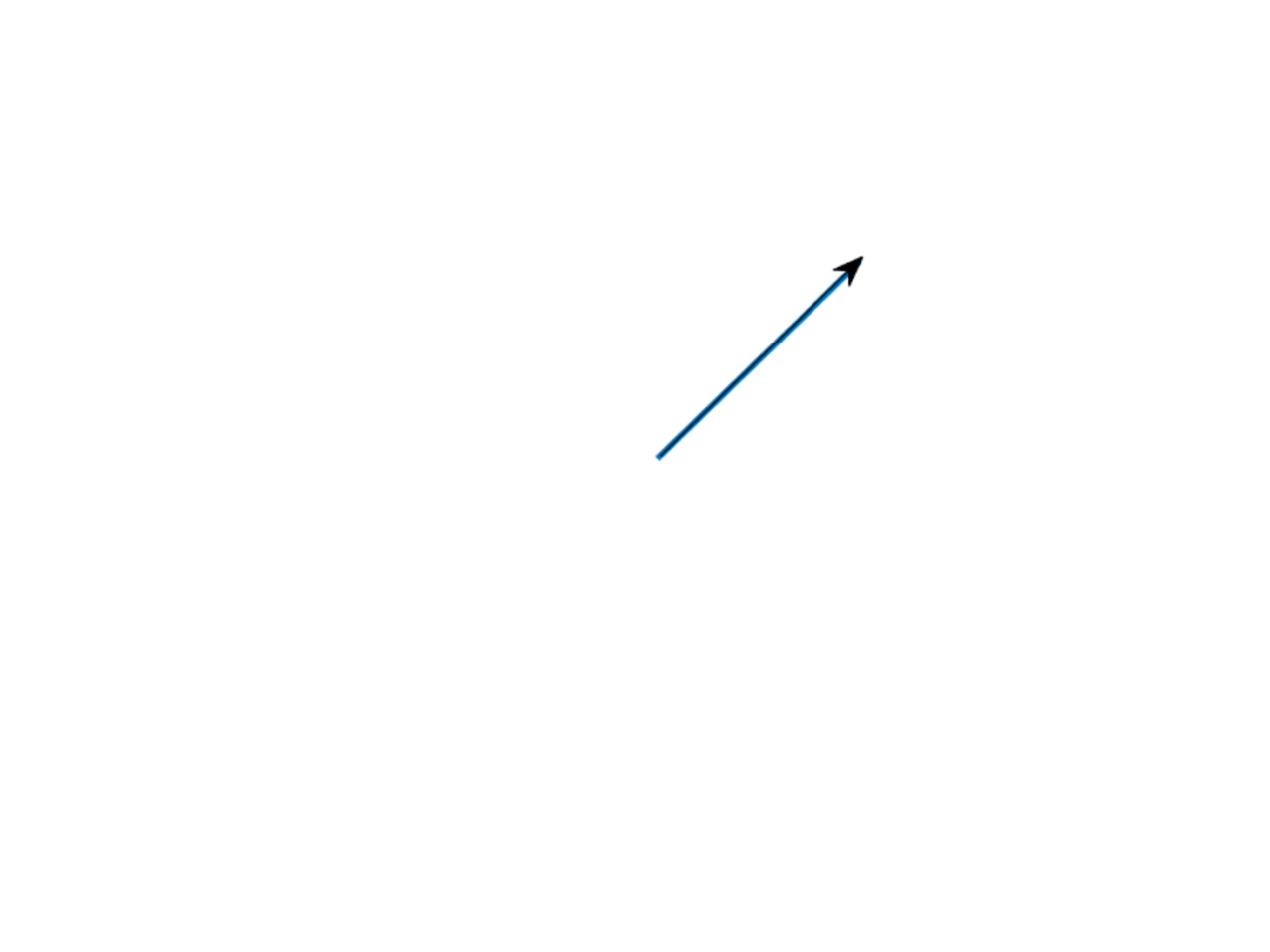}
\hspace{-2cm}
\includegraphics[width=5.4cm]{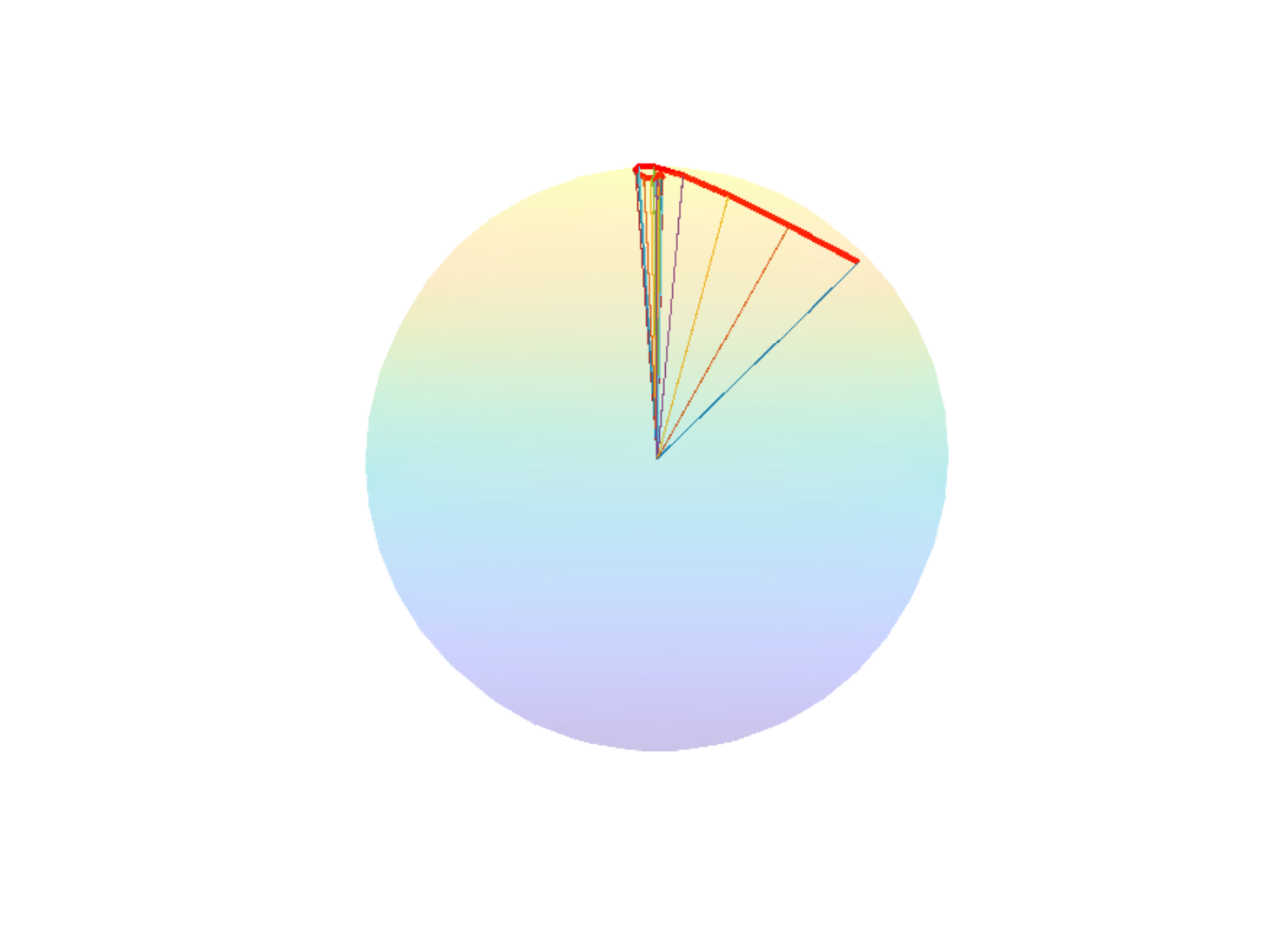}\\
\includegraphics[width=5.4cm]{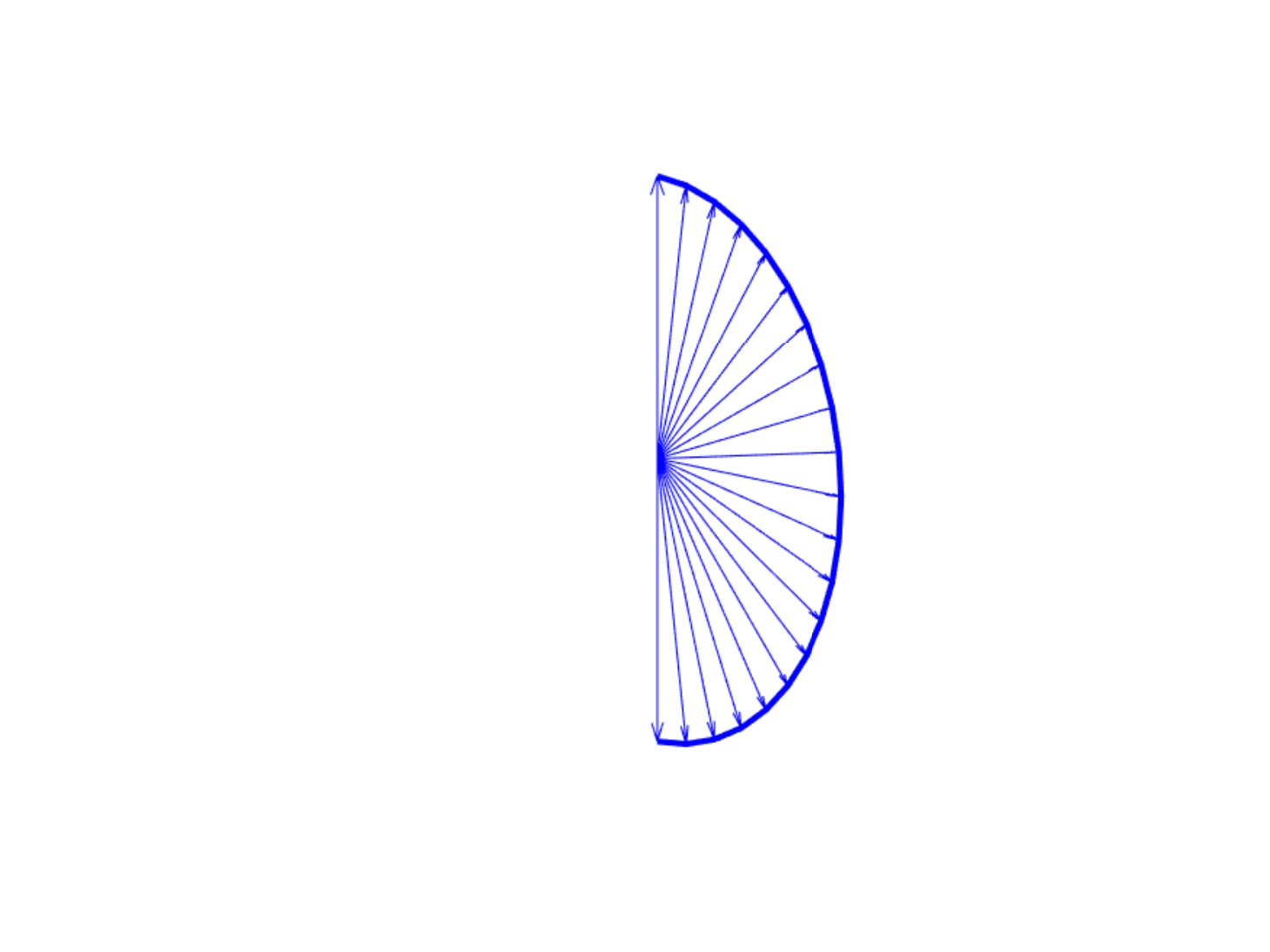}
\hspace{-2cm}
\includegraphics[width=5.4cm]{Figures/fig-part-test1-M0-eps-converted-to}
\hspace{-2cm}
\includegraphics[width=5.4cm]{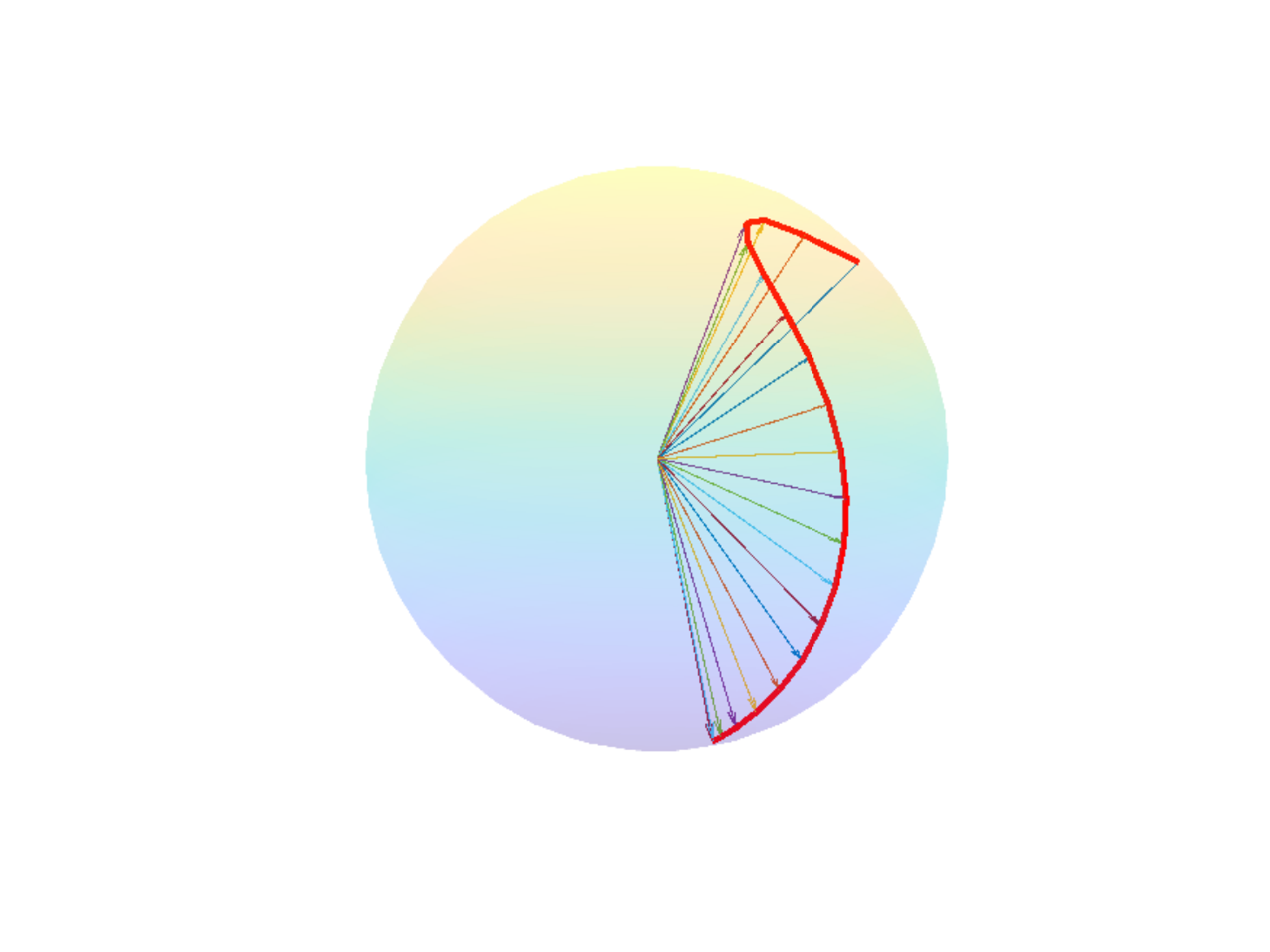}\\
\includegraphics[width=5.4cm]{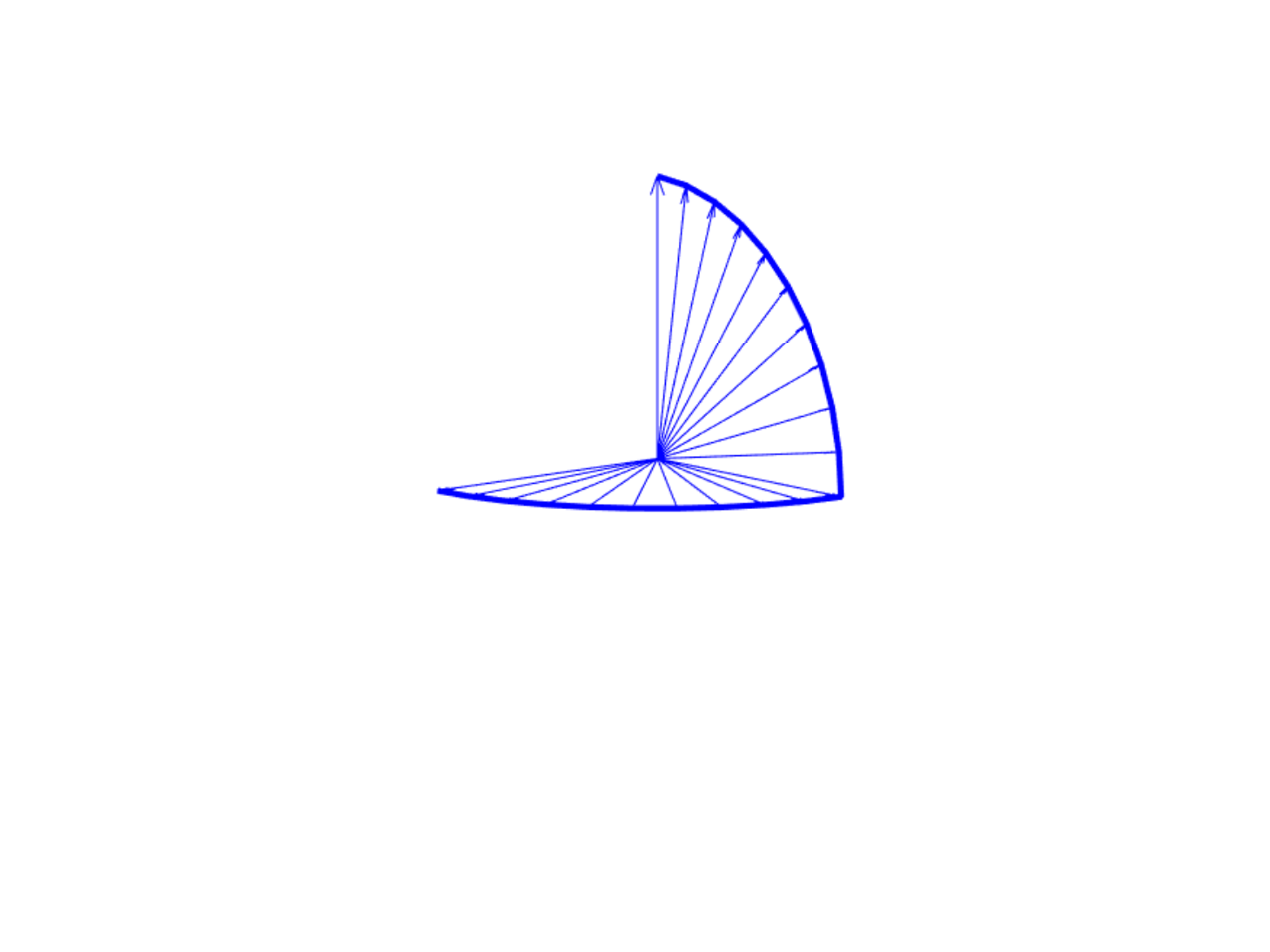}
\hspace{-2cm}
\includegraphics[width=5.4cm]{Figures/fig-part-test1-M0-eps-converted-to}
\hspace{-2cm}
\includegraphics[width=5.4cm]{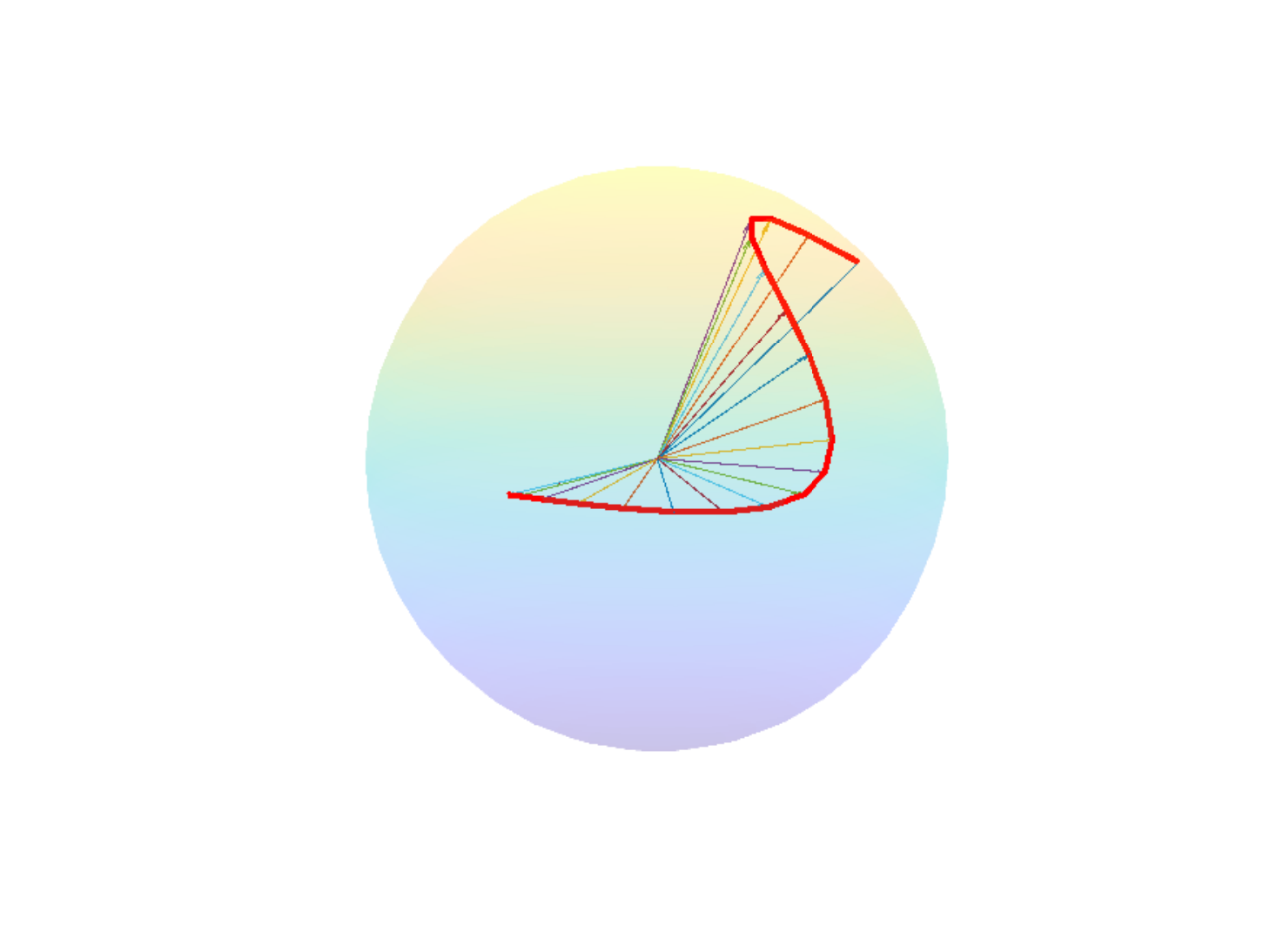}
\caption{Left: applied field $\h$. Middle: initial state $\widetilde{\m}_0$. Right: trajectory of $\widetilde{\m}(t)$.}
\label{fig-particle}
\end{figure}

\begin{figure}[!p] 
\vspace{1.5cm}
\centering
\includegraphics[width=5cm]{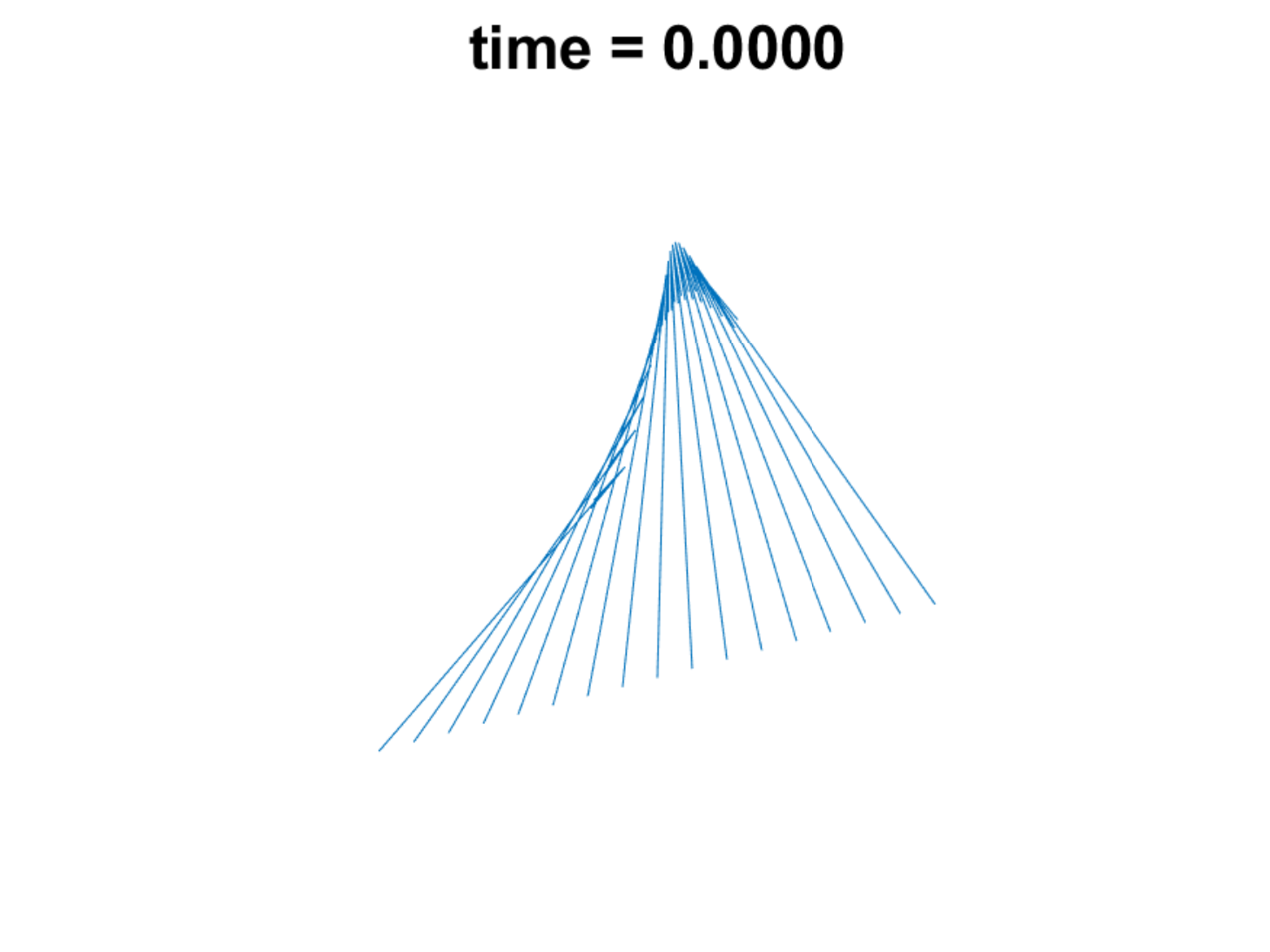}
\hspace{-1.8cm}
\includegraphics[width=5cm]{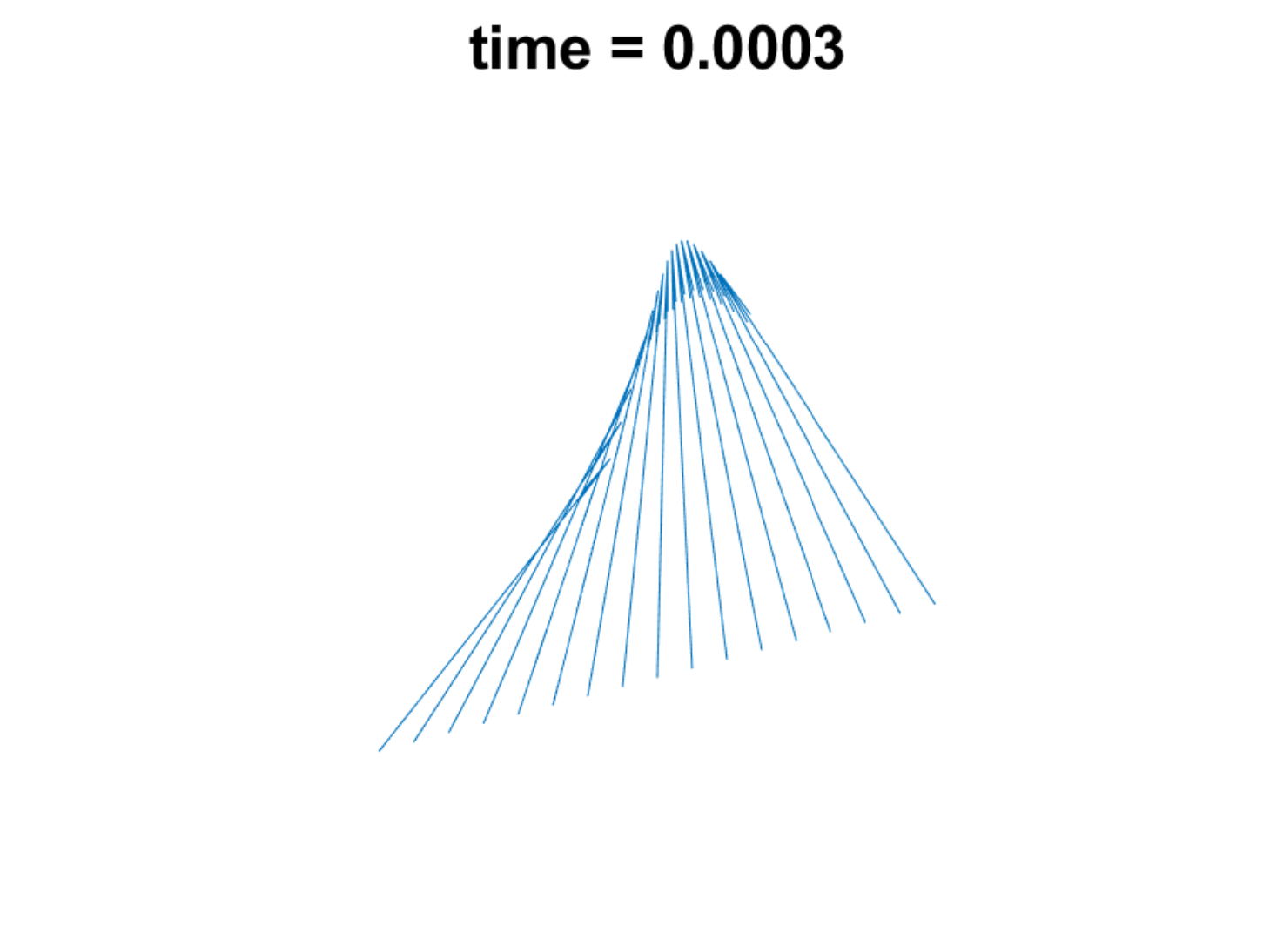}
\hspace{-1.8cm}
\includegraphics[width=5cm]{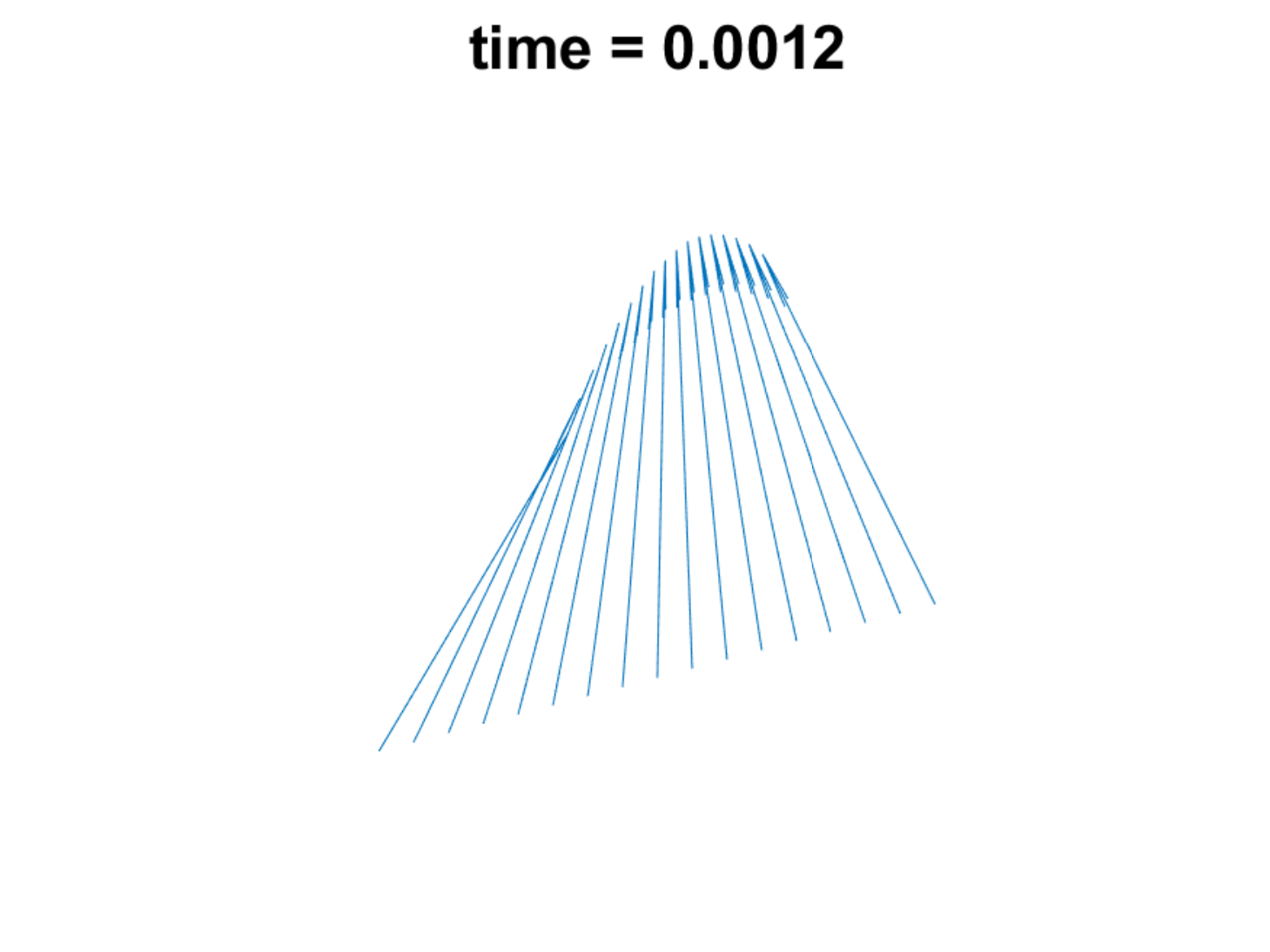}\\
\includegraphics[width=5cm]{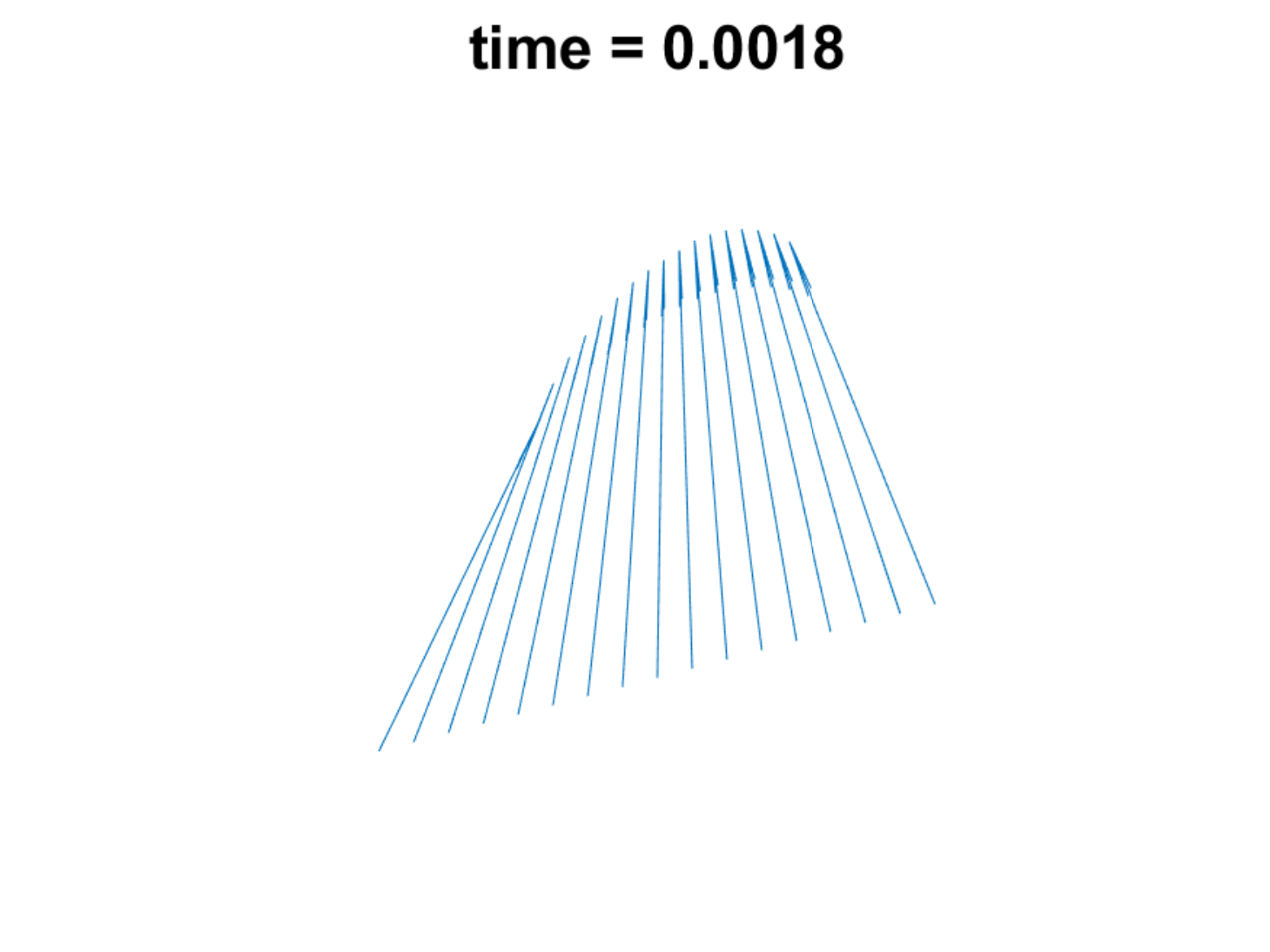}
\hspace{-1.8cm}
\includegraphics[width=5cm]{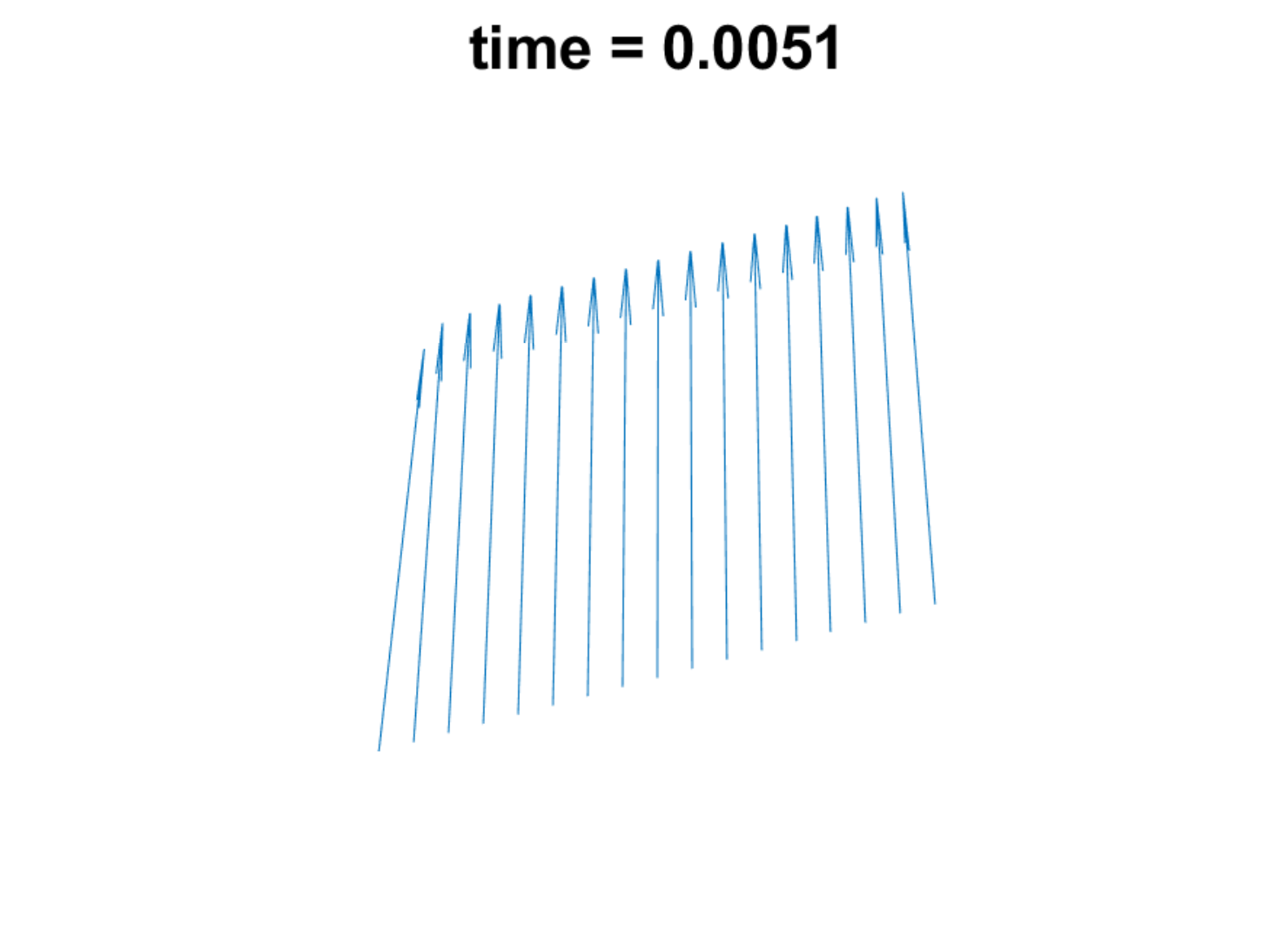}
\hspace{-1.8cm}
\includegraphics[width=5cm]{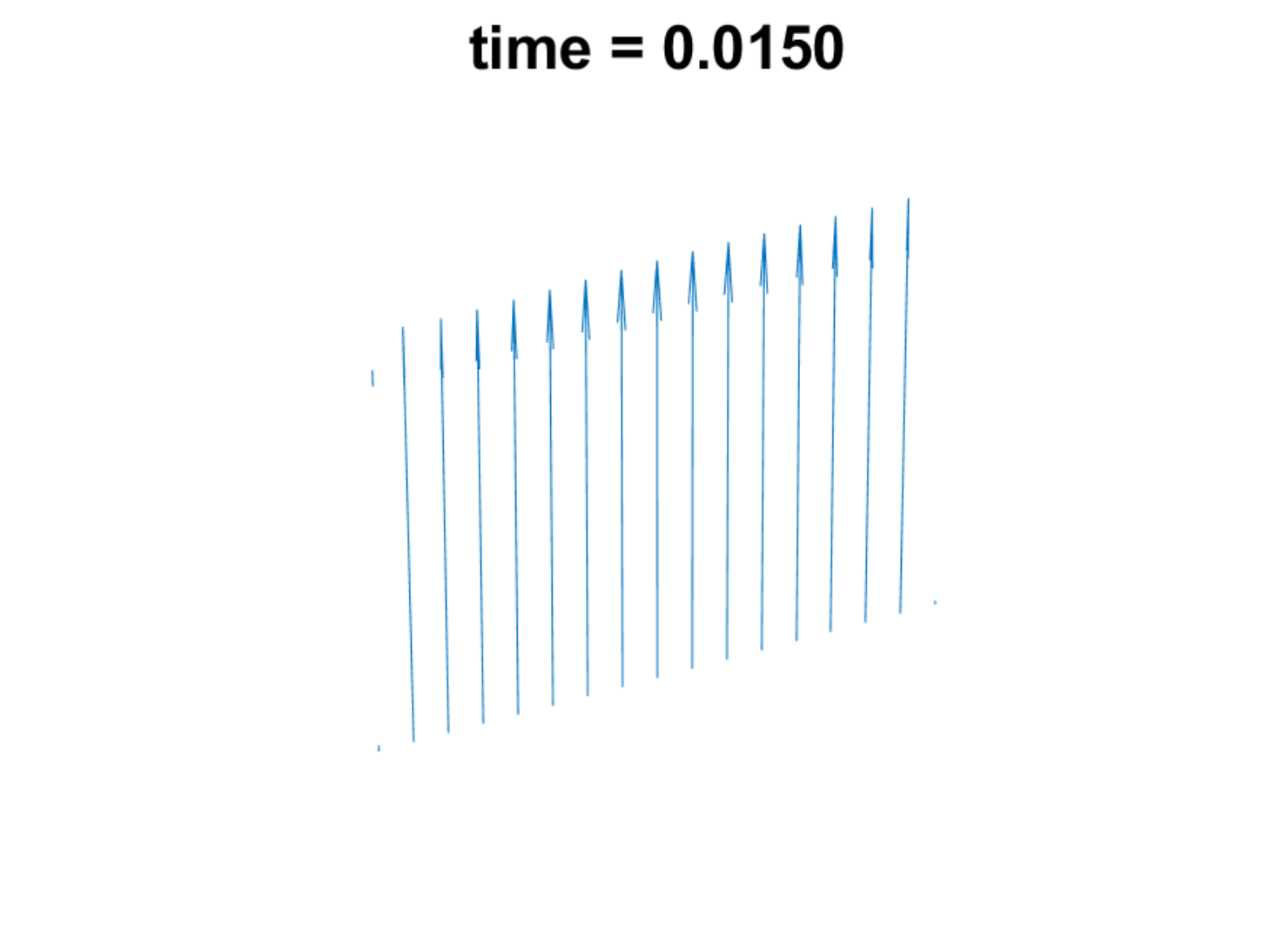}

\vspace{2cm}

\includegraphics[width=5.2cm]{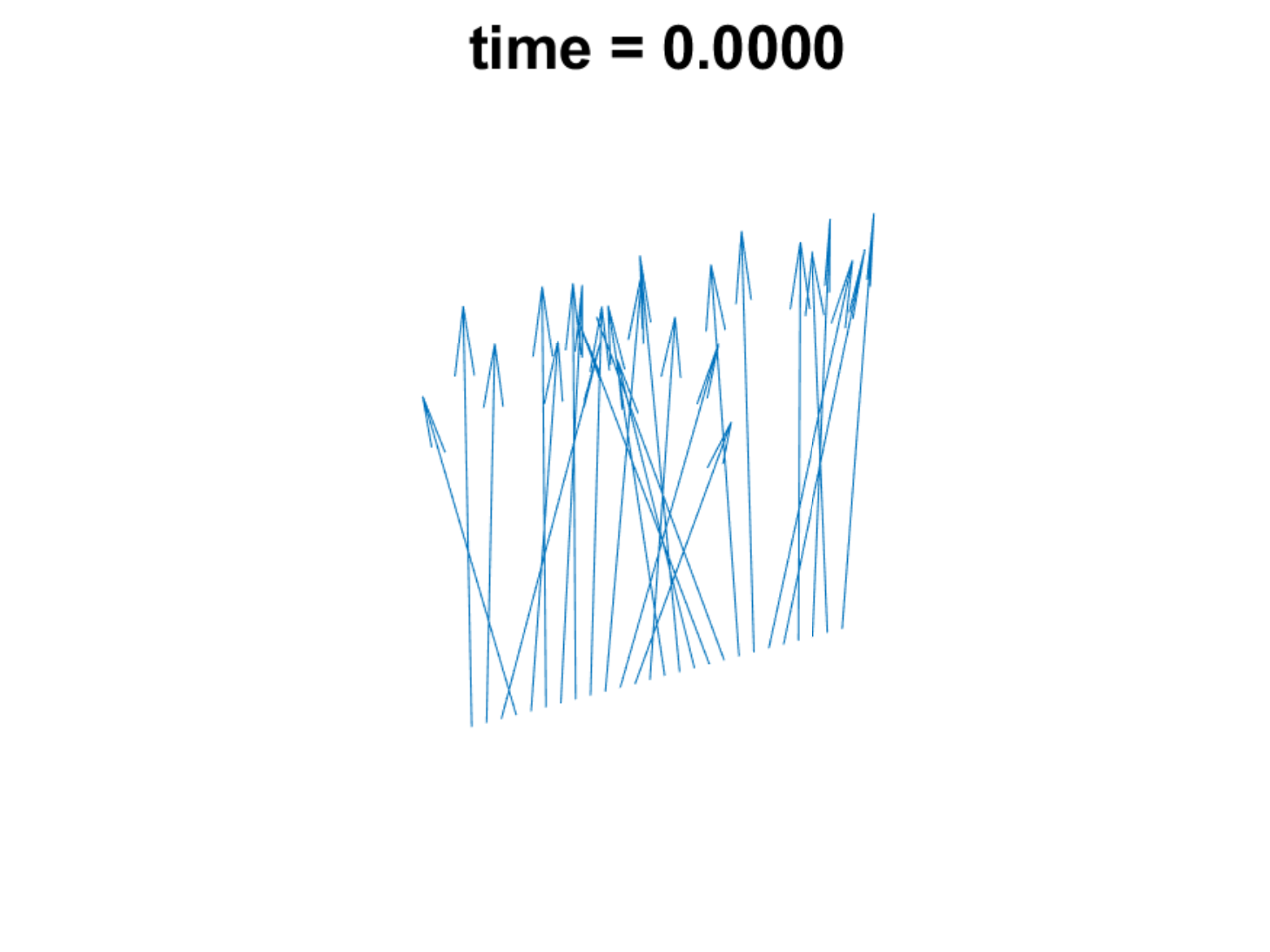}
\hspace{-2cm}
\includegraphics[width=5.2cm]{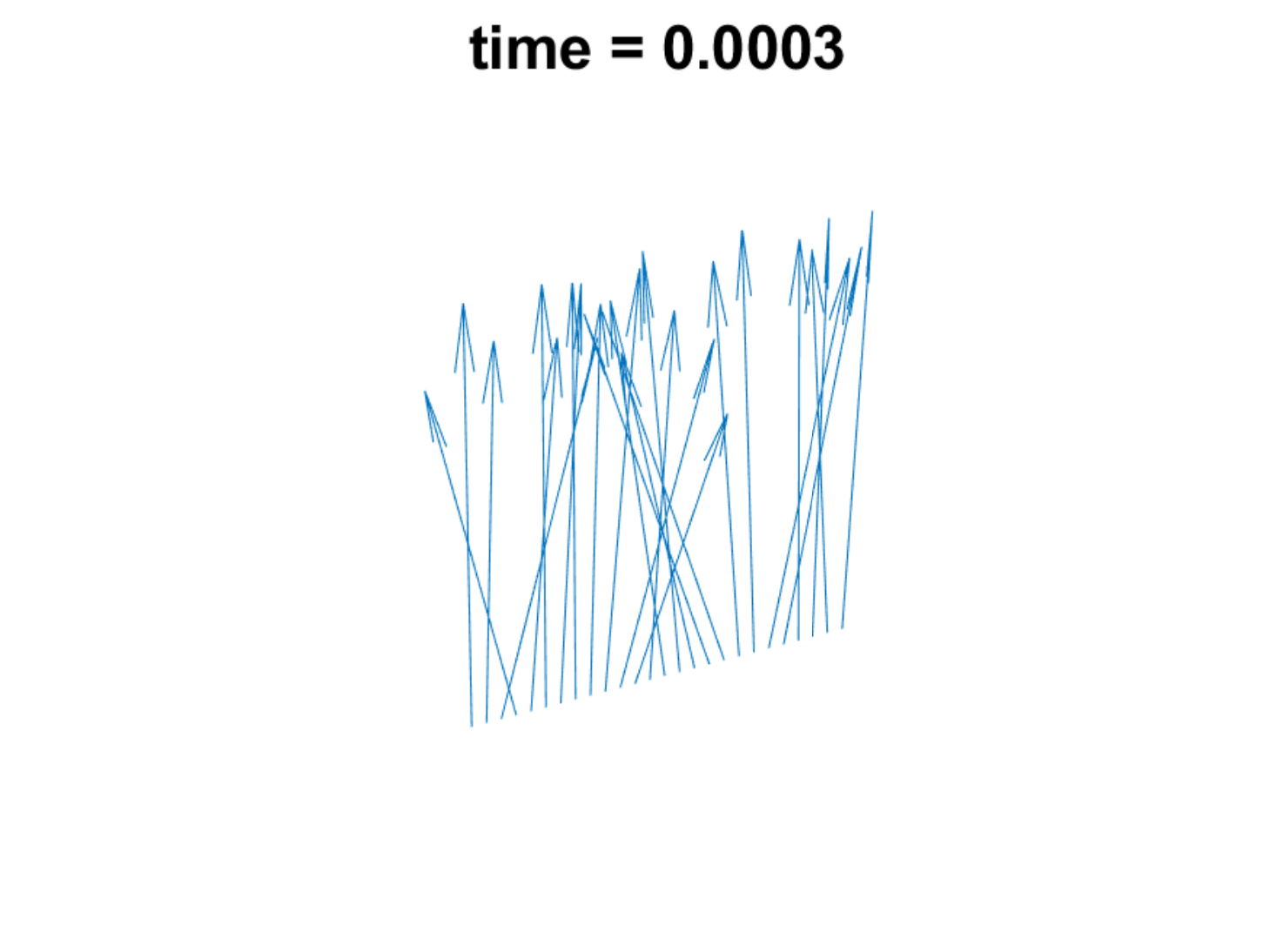}
\hspace{-2cm}
\includegraphics[width=5.2cm]{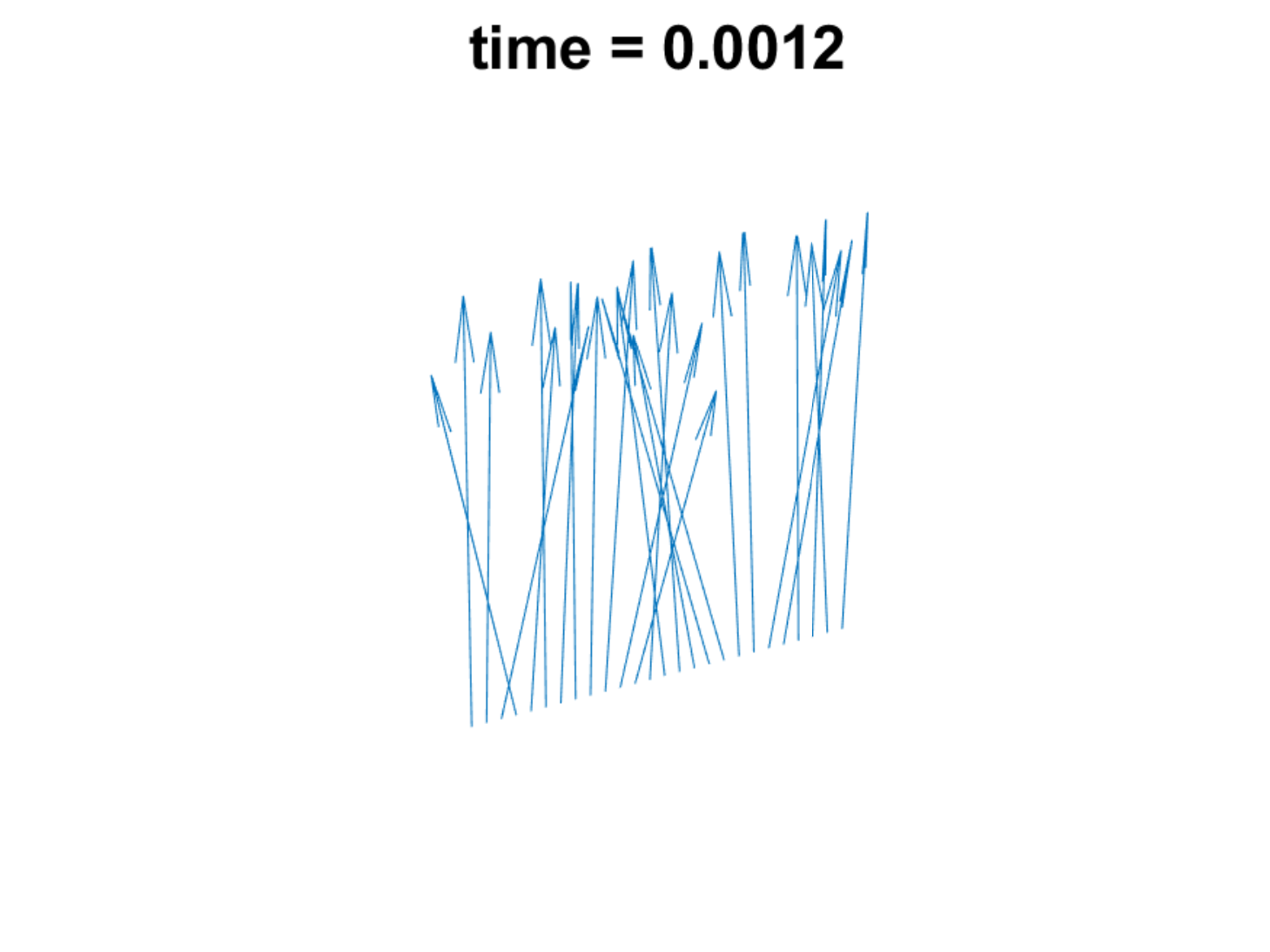}\\
\includegraphics[width=5.2cm]{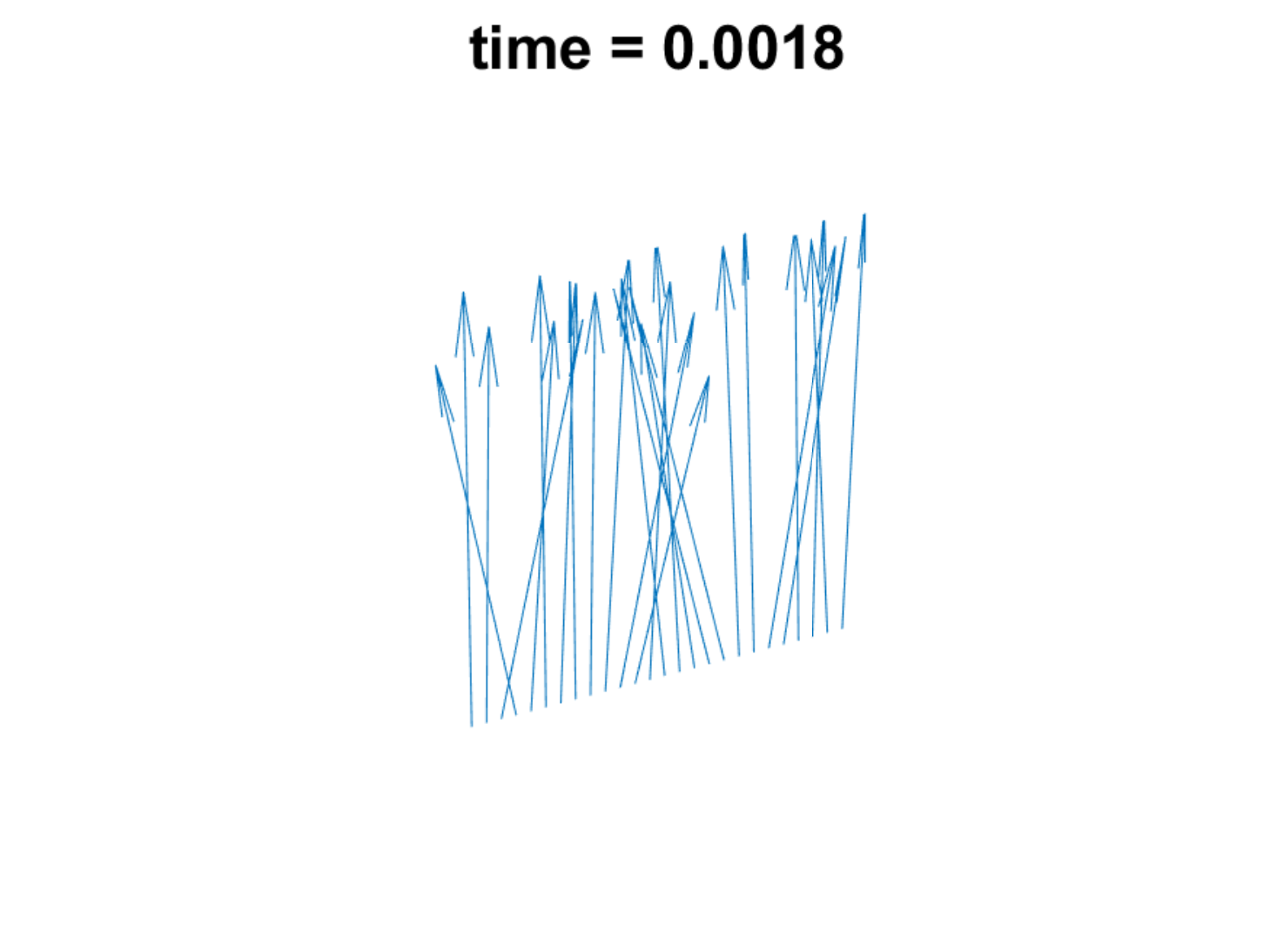}
\hspace{-2cm}
\includegraphics[width=5.2cm]{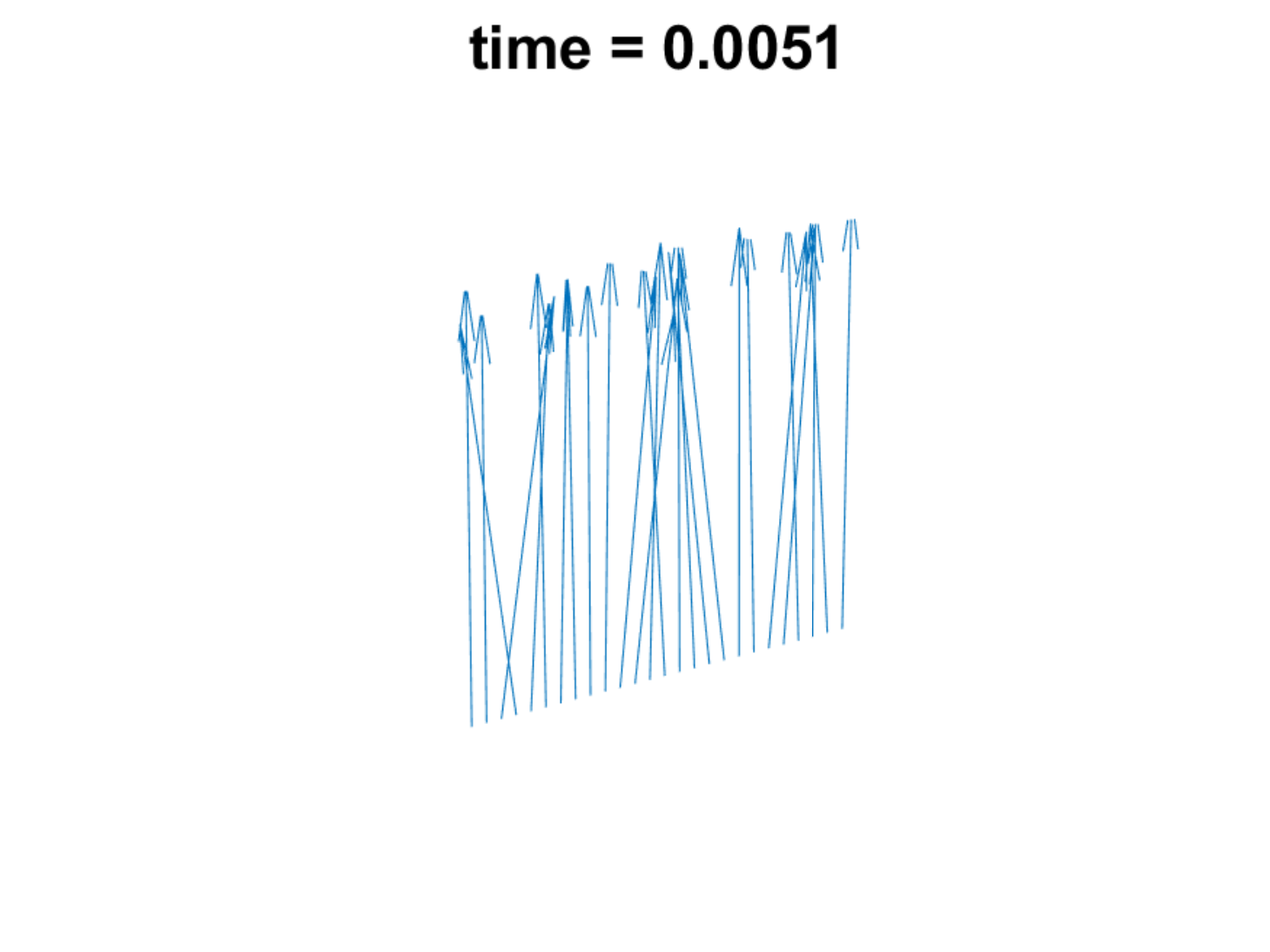}
\hspace{-2cm}
\includegraphics[width=5.2cm]{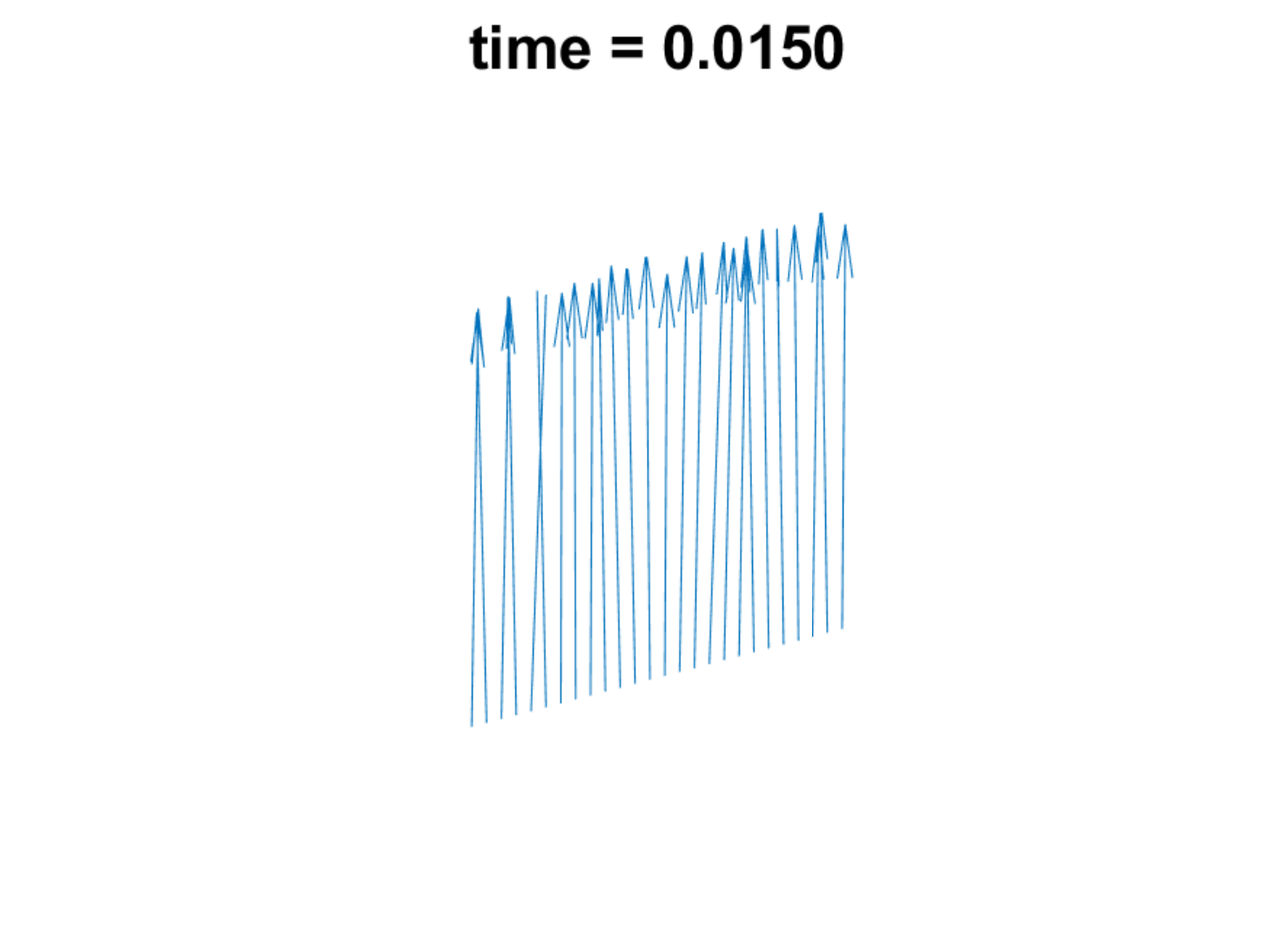}
\caption{Magnetization $\widetilde{\m}$ at different time instances.}
\label{fig-synthetic1}
\end{figure}

\begin{remark}
In the above equation, we use $m_\text{S}=474000$ \blue{J$\text{m}^{-3} \text{T}^{-1}$}. Without knowing $A$, which we expect to be very small, we are left with a very large coefficient in front of the Laplacian term and the nonlinear gradient term relative to the small coefficients $\hat{\alpha}_1\approx 5. 10^{-13}, \hat{\alpha}_2\approx 5. 10^{-12}$. Since we have not yet found a physically meaningful exchange constant $A$, particularly in view of our different approach to the model for the system function in comparison to, e.g., \cite{tktb12, tk18}, we simply set $A=0$ for Tests \ref{fig-particle}-\ref{fig-synthetic1} in this section. Even with $A=0$, we observe the dominating relaxation effect. However, we wish to emphasize that Tests 1-3, Figures \ref{test1}-\ref{test33} were run with nonzero diffusion and gradient terms (c.f. \eqref{LLG-dropMs}).
\end{remark}

In Figure \ref{fig-particle}, the left columns display three states of the applied magnetic field $\h$, namely, one static field in $e_3$-direction and two different time-dependent fields. Starting from an initial state (homogeneous in space), the magnetization vector $\widetilde{\m}$ (right columns) follows the trajectory of $\h$ with a delay known as the relaxation effect. The length of $\widetilde{\m}$ is not fully preserved, but it is obviously not far from the unit length. Figure \ref{fig-synthetic1} displays the progression of the magnetization $\widetilde{\m}$ when the static field in $e_3$-direction is applied to a space-inhomogeneous initial state $\widetilde{\m}_0$ (top six plots) and another $\widetilde{\m}_0$ distributed randomly in space (bottom six plots).

\newpage
\begin{table}
\bigskip\bigskip
\caption{Common physical parameters.}
\centering
\begin{tabular}{ l  c l l}
\hline
\rule{0pt}{15pt}
Parameter& \qquad
\qquad\qquad& Value & Unit\\[1ex]
\hline
\rule{0pt}{12pt}
Magnetic permeability & $\mu_0$ & 4$\pi\times 10^{-7}$ &H $\text{m}^{-1}$\\
\rule{0pt}{12pt}
Sat. magnetization& $m_{\mathrm{S}}$ & 474 000 & \blue{J $\text{m}^{-3} \text{T}^{-1}$}\\
\rule{0pt}{12pt}
Gyromagnetic ratio &$\gamma$ & 1.75$\times 10^{11}$ &rad $\text{s}^{-1}$\\
\rule{0pt}{12pt}
Damping parameter &$ \alpha_\text{D}$& 0.1\\
\hline
\rule{0pt}{12pt}
Field of view & $\Omega$ & [-0.006, 0.006] &m\\
\rule{0pt}{12pt}
Max observation time & T & 0.03$\times 10^{-3}$ &s\\
\rule{0pt}{12pt}
External field strength & $|\h|$ & $10^{-4}$& T\\
\hline
\end{tabular}
\label{tab-physic}
\end{table}

\clearpage

\clearpage

\subsection{Reconstruction in all-at-once and reduced settings} \label{Sec:recos}
In this section, we compare the reconstruction outcome from the all-at-once and reduced settings in case of exact measured data. Staying with the implementation method from Section \ref{sec-LLGsolver}, we carry out the Landweber approximation for Test 2 in Table \ref{tab-1}. For the measurement process, we define the observation operator via: $\mu_0=1, \tilde{a}_l=1, c_k=1, \pr=(1,1,1)$. The initial guess for the state is $\m_\text{init}=\m_0$.

\begin{figure}[!htb] 
\vspace{1.5cm}
\centering
\includegraphics[width=4.3cm]{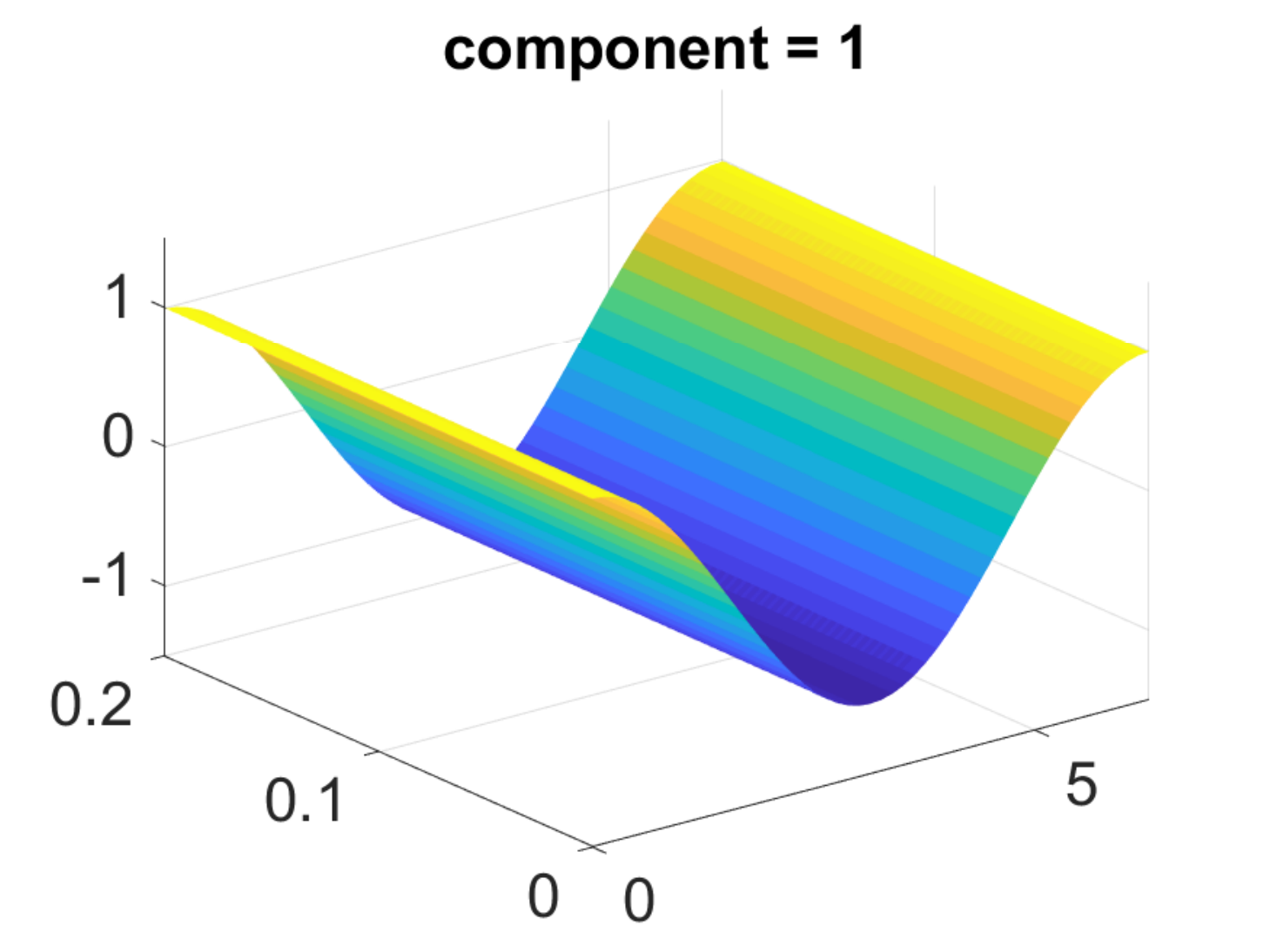}
\hspace{-0.5cm}
\includegraphics[width=4.3cm]{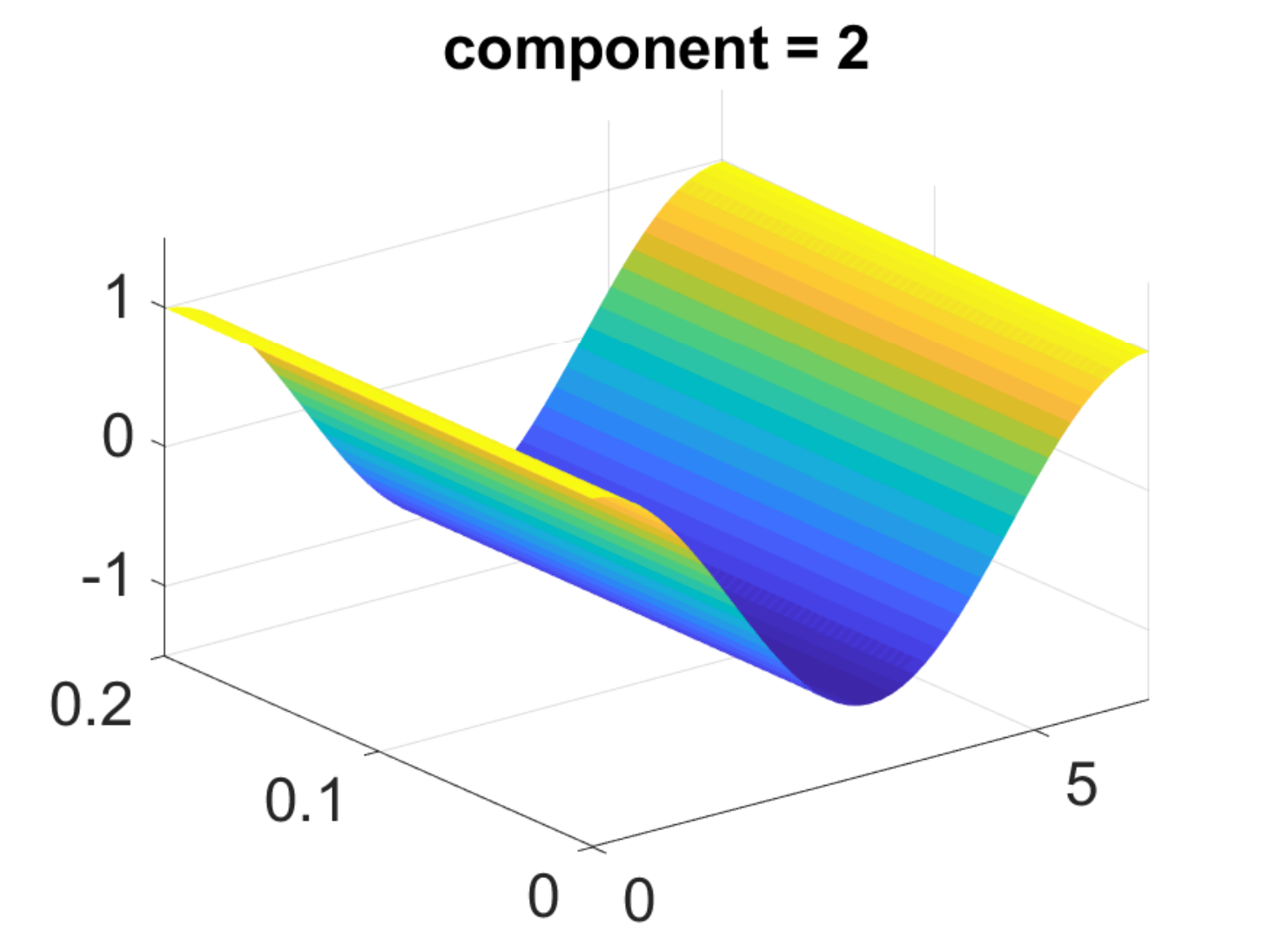}
\hspace{-0.5cm}
\includegraphics[width=4.3cm]{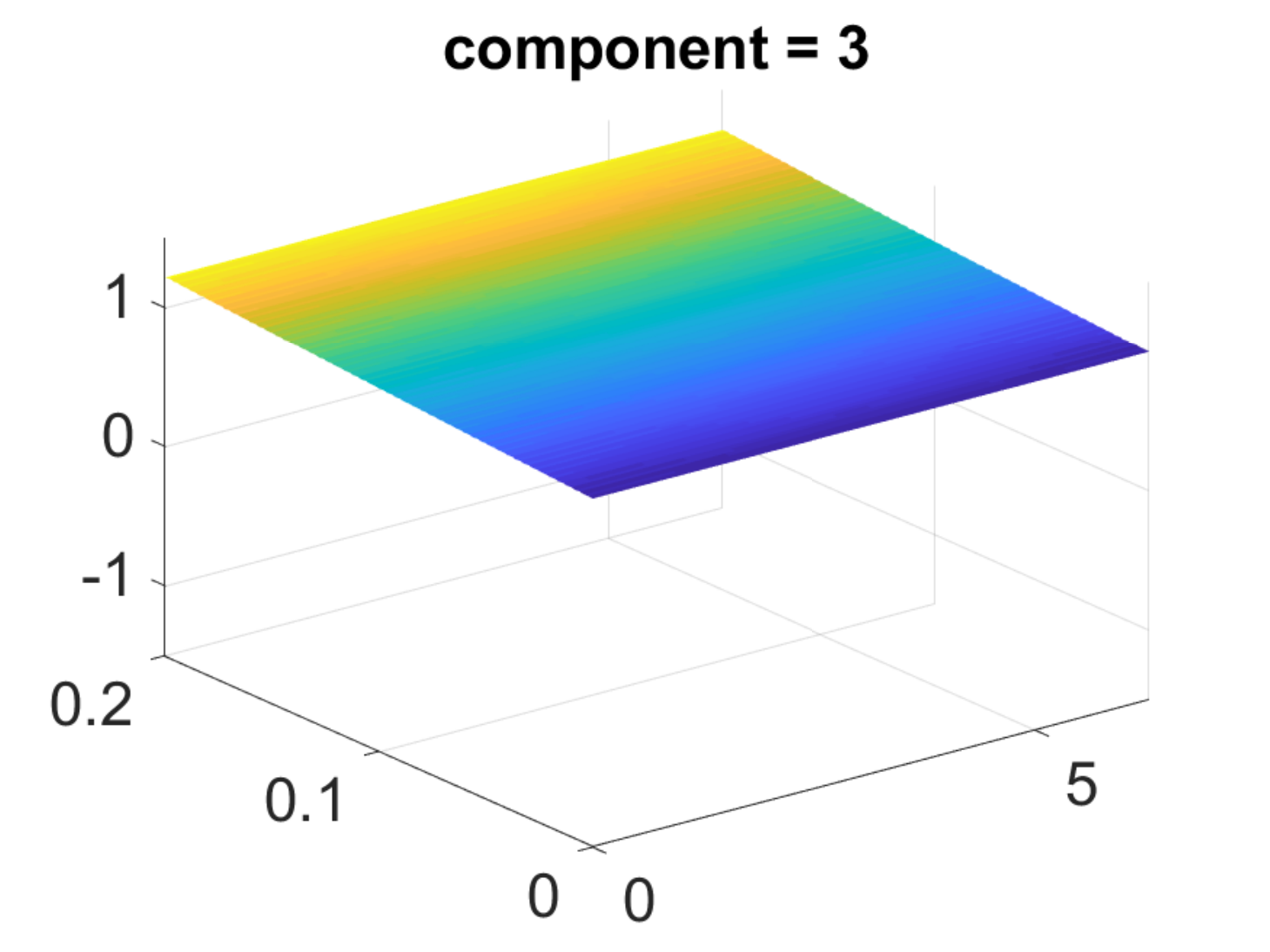}\\[1ex]
\includegraphics[width=4.3cm]{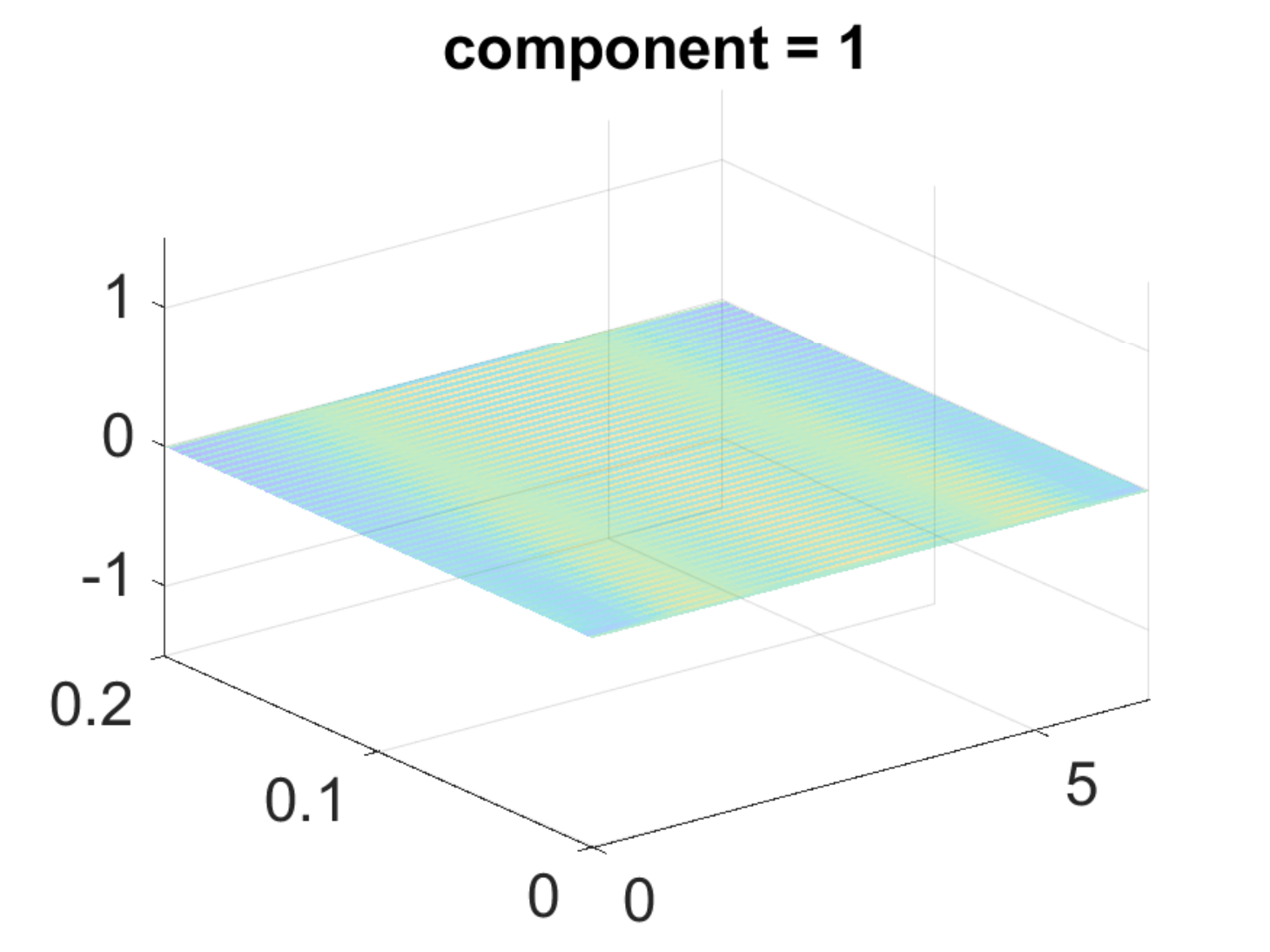}
\hspace{-0.5cm}
\includegraphics[width=4.3cm]{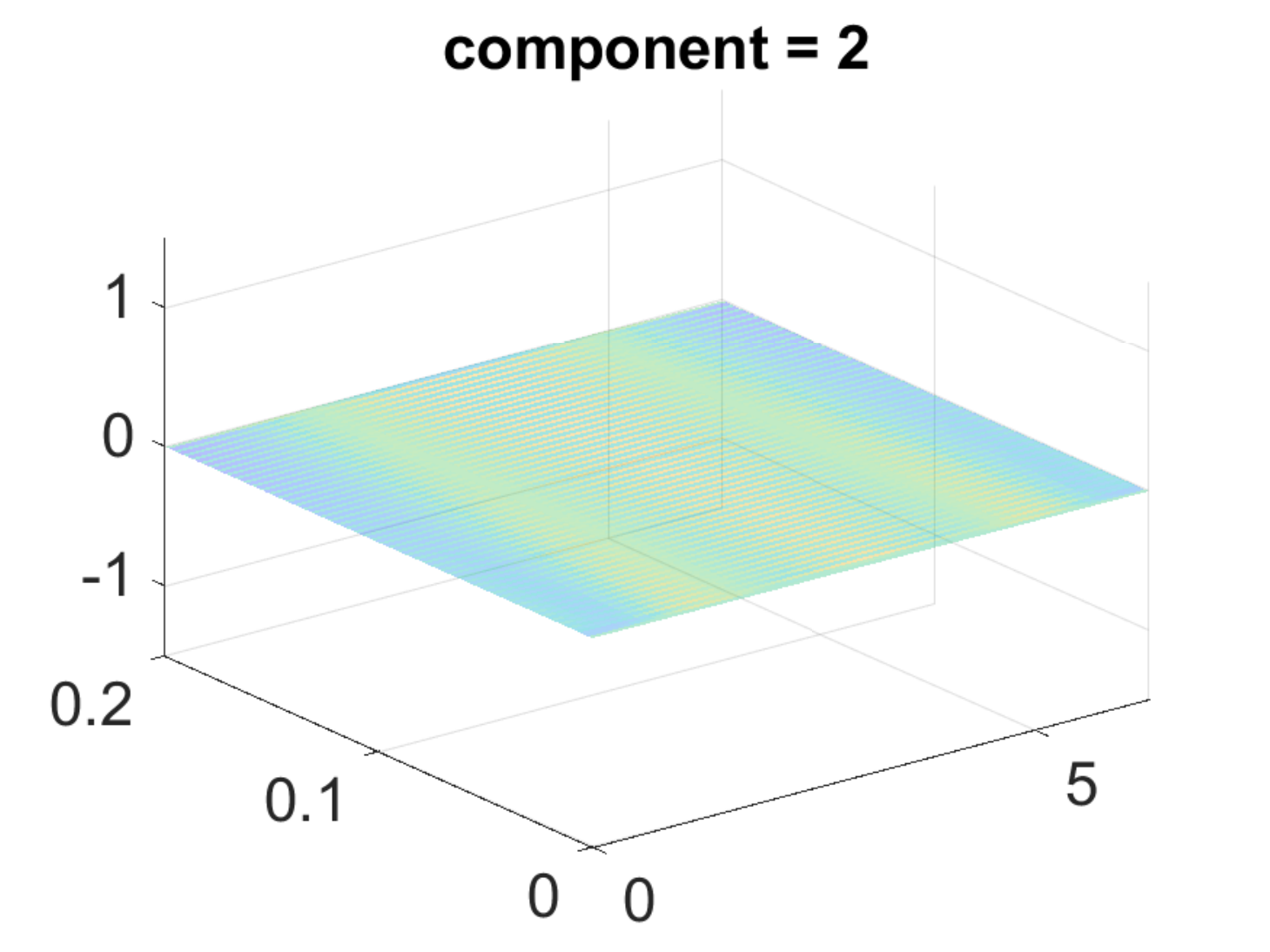}
\hspace{-0.5cm}
\includegraphics[width=4.3cm]{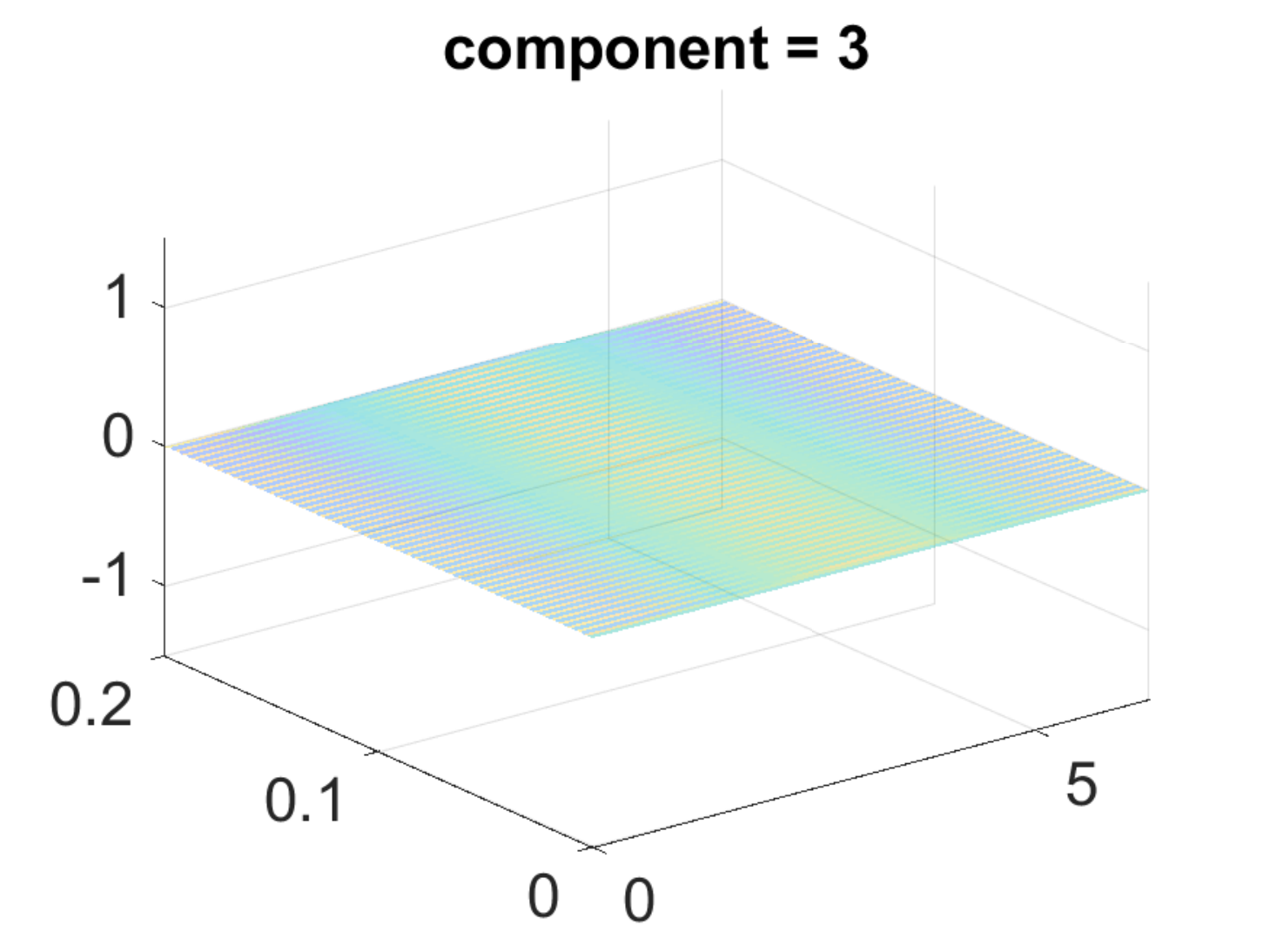}
\caption{Test 2, all-at-once setting. Reconstructed $\mv$ (top) and $\mv-\mv_{\text{exact}}$ (bottom) plotted against space (x-axis) and time (y-axis).}
\label{test-res-aao-m}
\end{figure}

\begin{figure}[!htb] 
\centering
\includegraphics[width=4.3cm]{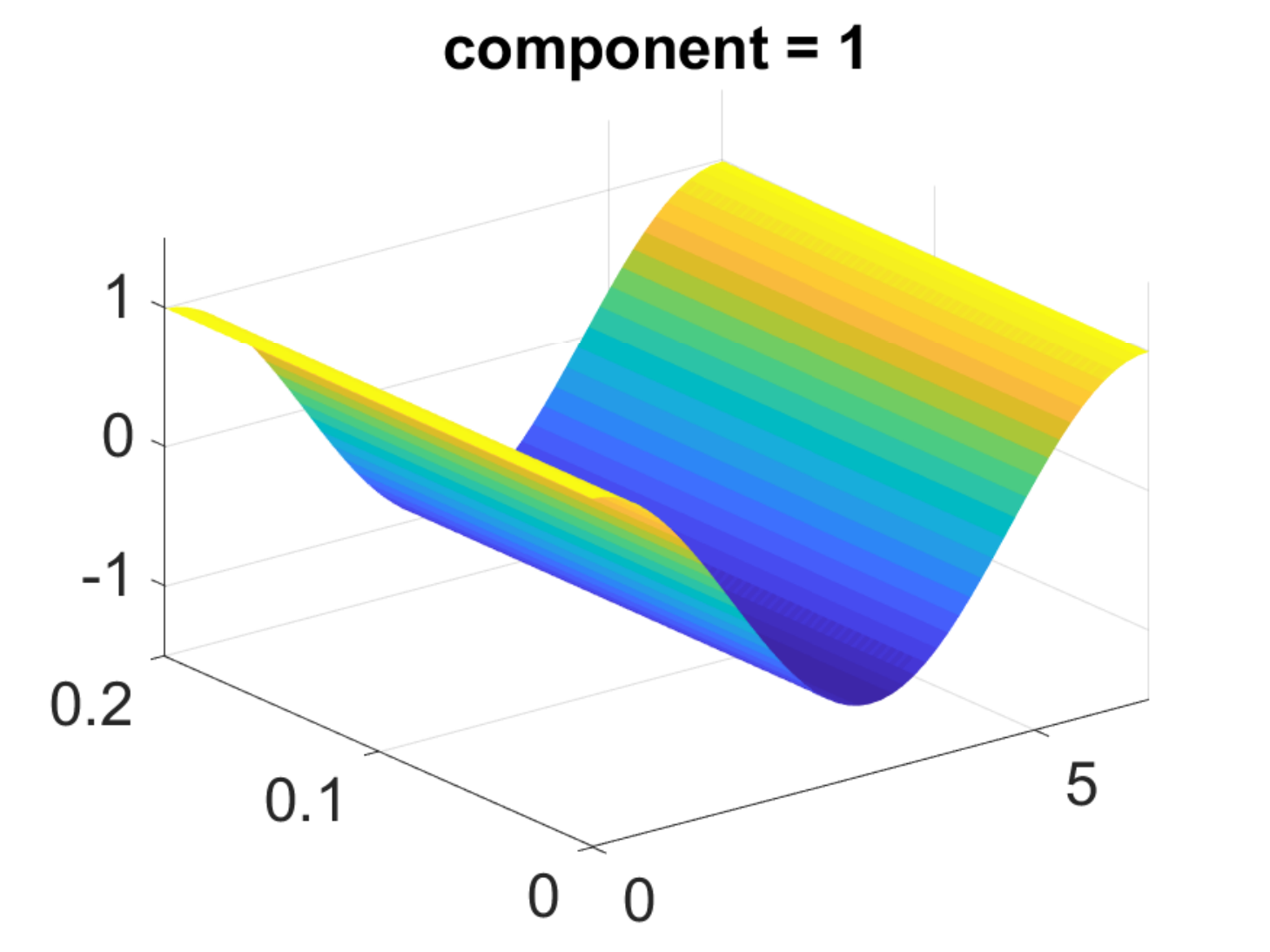}
\hspace{-0.5cm}
\includegraphics[width=4.3cm]{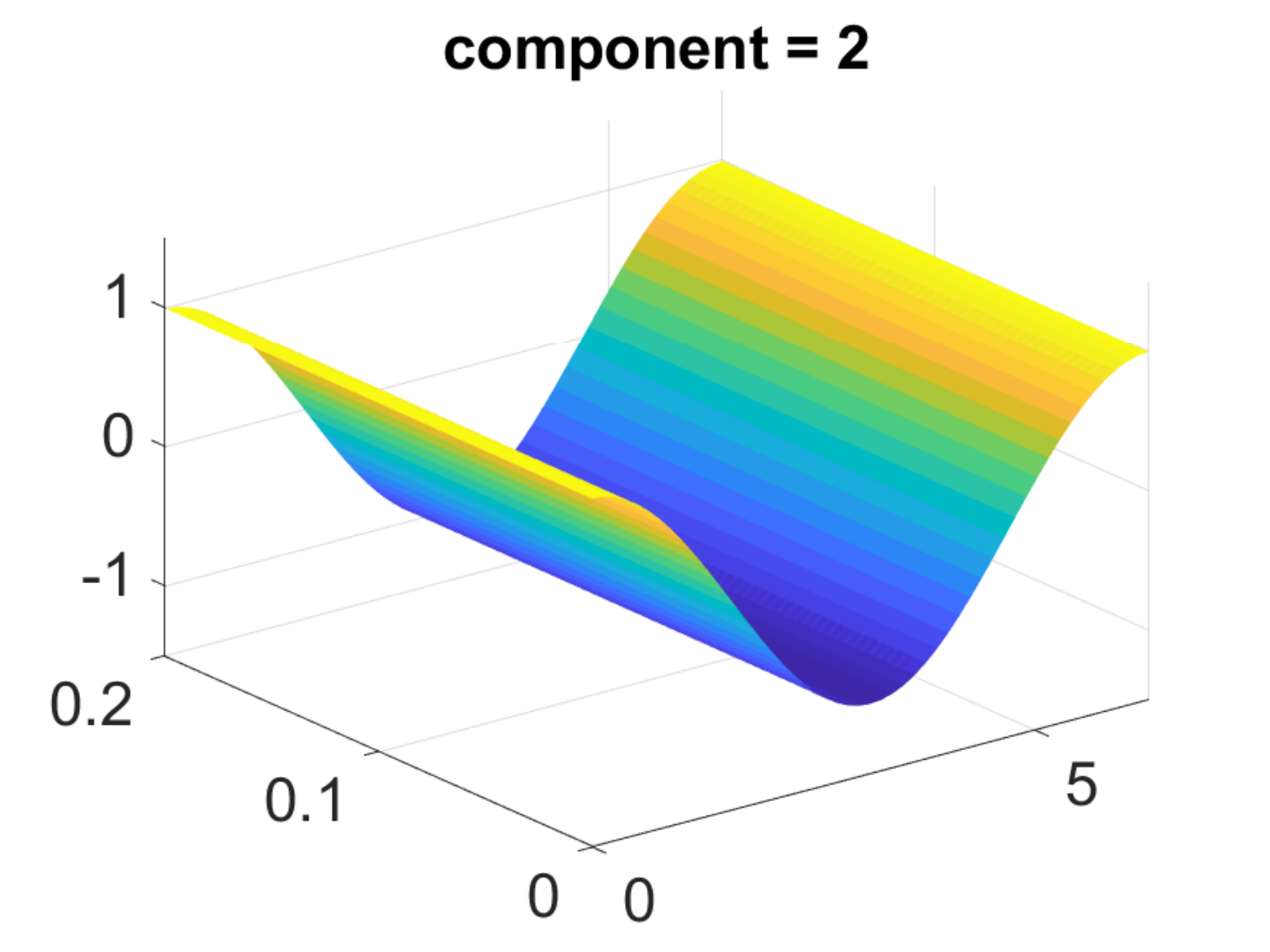}
\hspace{-0.5cm}
\includegraphics[width=4.3cm]{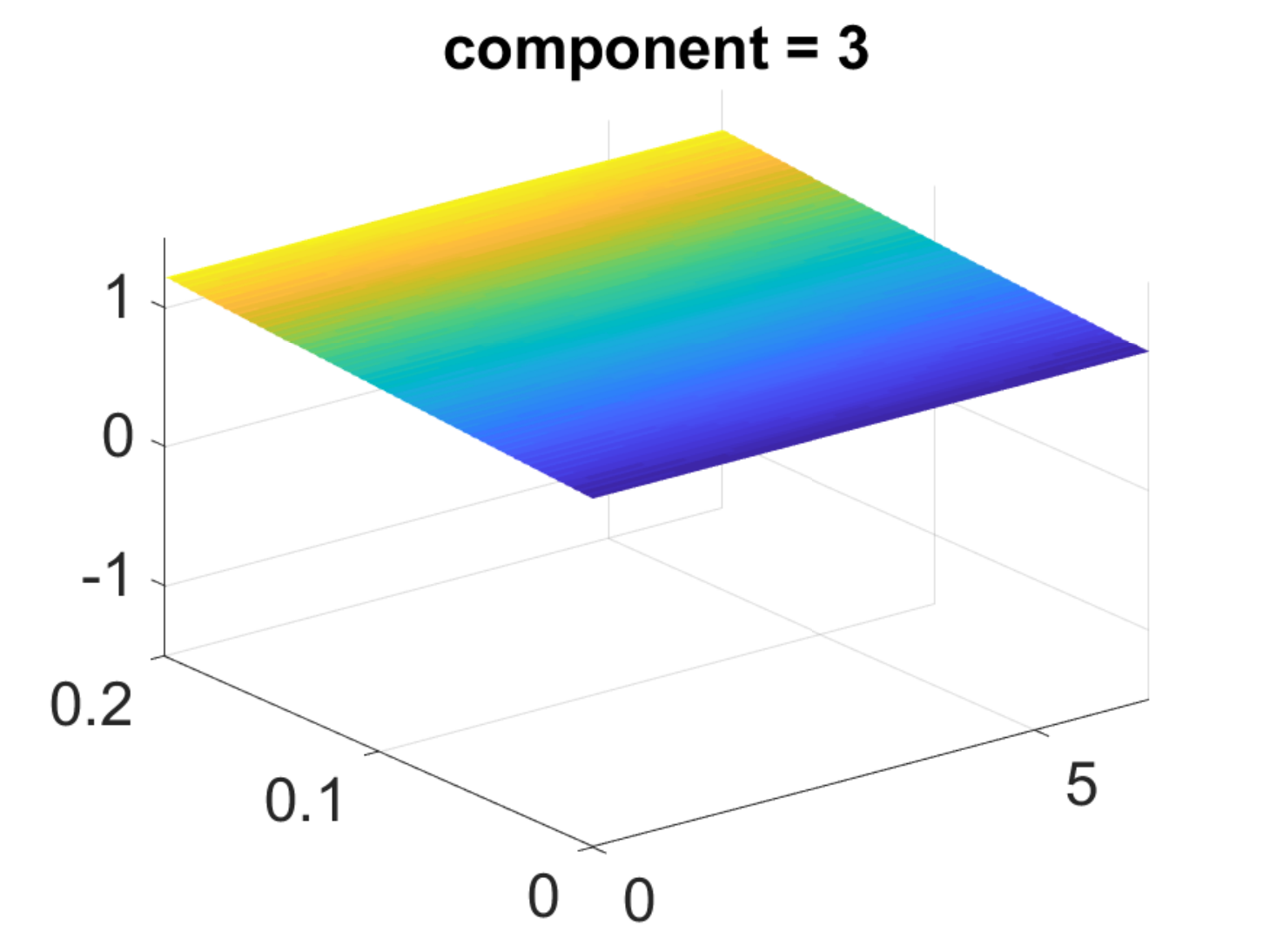}\\[1ex]
\includegraphics[width=4.3cm]{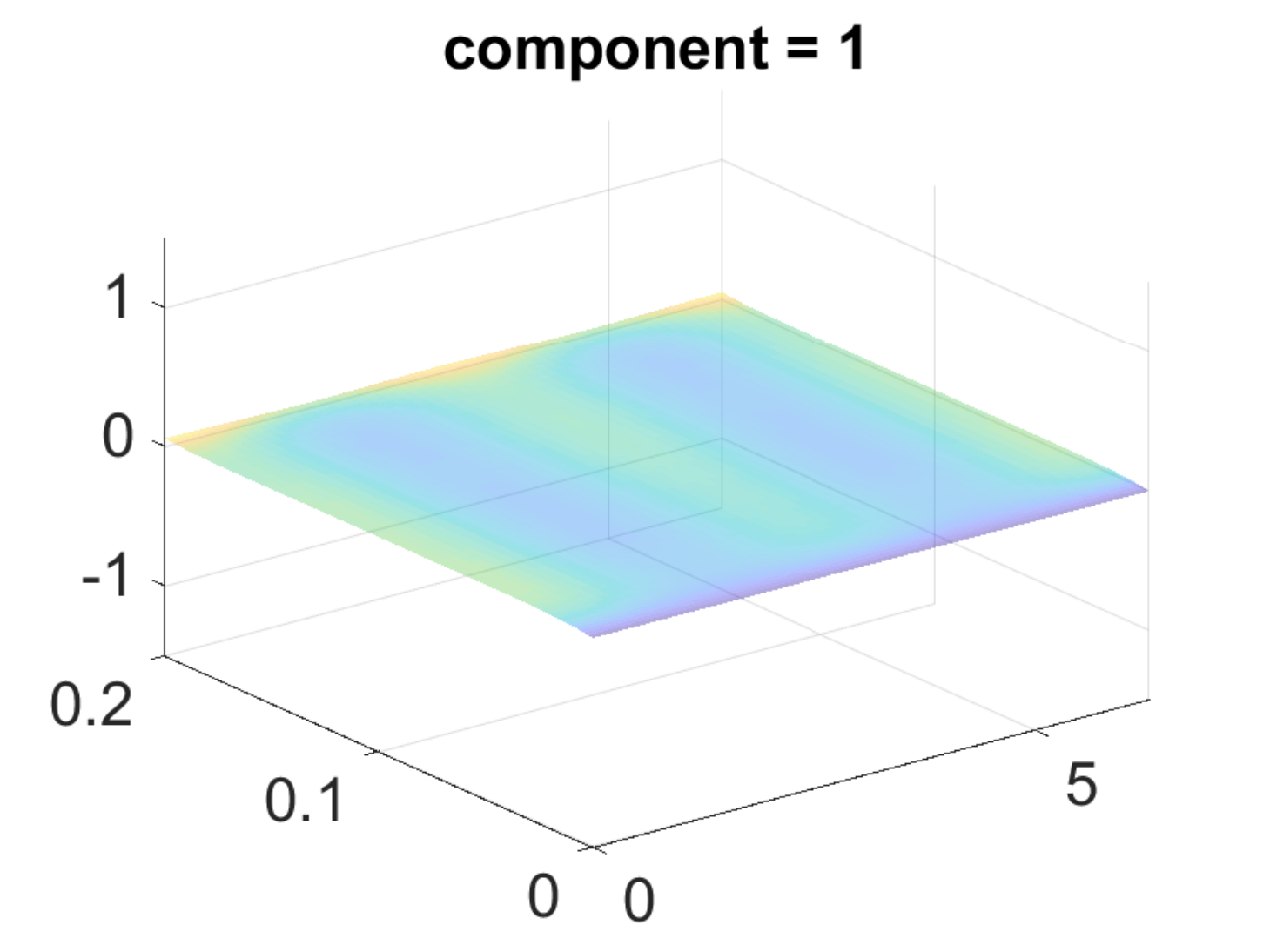}
\hspace{-0.5cm}
\includegraphics[width=4.3cm]{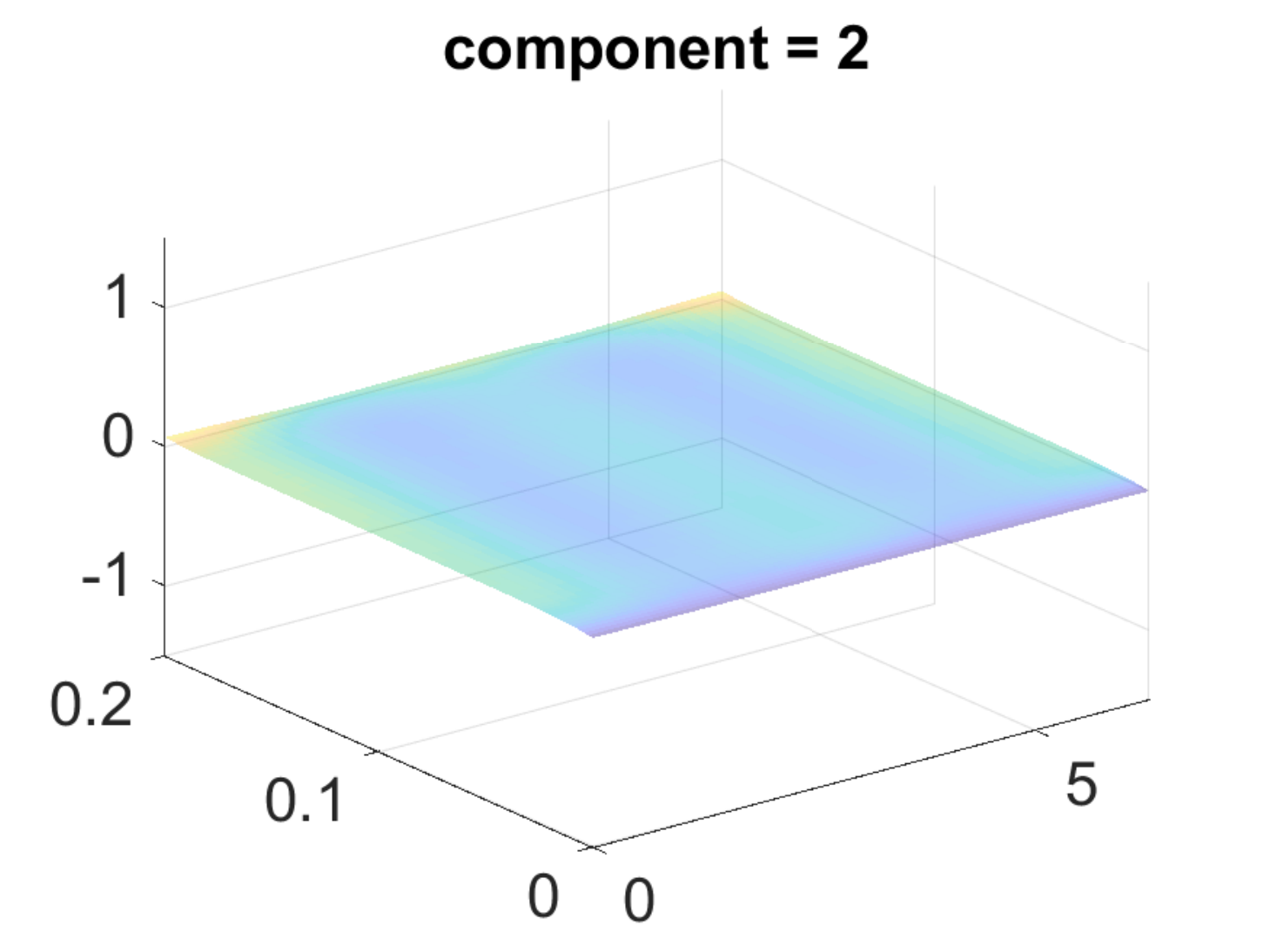}
\hspace{-0.5cm}
\includegraphics[width=4.3cm]{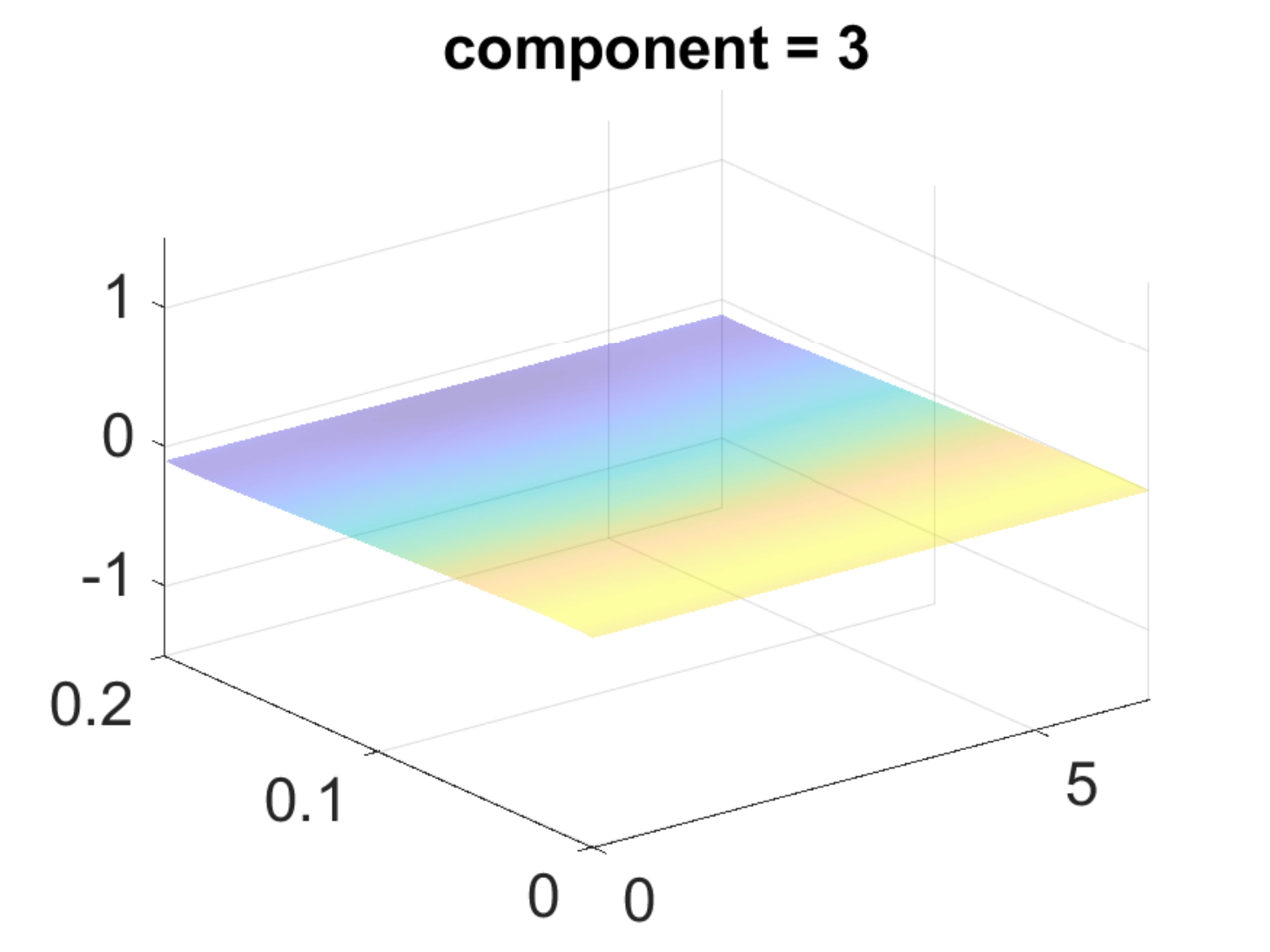}
\caption{Test 2, reduced setting. Reconstructed $\mv$ (top) and $\mv-\mv_{\text{exact}}$ (bottom) plotted against space (x-axis) and time (y-axis).}
\label{test-res-red-m}
\end{figure}

\begin{figure}[!htb] 
\vspace{1cm}
\centering
\includegraphics[width=5.5cm]{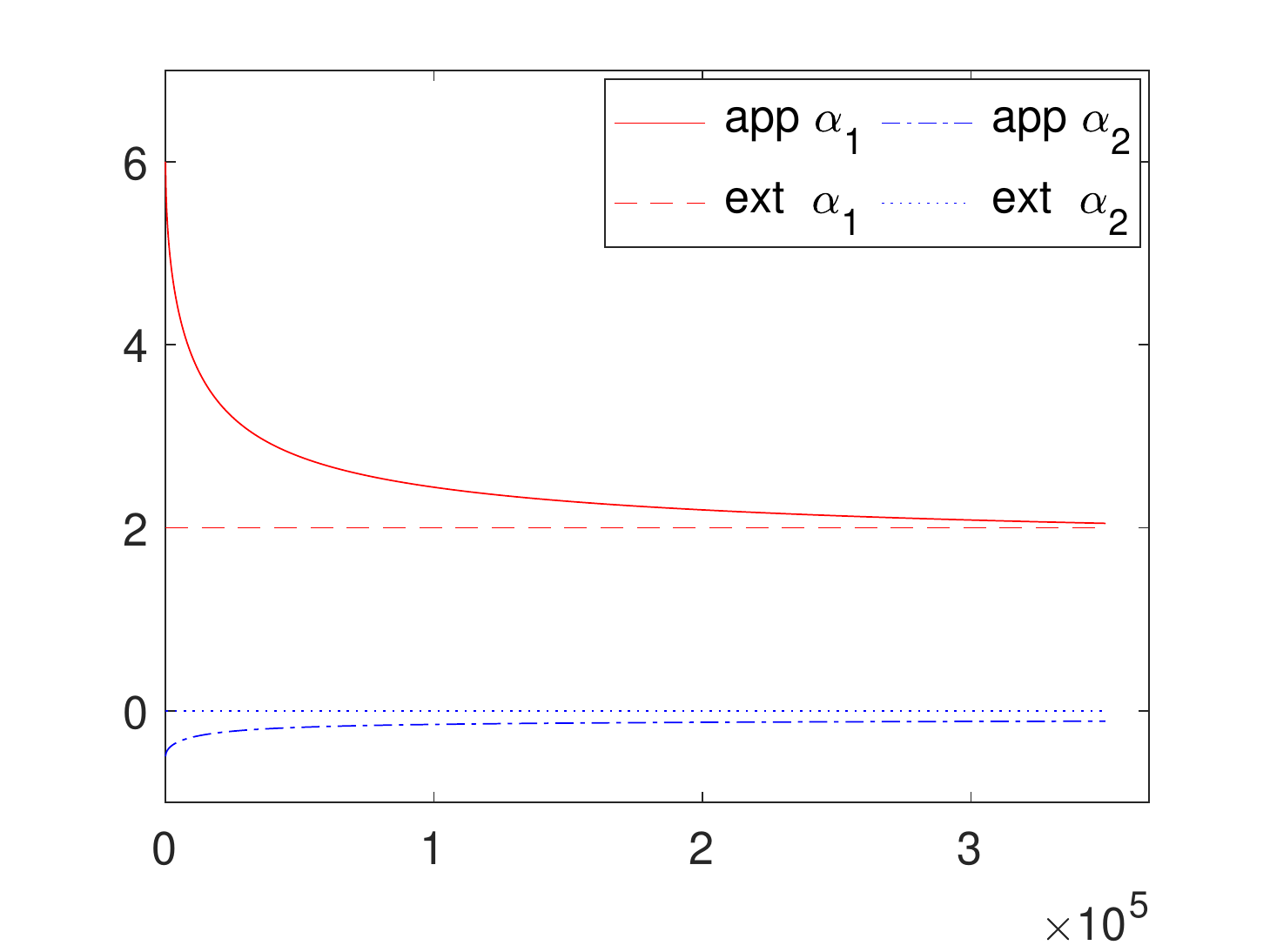}
\includegraphics[width=5.5cm]{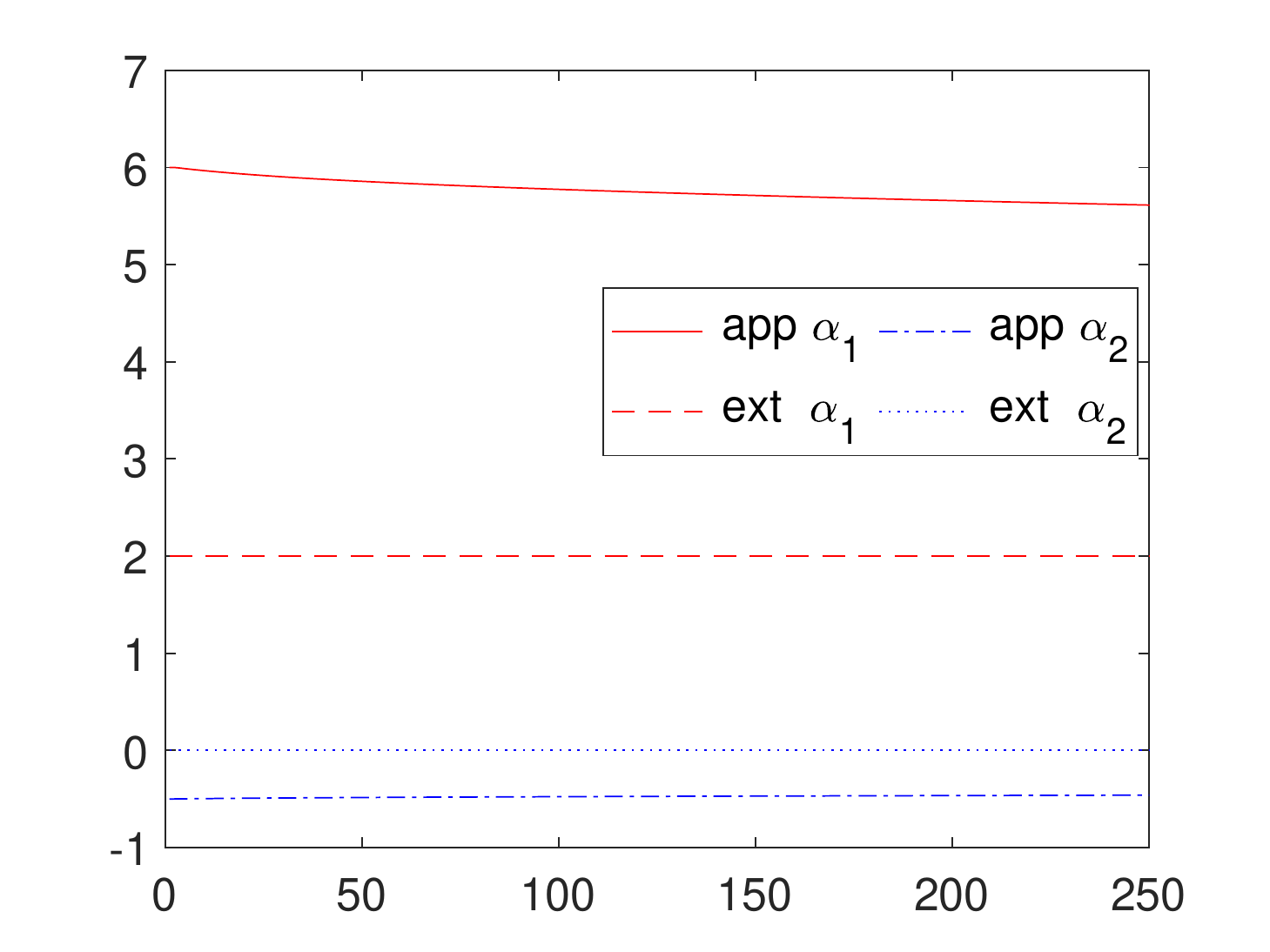}
\caption{Test 2, all-at-once setting: reconstructed parameter over iteration index (left) and zoom of first 250 iterations (right).}
\label{test-res-aao-alp}
\bigskip

\includegraphics[width=5.5cm]{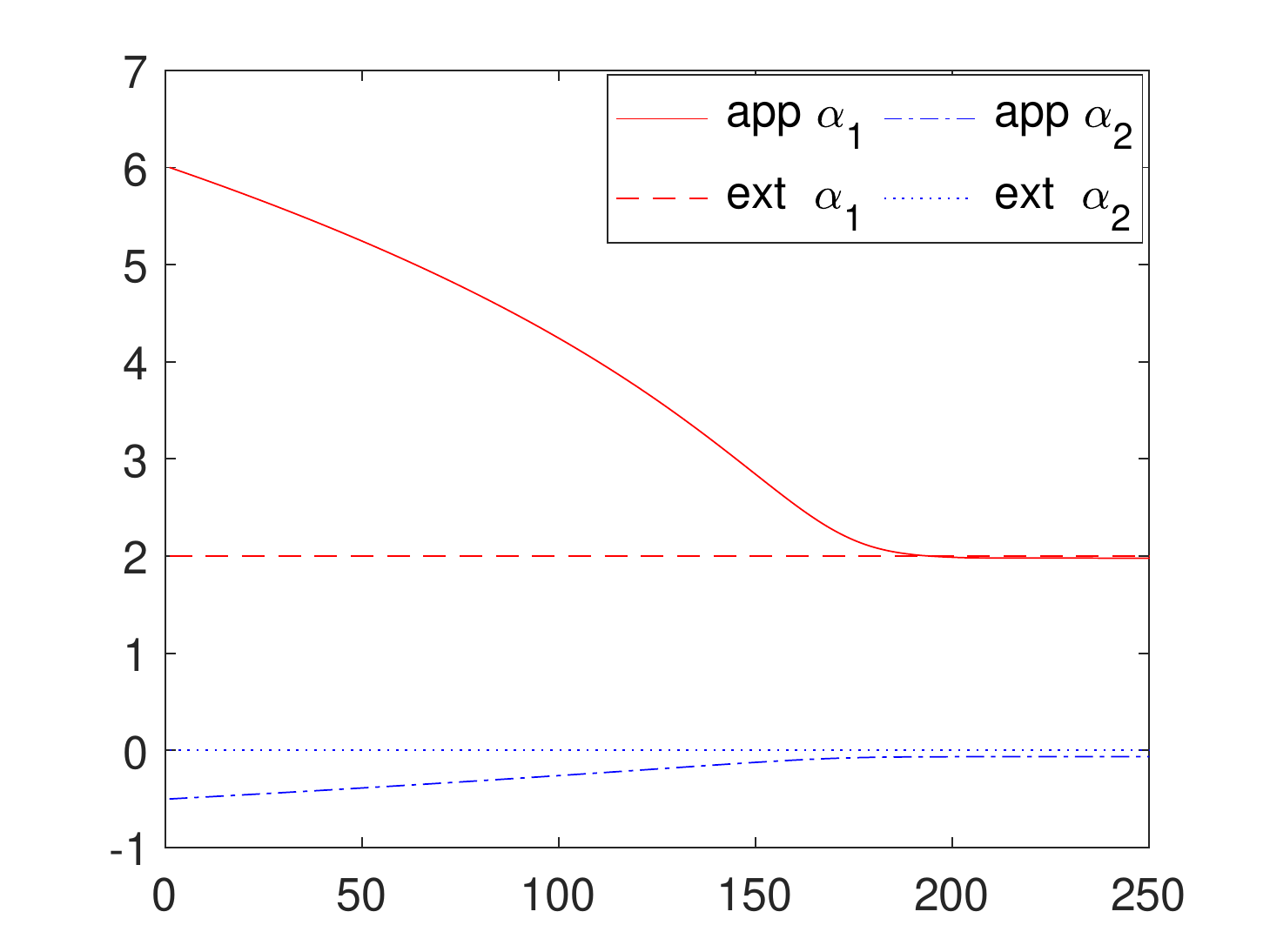}
\includegraphics[width=5.5cm]{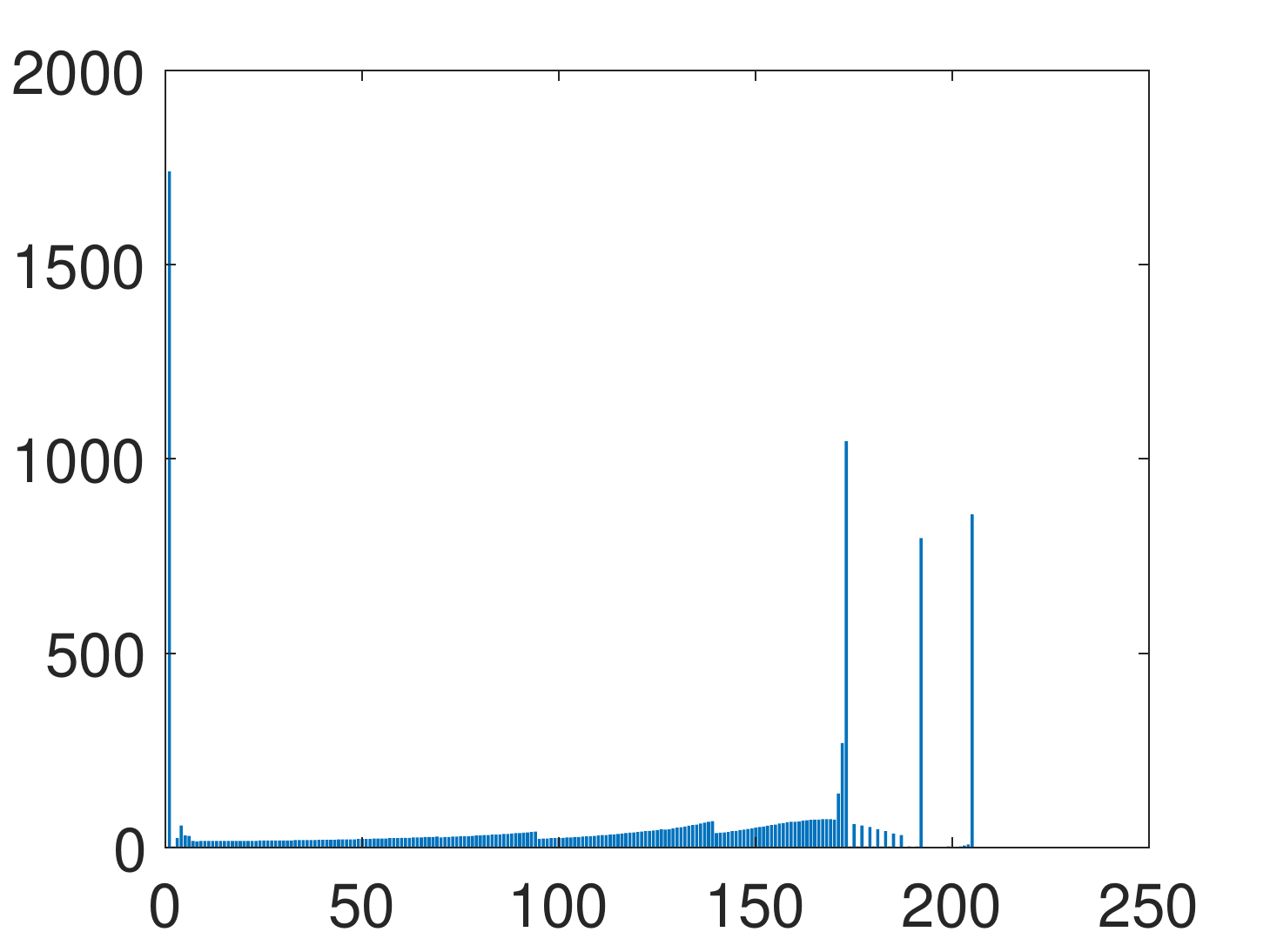}
\caption{Test 2, reduced setting: reconstructed parameter (left) and number of internal loops (right) in each Landweber iteration.}
\label{test-res-red-alp}
\end{figure}

Figures \ref{test-res-aao-m} and \ref{test-res-red-m} respectively present the reconstructed state in the all-at-once setting and in the reduced setting. The all-at-once Landweber ran with $350,000$ iterations and Landweber step size $\mu=1$, while the reduced one ran with $250$ iterations of step size $\mu=1$. The reconstructed parameters $\altil_1, \altil_2$ are depicted in Figures \ref{test-res-aao-alp}, \ref{test-res-red-alp} (left), where both settings confirm the results at acceptable error level.

In the all-at-once setting, the Landweber approximation is applied to both state $\m$ and parameter $\altil$. In the reduced setting, only the parameter is approximated, while the state is calculated exactly with the help of the parameter-to-state map. Hence, the reduced version requires more Landweber iterations to reach an acceptable tolerance compared to the reduced one. Despite the lower number of Landweber iterations, the reduced setting, on the other hand, executes another amount of internal loops each time it calls the LLG solver. Figure \ref{test-res-red-alp} (right) shows the number of internal loops in 250 reduced Landweber iterations. The total number of internal loops is 11095.

With the same step size $\mu=1$, the runtime reports: 
79,000 seconds for 350,000 all-at-once iterations and  90,000 seconds for 250 reduced iterations. The reduced Landweber can be sped up by using a larger step size $\mu=10$, which is feasible in this setting; however, it is not feasible in the all-at-once one. Indeed, the step sizes in the two settings are not correlated since they are chosen subject to the criterion $\mu\in\big(0,\frac{1}{\|F'(x)\|^2}\big]$ for all $x\in\cB_\rho(x^0)$, and according to the problem formulation, the forward operators, thus their derivatives, are different in each setting.

\begin{figure}[!htb] 
\centering
\includegraphics[width=4.3cm]{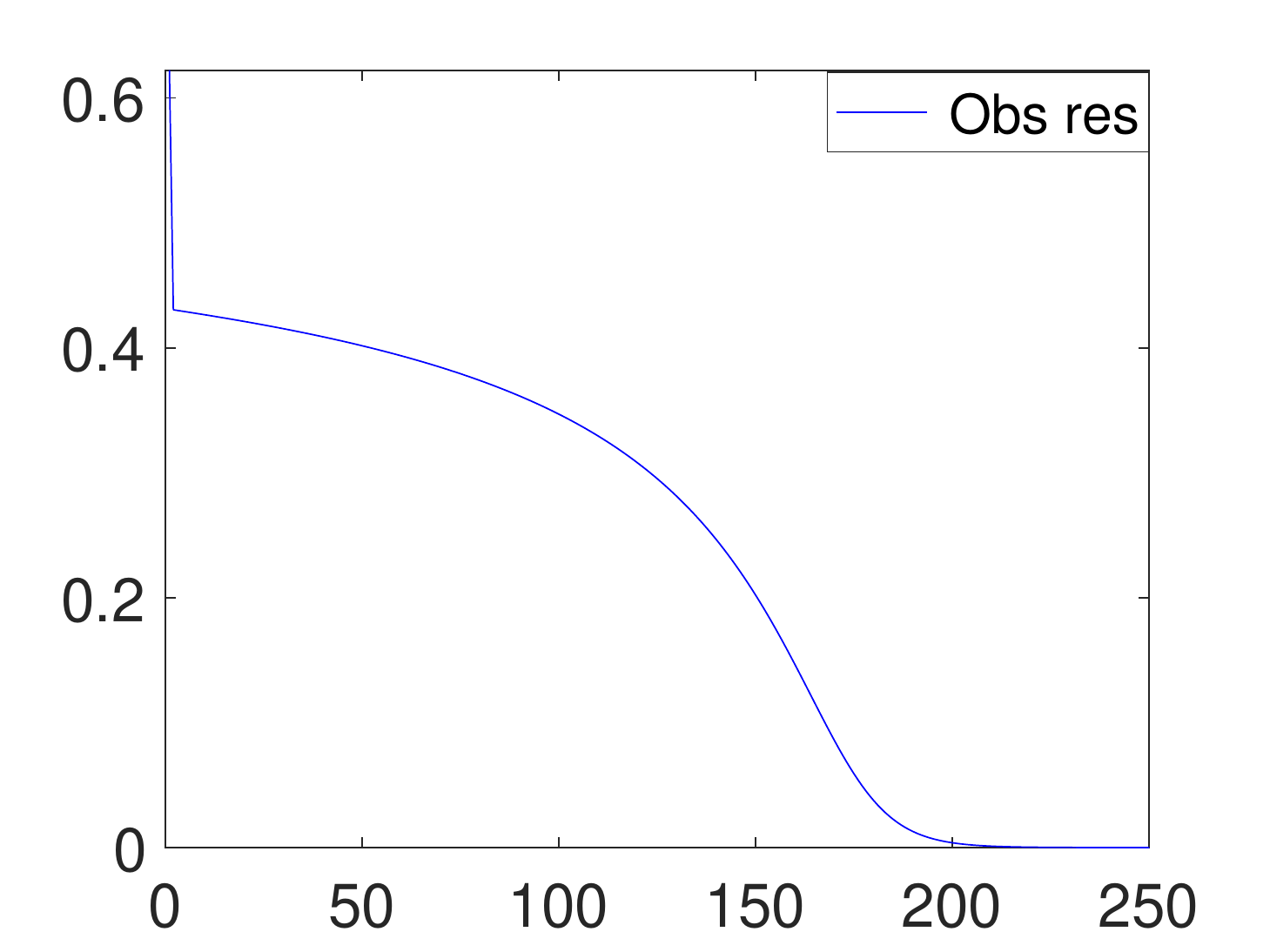}
\hspace{-0.5cm}
\includegraphics[width=4.3cm]{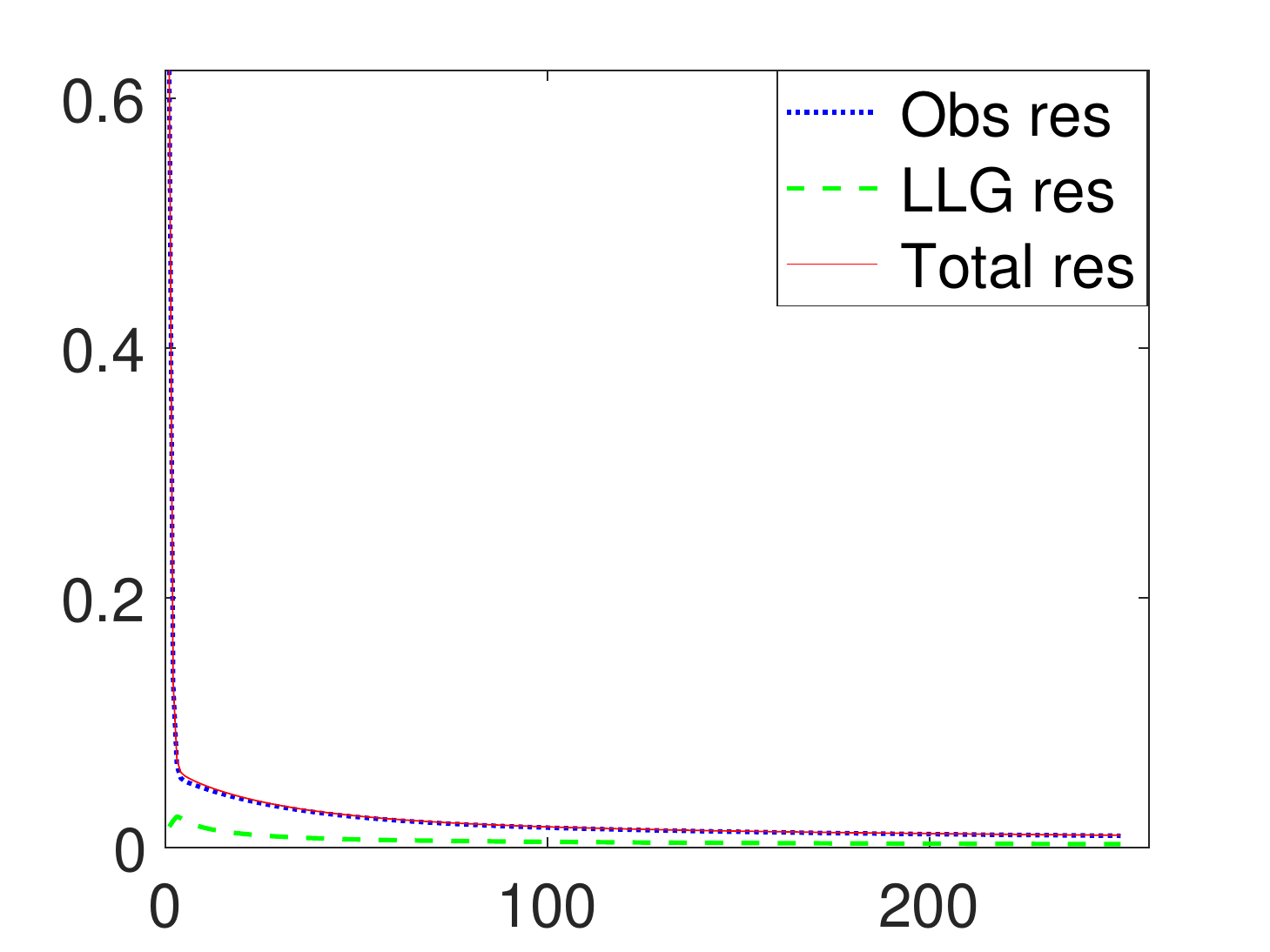}
\hspace{-0.5cm}
\includegraphics[width=4.3cm]{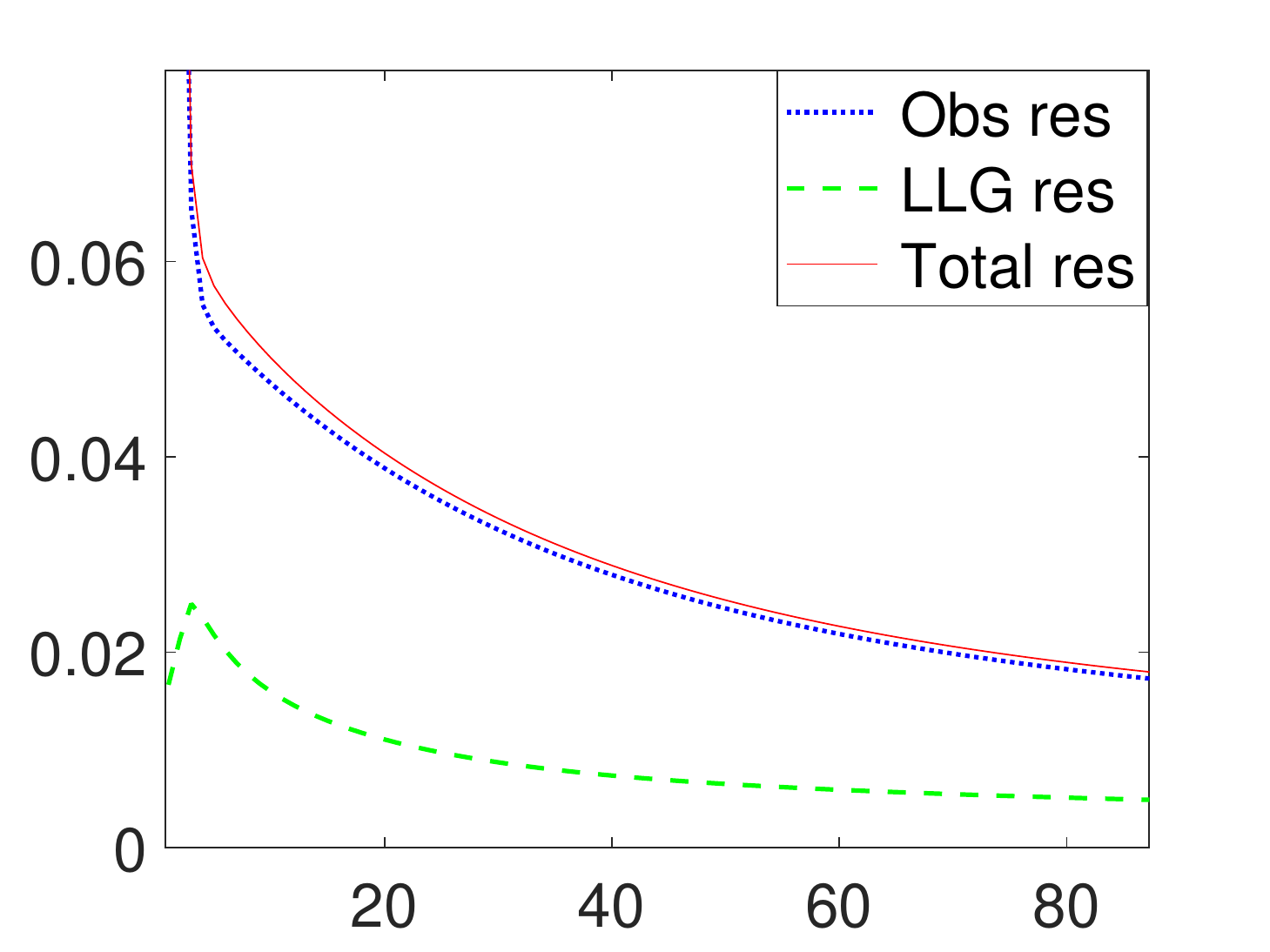}
\caption{Test 2, residual over iteration index: reduced setting (left), first 250 iterations for the all-at-once setting (middle), and a zoom of the all-at-once residual plot (right).}
\label{test-res-redaao-error}
\end{figure}

In Figure \ref{test-res-redaao-error}, the left plot displays the observation residual $\|\cK\m_t-y\|_{L^2(0,T)}$ of each iteration in the reduced version.
The middle and right plots together display the observation residual, the LLG residual and the $\ell^2$-norm total residual in the all-at-once version. The LLG residual is measured in the $H^1(0,T,H^1(\Om2R))^*$-norm.

\begin{figure}[!htb] 
\vspace{1.5cm}
\centering
\includegraphics[width=4.3cm]{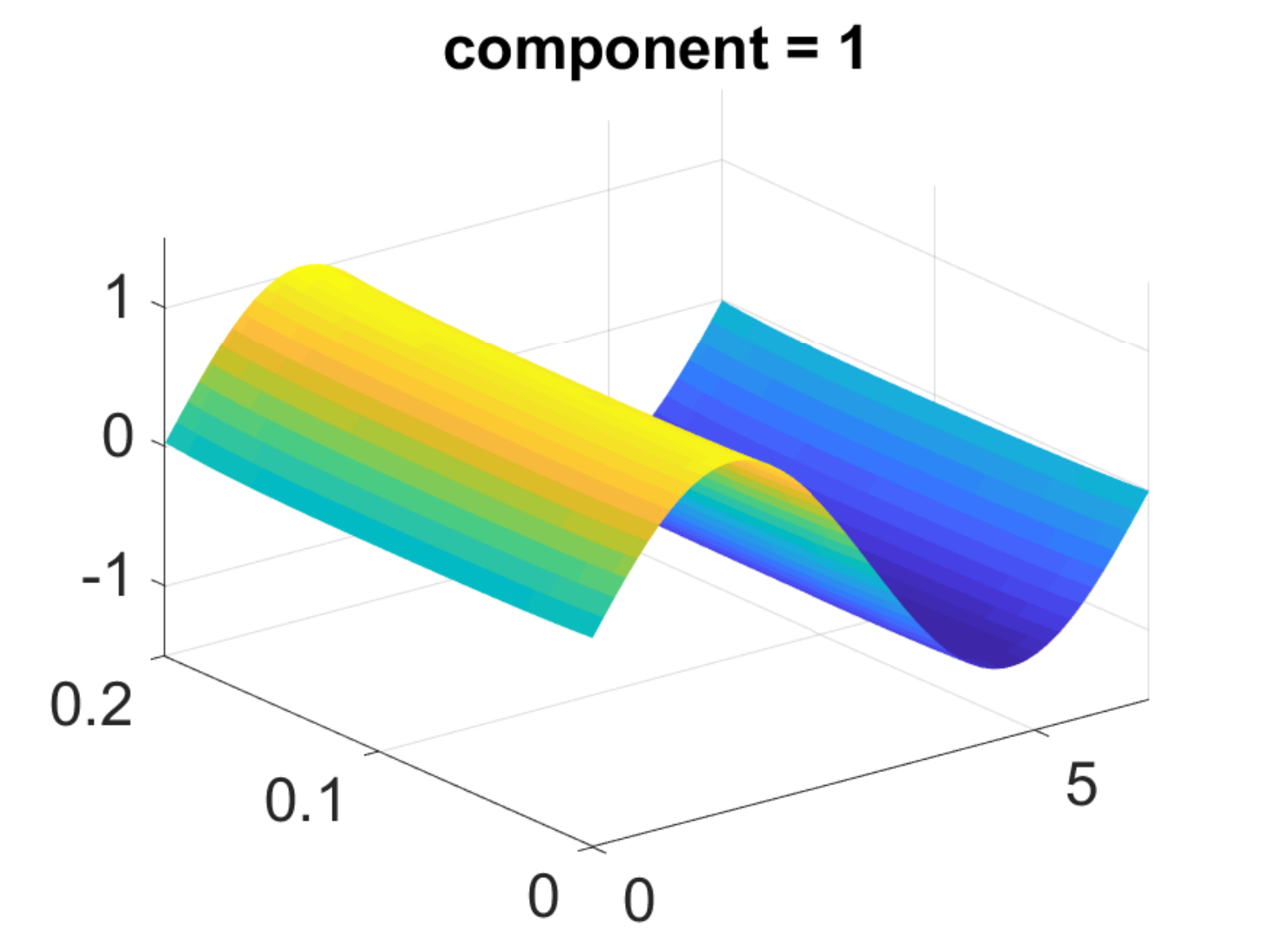}
\hspace{-0.5cm}
\includegraphics[width=4.3cm]{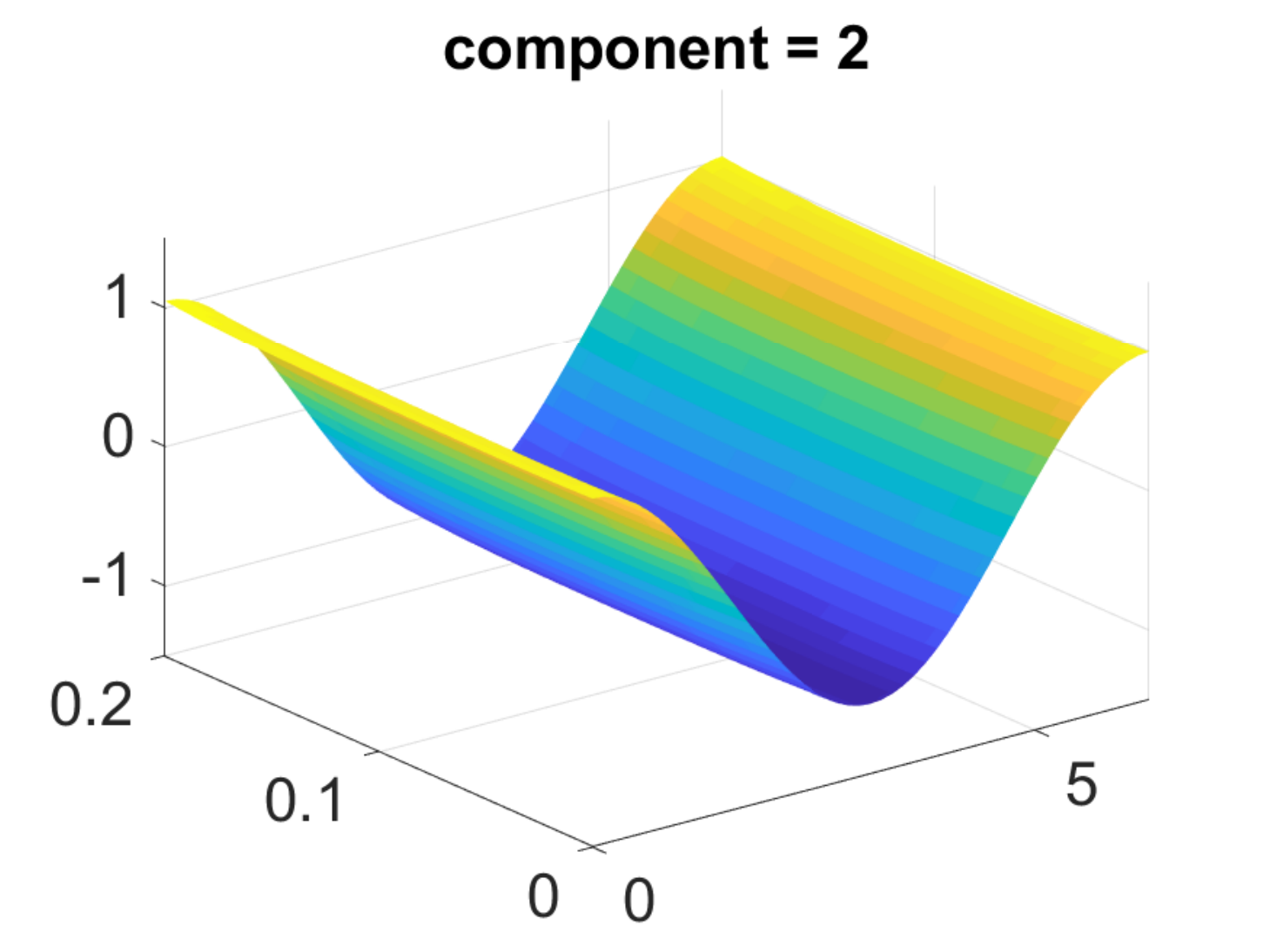}
\hspace{-0.5cm}
\includegraphics[width=4.3cm]{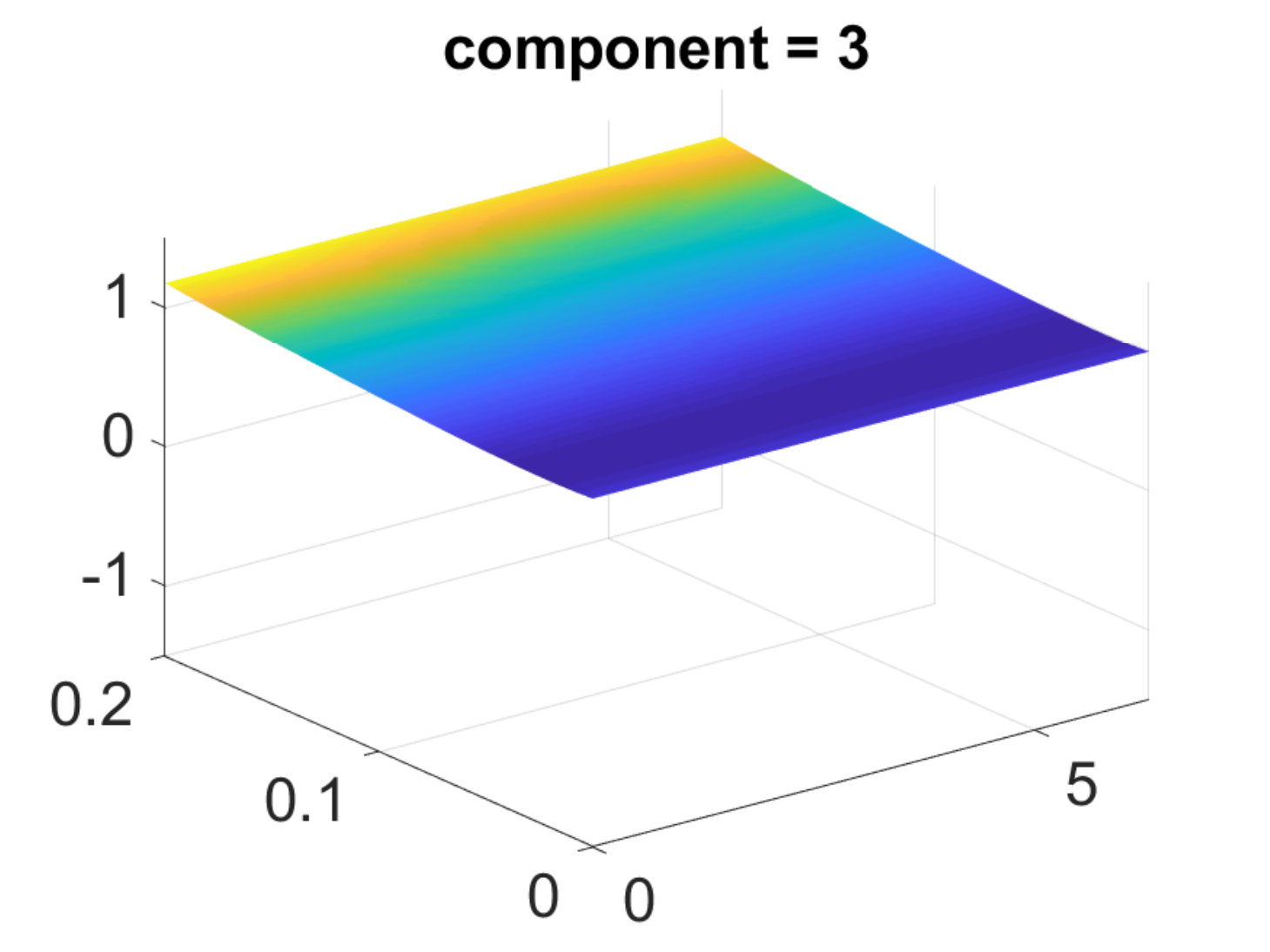}\\[1ex]
\includegraphics[width=4.3cm]{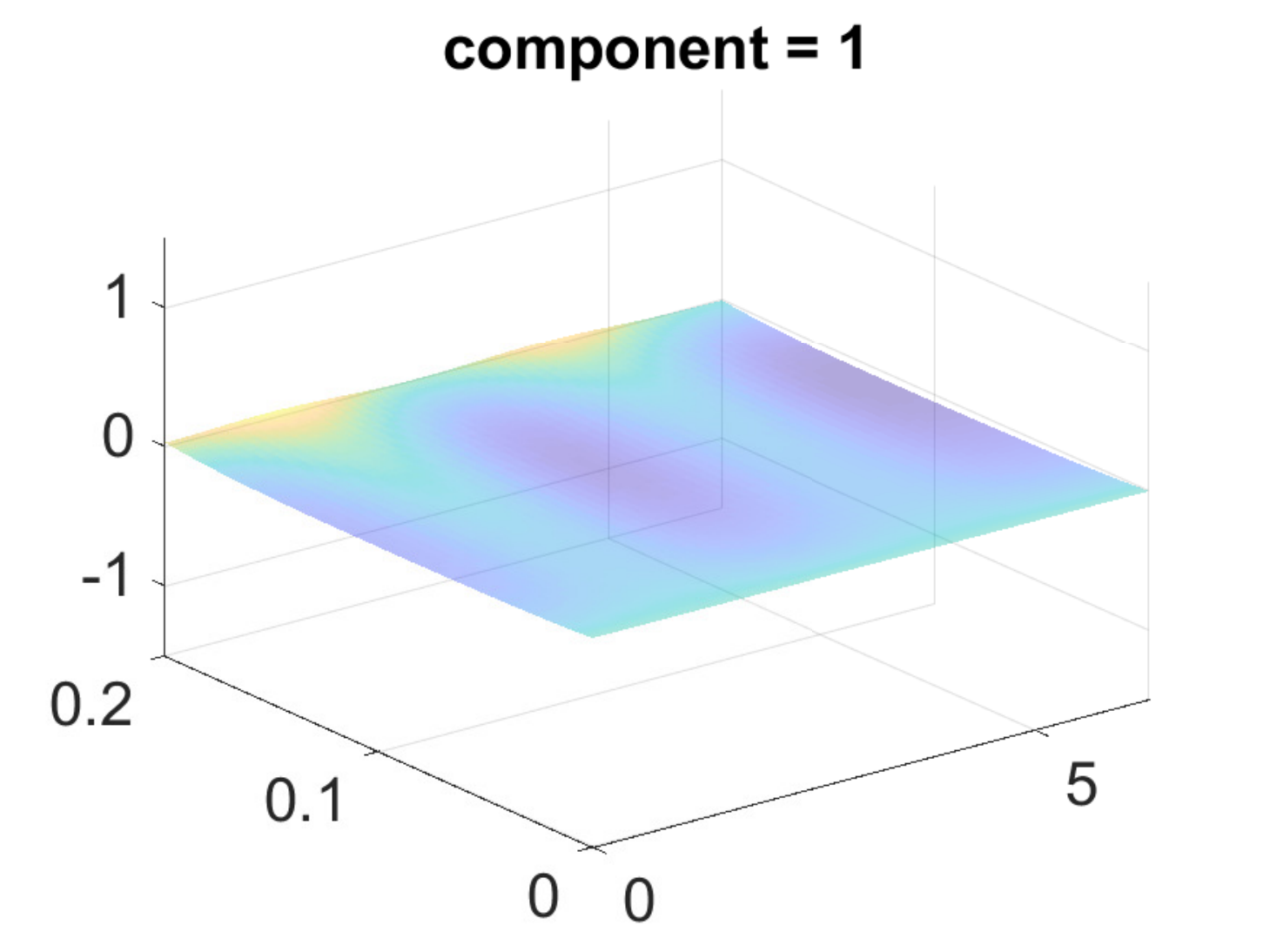}
\hspace{-0.5cm}
\includegraphics[width=4.3cm]{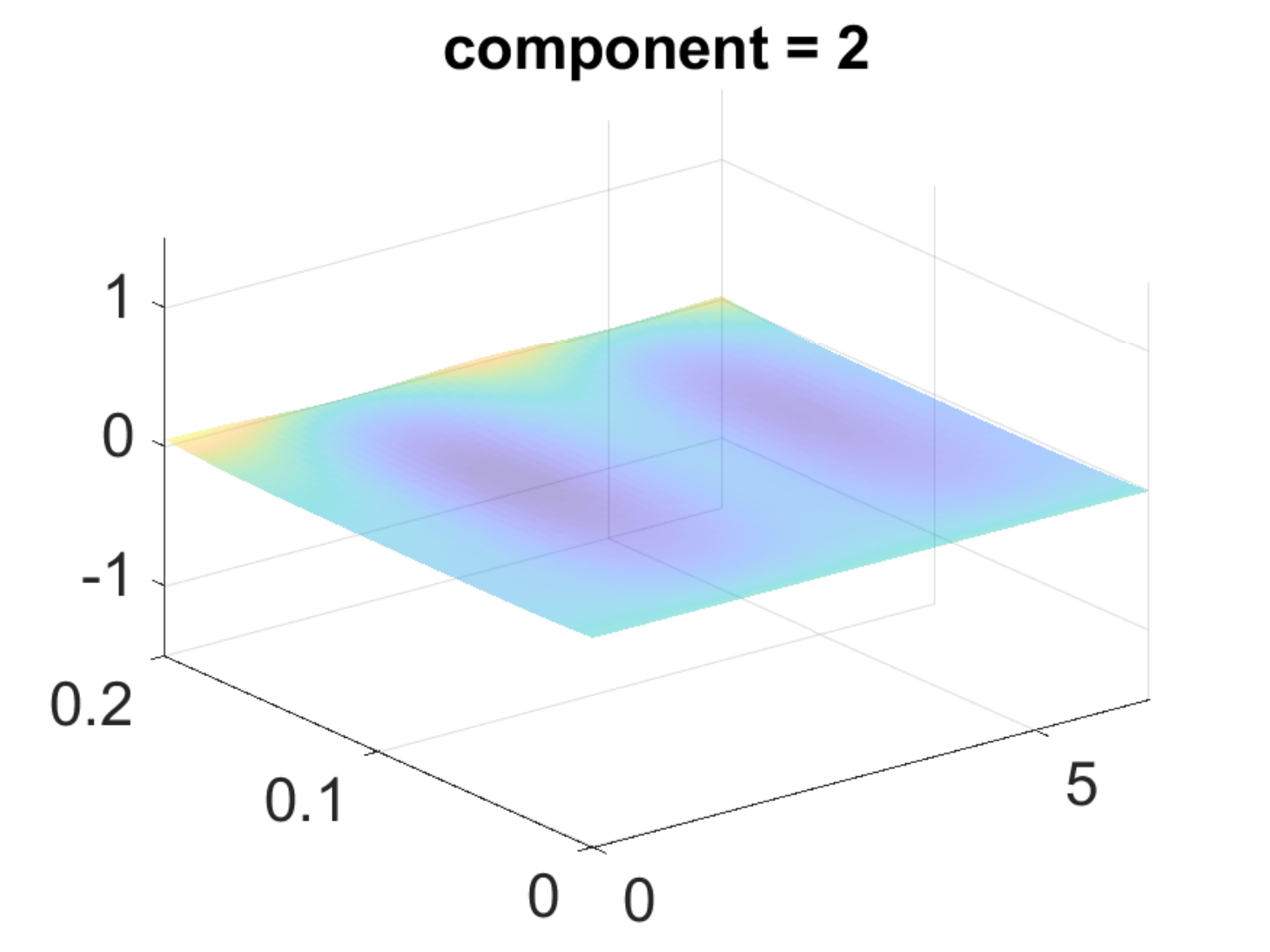}
\hspace{-0.5cm}
\includegraphics[width=4.3cm]{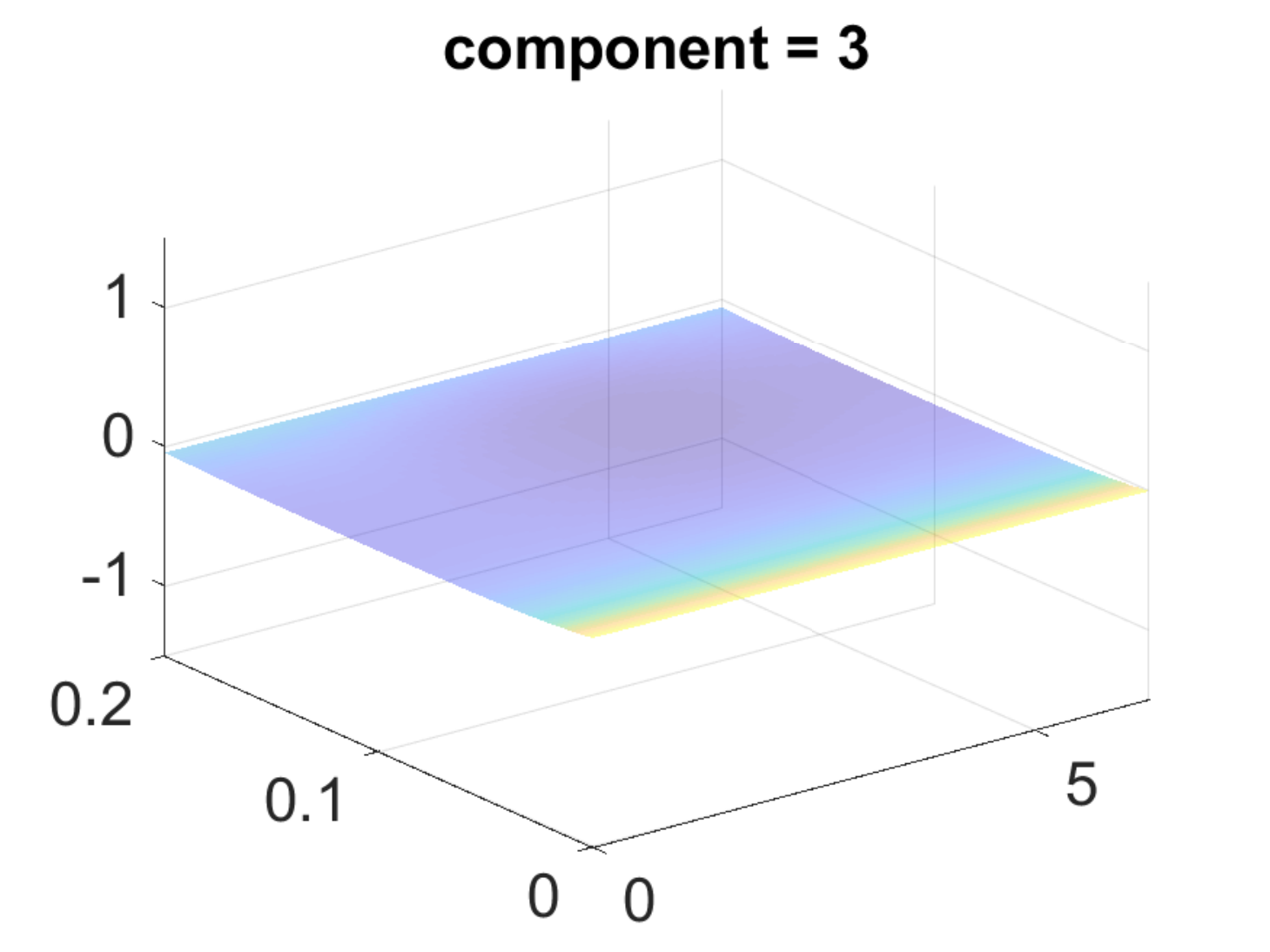}
\caption{Test 3, 3\% noise, all-at-once setting. Reconstructed $\mv$ (top) and $\mv-\mv_{\text{exact}}$ (bottom) plotted against space (x-axis) and time (y-axis).}
\label{test-noise03-aao-m}
\end{figure}

\begin{figure}[!htb] 
\centering
\includegraphics[width=4.3cm]{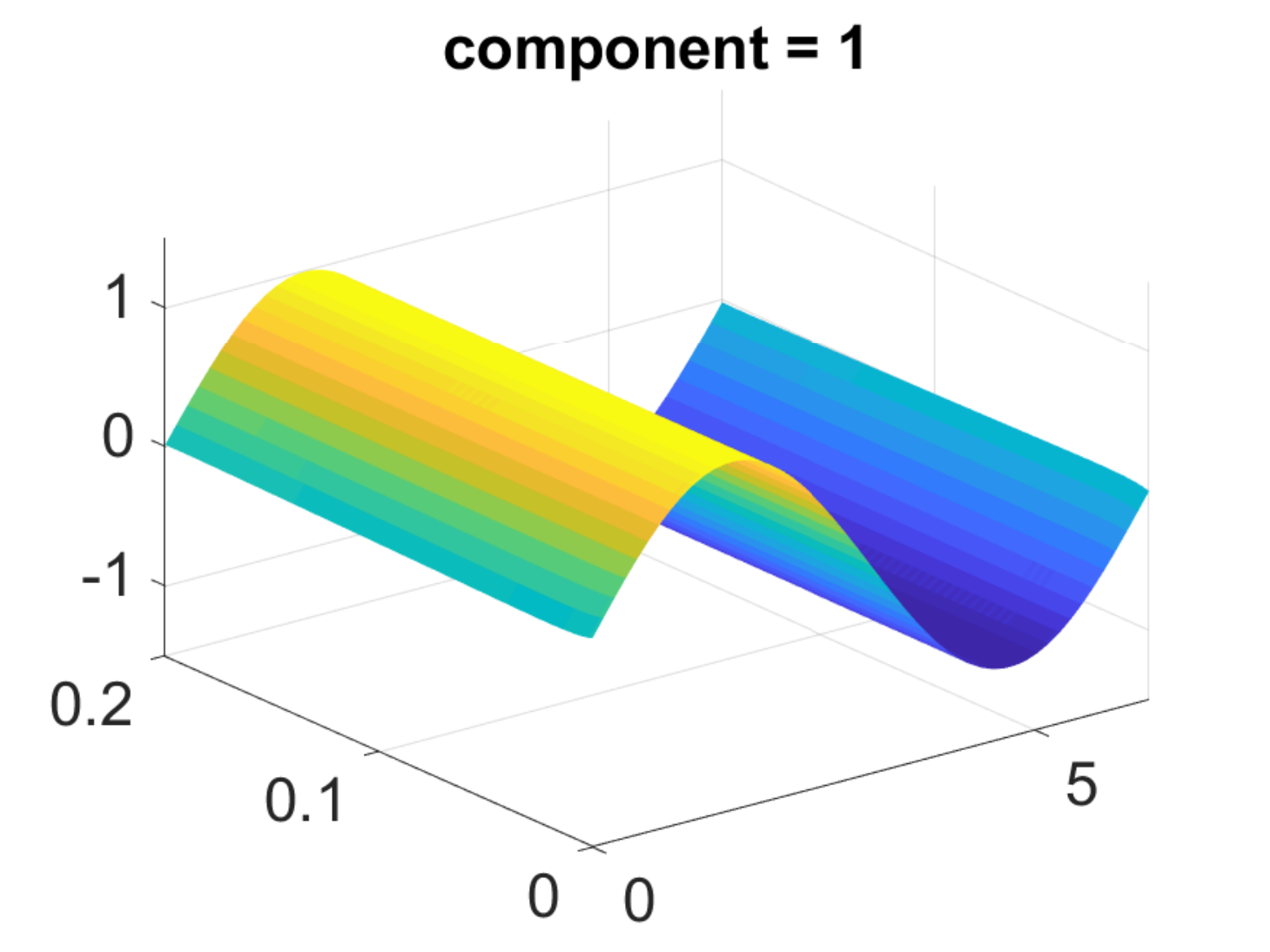}
\hspace{-0.5cm}
\includegraphics[width=4.3cm]{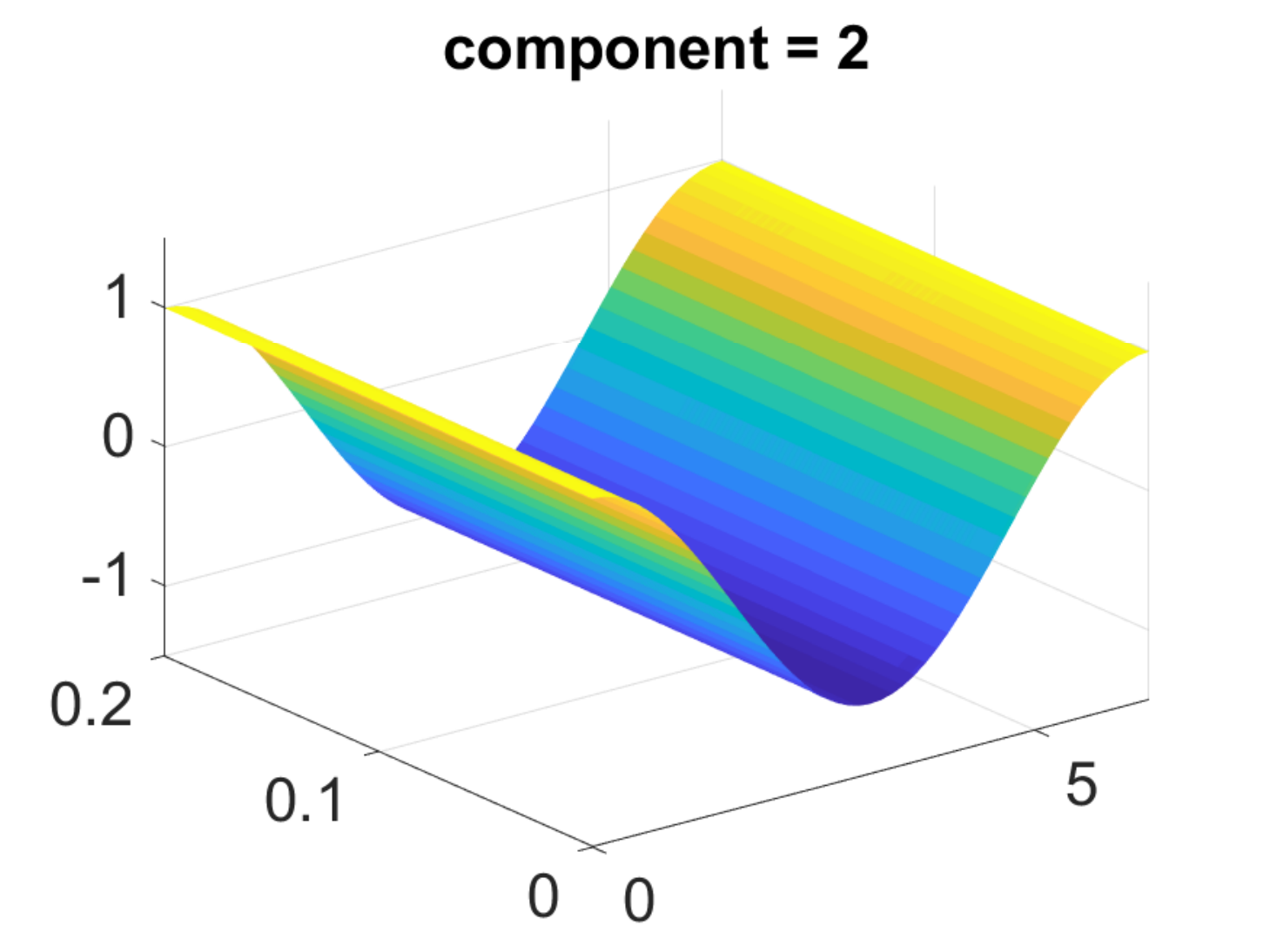}
\hspace{-0.5cm}
\includegraphics[width=4.3cm]{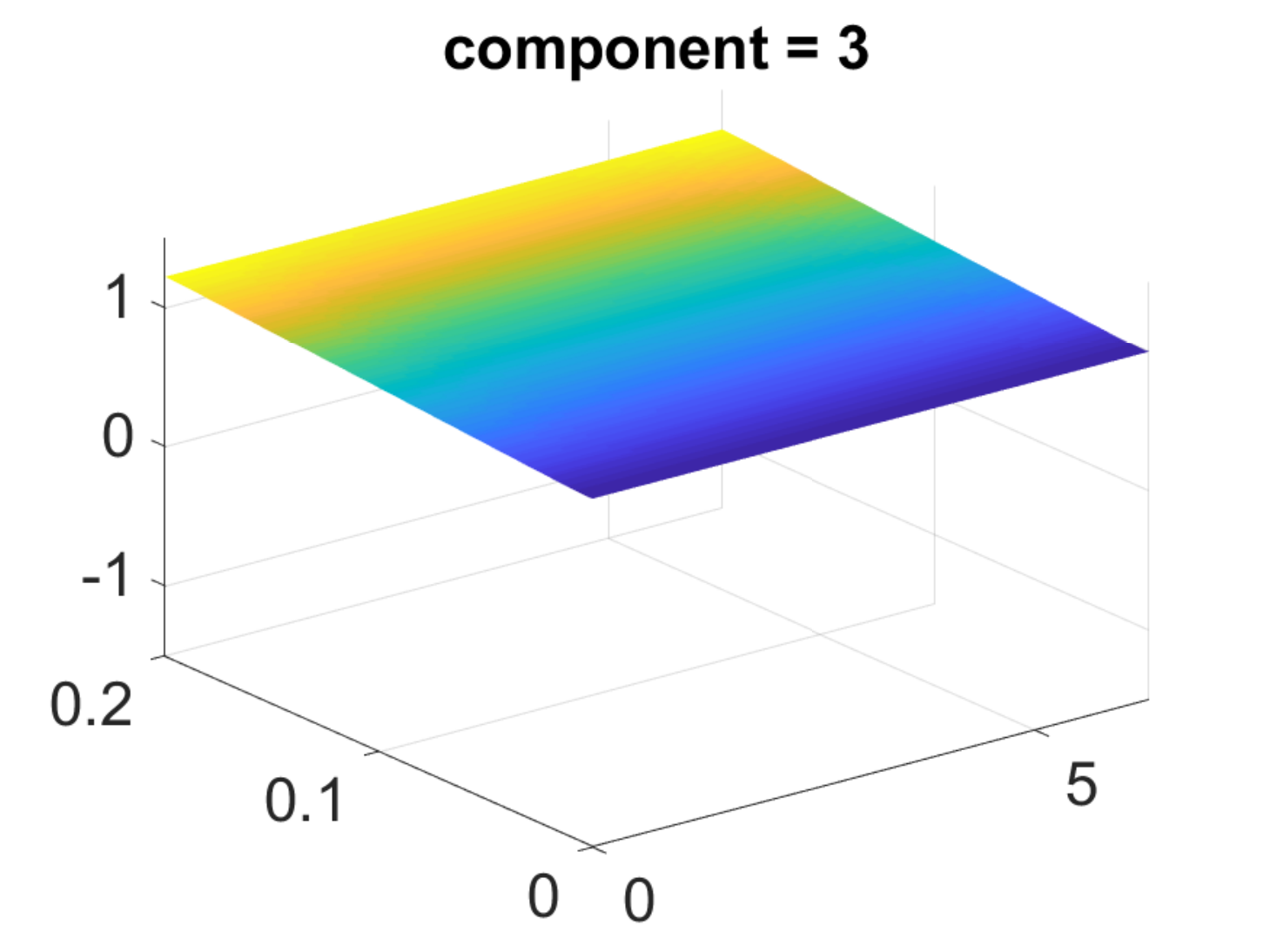}\\[1ex]
\includegraphics[width=4.3cm]{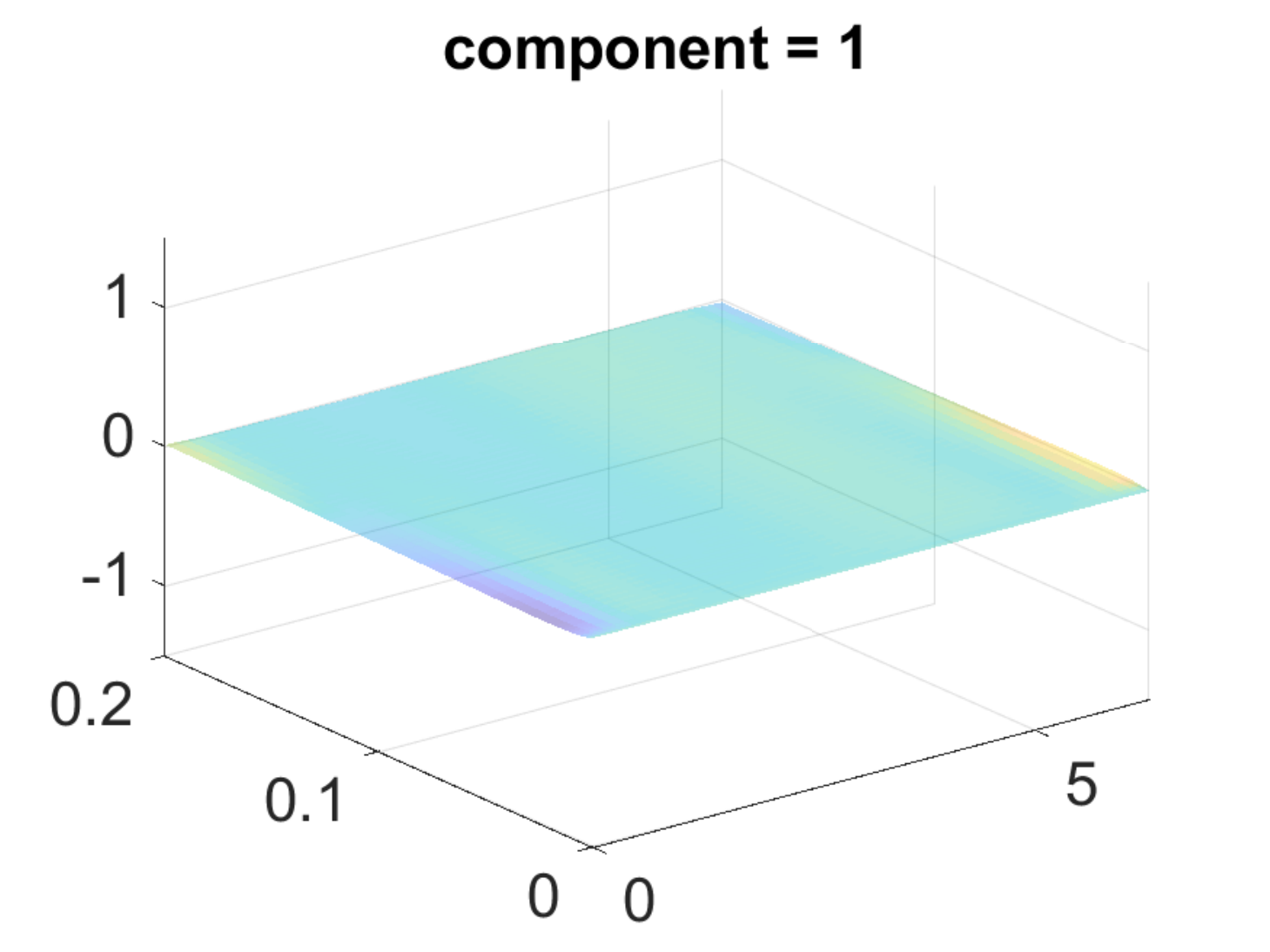}
\hspace{-0.5cm}
\includegraphics[width=4.3cm]{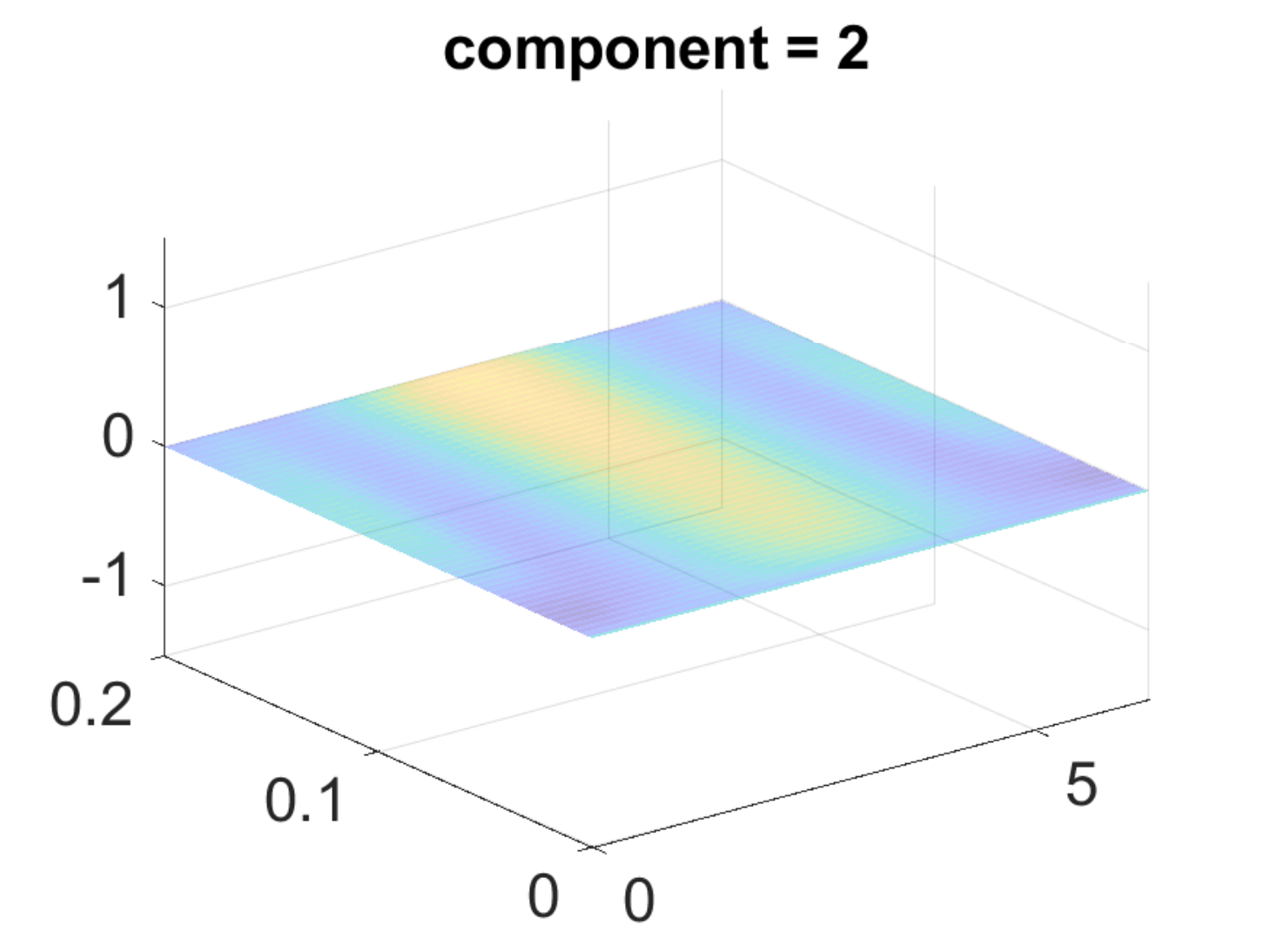}
\hspace{-0.5cm}
\includegraphics[width=4.3cm]{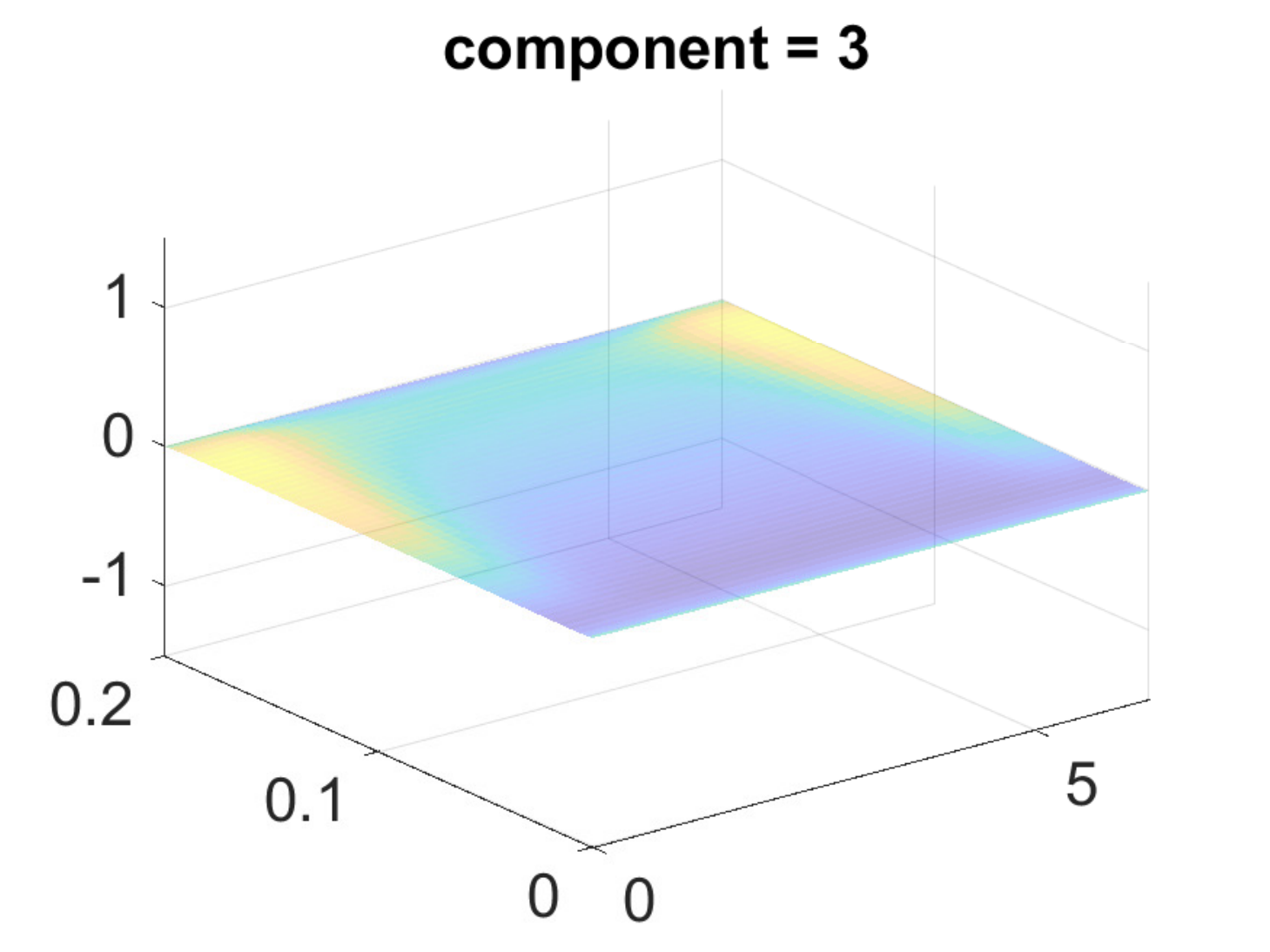}
\caption{Test 3, 3\% noise, reduced setting. Reconstructed $\mv$ (top) and $\mv-\mv_{\text{exact}}$ (bottom) plotted against space (x-axis) and time (y-axis).}
\label{test-noise03-red-m}
\end{figure}

\begin{figure}[!htb] 
\vspace{1cm}
\centering
\includegraphics[width=6cm,height=5cm]{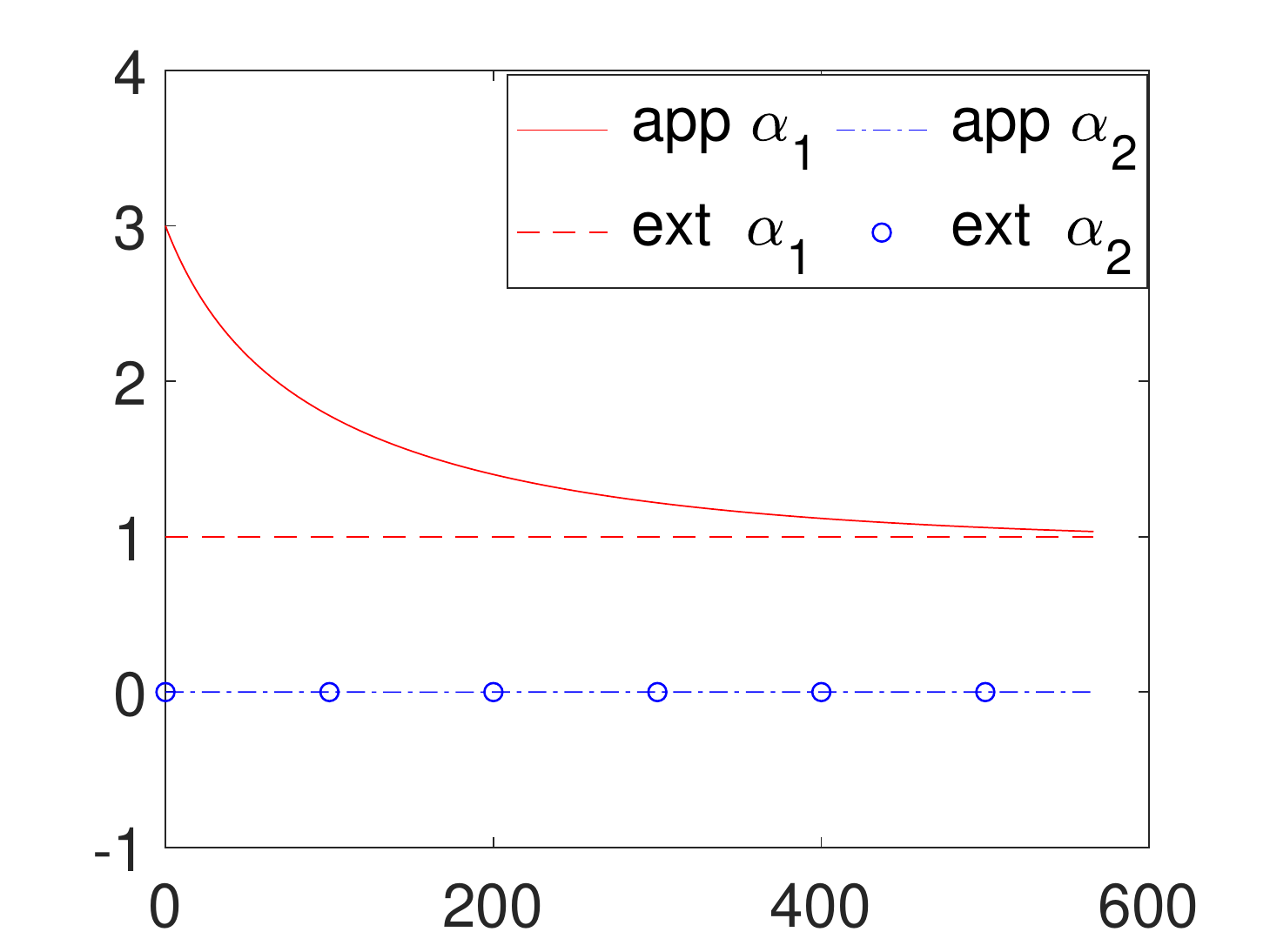}
\includegraphics[width=6cm,height=5cm]{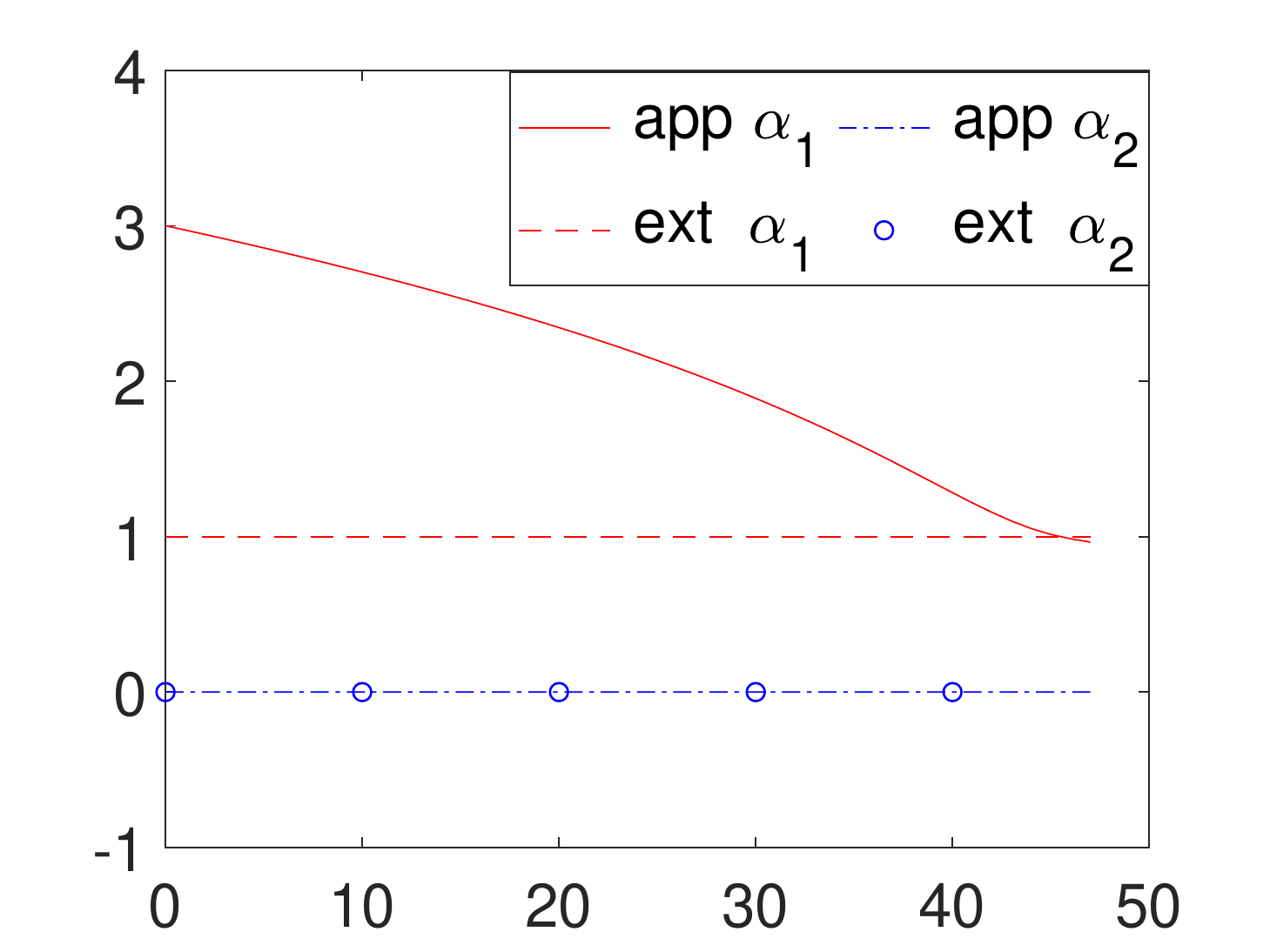}
\caption{Test 3, 3\% noise, reconstructed parameter over iteration index. Left: all-at-once setting. Right: reduced setting.}
\label{test-noise03-alp}

\vspace{1cm}

\includegraphics[width=5.5cm]{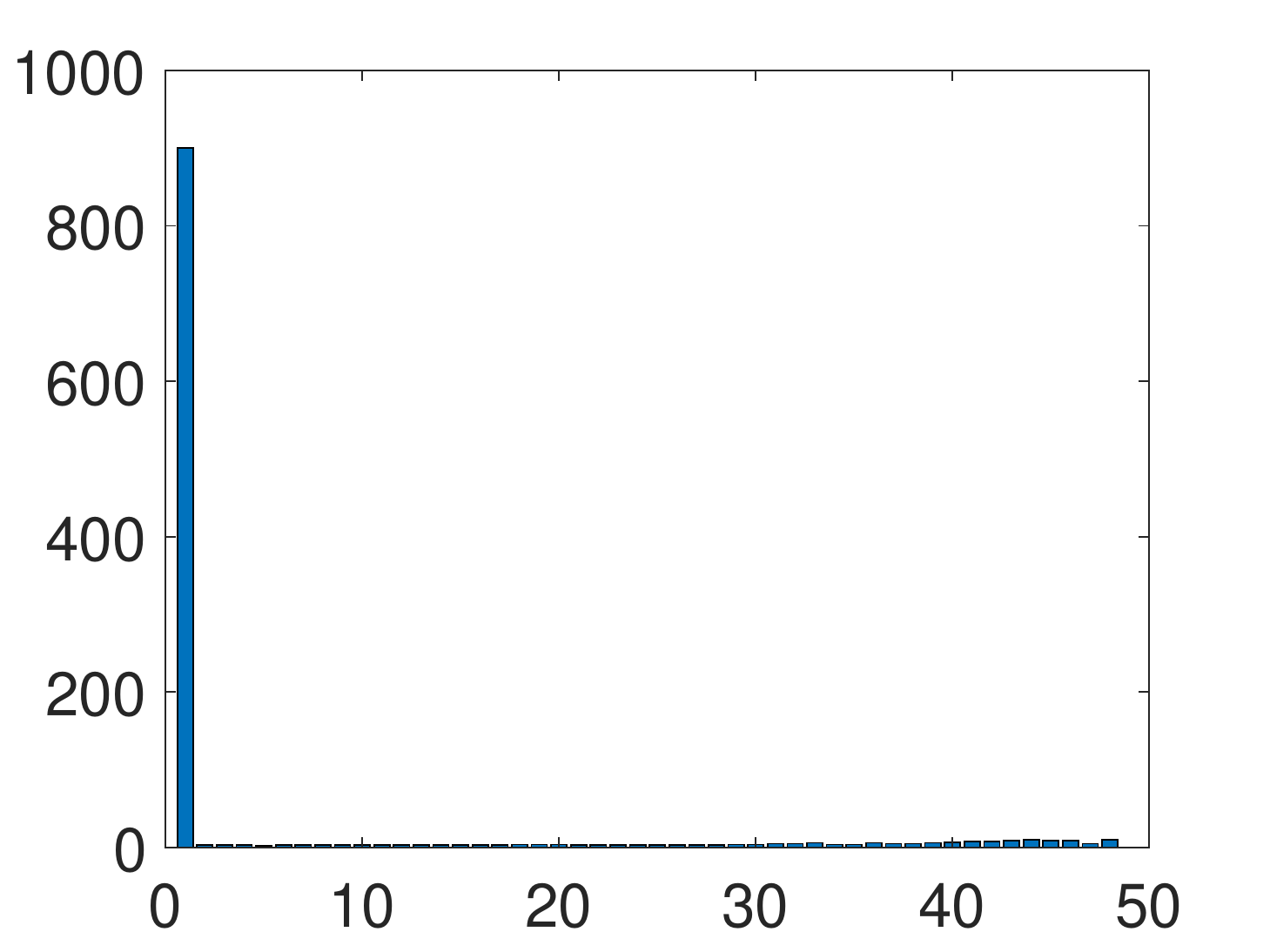}
\caption{Test 3, 3\% noise, reduced setting: number of internal loops in each Landweber iteration.}
\label{test-noise03-iter}

\vspace{1.5cm}

\includegraphics[width=5.5cm]{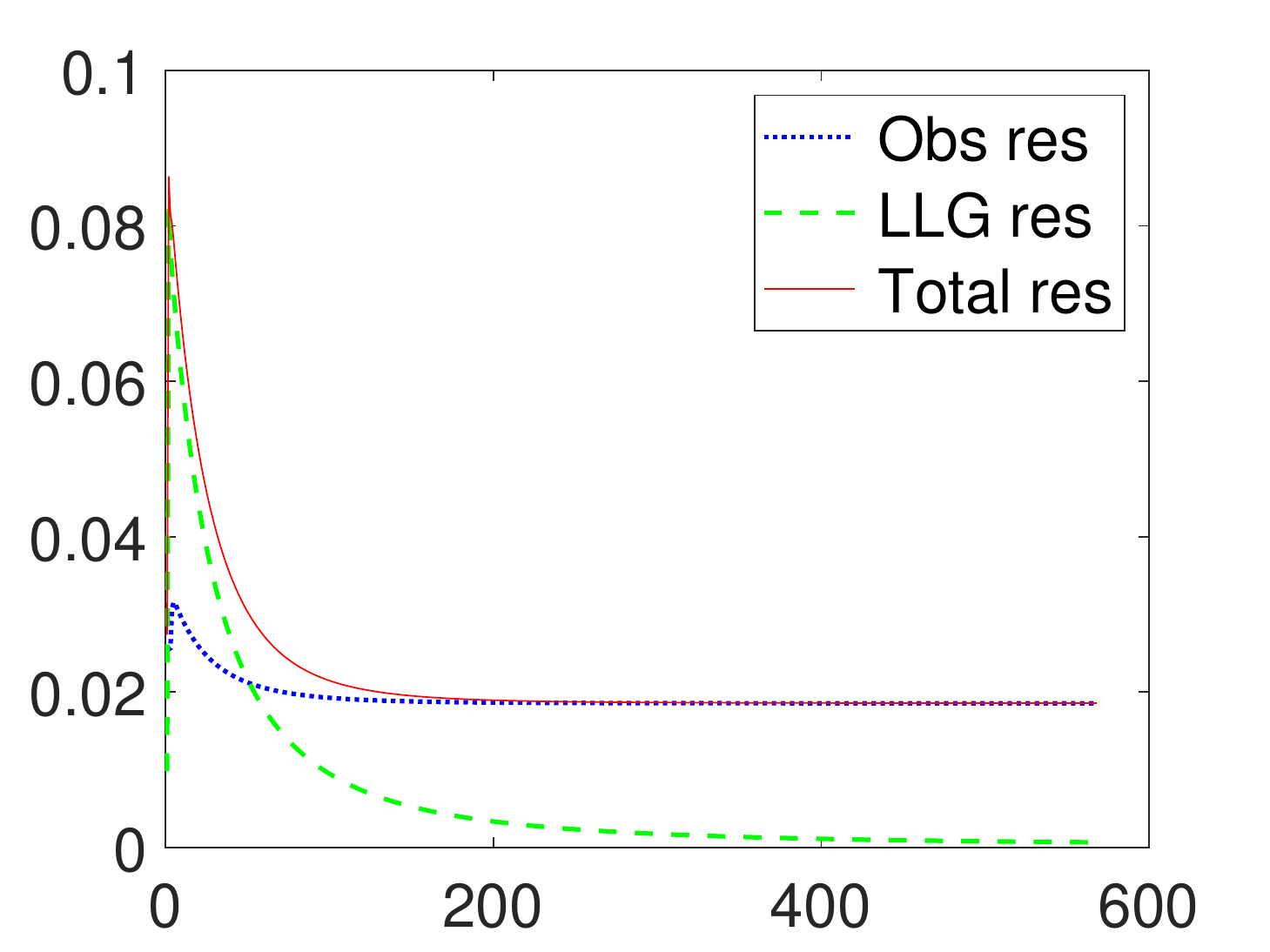}
\includegraphics[width=5.5cm]{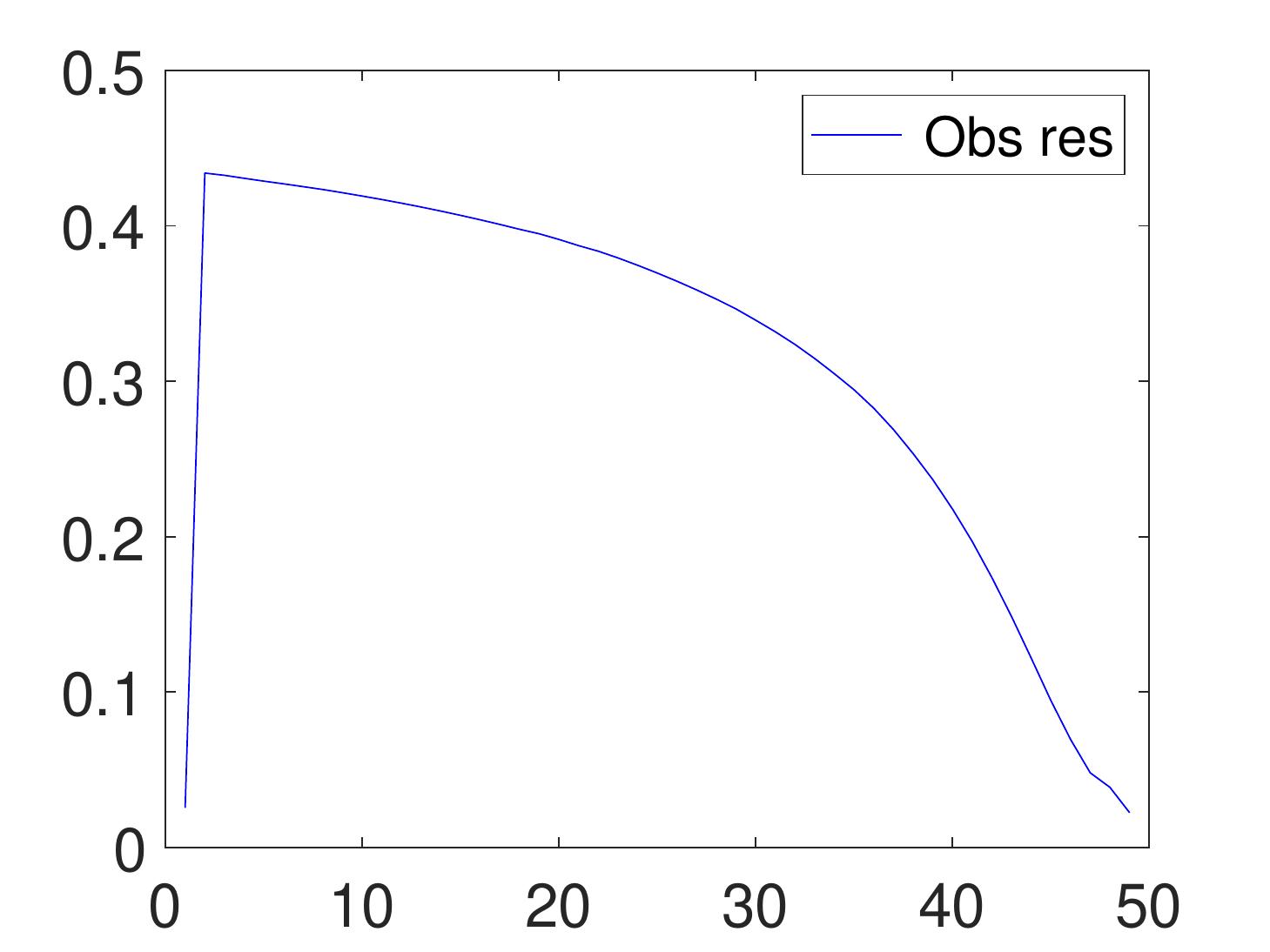}
\caption{Test 3, 3\% noise, residual over iteration index. Left: all-at-once setting. Right: reduced setting.}
\label{test-noise03-error}
\end{figure}

We now examine the reconstruction for Test 3 in Table \ref{tab-1} in case of noisy data. We perturb the exact measured voltage with 3\%, 5\% and 10\% random noise. Also, instead of initializing the algorithm at $\m_0$, we choose a perturbed version of it, namely, $\m_\text{init}=\m_0-0.1\sin(20 t)$. The iteration is stopped according to the discrepancy principle. Figures \ref{test-noise03-aao-m}-\ref{test-noise03-error} present the details of the test in the same fashion with Test 2 mentioned above, but with the involvement of $3\%$ data noise.
Table \ref{noise} respectively reports, for each noise level 3\%, 5\% and 10\%, the number of iterations (\#it), the LLG residual ($r_{llg}$), the observation residual ($r_{obs}$) and the reconstruction error $e_{\alpha} := |\alpha-\alpha_{ext}|$ in both settings. The LLG residuals in the reduced setting are typically smaller than the ones in the all-at-once setting, since there the states are solved exactly.

\newpage

\vspace{1cm}

\clearpage

\vspace{1cm}

\clearpage

\newpage

\vspace{1cm}

\clearpage

\clearpage


\begin{table}
\renewcommand{\arraystretch}{1.5}
\caption{Reconstruction with noisy data.}
\bigskip
\centering
\begin{tabular}{r|l|l}
\hline
& \qquad\qquad\qquad All-at-once &\qquad\qquad\qquad Reduced\\
 \hline
$\delta$ & 
\#it\quad $r_{llg}$ \qquad $r_{obs}$ \qquad $e_{\alpha_1}$ \quad $e_{\alpha_2}$ \quad 
&\#it  \quad$r_{llg}$ \qquad $r_{obs}$ \qquad $e_{\alpha_1}$ \quad $e_{\alpha_2}$ \\[1ex]
\hline
10\% &
259\,\,	 	0.0022 		\,\, 0.0619  \,\, 0.292   \,\, 0.090&
49\,\, $3\times10^{-6}$\,\, 0.0703  \,\, 0.030   \,\,0.034 
\\
5\% &
401\,\, 		0.0011      \,\, 0.0309 \,\, 0.125 \,\,  0.072&
49\,\, $3\times10^{-6}$\,\, 0.0321 \,\, 0.040 \,\,0.033
\\
3\% &
564\,\,    	0.0007      \,\, 0.0186  \,\, 0.040   \,\, 0.062&
49\,\, $3\times10^{-6}$	\,\,0.0200  \,\, 0.044   \,\,0.033\\
\hline
\end{tabular}
\label{noise}
\end{table}

\section*{Conclusion and outlook}
In \cite{KNSW}, it has been discussed that the modelling of the response of the magnetization vector field in response to an externally applied dynamic magnetic field is essential to enable a model-based calibration in magnetic particle imaging. The relation between magnetization $\m$ and external field $\h$ is nonlinear and can be described by the Landau-Lifshitz-Gilbert equation, which particularly models the relaxation effect, i.e., the effect that the magnetization vector aligns with the external field only with a certain delay. 

\vspace*{1ex}

This paper investigates the numerical approximation of the magnetization vector $\m$ as well as the physical parameters $\alpha_1, \alpha_2$ by means of the LLG equation as the underlying model and the signals measured by the scanner. The numerical results show that our method is robust. Moreover, the computational schemes are preferable in practice since all the steps involve solving only linear PDEs. In addition, we provide multiple choices to the users: an all-at-once version and a reduced version. 
%

\vspace*{1ex}

\blue{By integrating more sophisticated solvers for the occurring linear partial differential equations, the presented methods shall be extended to efficiently evaluate problems in two and three dimensions.} In addition, they shall be adapted and applied to real measured data. Here, we want to emphasize that this involves not only noise in the measured data, but also inexactness in the model itself: the coil sensitivities can be measured, hence subject to noise, or described analytically and thus subject to model inexactness in comparison to the coil sensitivities of a real receive coil. Furthermore, there are differences between the analytic description of the external magnetic field and the actual applied field. As a consequence, there is no ground truth available to evaluate the proposed methods, which emphasizes the importance of the numerical experiments presented in this work. 

The results yield a good starting point for a data-driven model-based determination of the system function in MPI. To this end, the effective field used in the LLG equation shall be complemented by an anisotropy term, which has been shown to influence the magnetization model as well, see, e.g., \cite{tk18}.
Finally, the use of a system function relying on the LLG model shall be evaluated in the imaging process of MPI.

\section*{Acknowledgments} Tram Nguyen wishes to thank her former supervisor Barbara Kaltenbacher, Alpen-Adria Universit\"at Klagenfurt, for fruitful and inspiring discussions. The work of Anne Wald was funded by the German Federal Ministry of Education and Research (Bundesministerium f\"ur Bildung und Forschung, BMBF) under 05M16TSA. 


\medskip
Received January 2021; revised April 2021.
\medskip

\end{document}